\documentclass[12pt]{article}

\usepackage{epic}
\usepackage{eepic}
\usepackage{epsf} 
\usepackage{ amssymb, amscd}
\usepackage{amsmath}
\usepackage{color}
\numberwithin{equation}{section}
\usepackage{epsfig}
\usepackage{graphicx}

\newcommand{\lori}{\longrightarrow}

\def\cB            {{\mathcal{B}}}

\def\cD            {{\mathcal{D}}}

\def\cF            {{\mathcal{F}}}

\def\cH            {{\mathcal{H}}}

\def\cK            {{\mathcal{K}}}
\def\cL            {{\mathcal{L}}}
\def\cM            {{\mathcal{M}}}
\def\cN            {{\mathcal{N}}}

\def\cP            {{\mathcal{P}}}

\def\cV            {{\mathcal{V}}}
\def\cW            {{\mathcal{W}}}

\def\cZ            {{\mathcal{Z}}}
\def\bbA           {\mathbb{A}}
\def\bbC           {\mathbb{C}}
\def\bbD           {\mathbb{D}}
\def\bbE           {\mathbb{E}}

\def\bbQ           {\mathbb{Q}}
\def\bbR           {\mathbb{R}}
\def\bbT           {\mathbb{T}}
\def\bbZ           {\mathbb{Z}}

\def\MXN           {{}_M {\mathcal{X}}_N}
\def\MXM           {{}_M {\mathcal{X}}_M}

\def\MXMo          {{}_M^{} {\mathcal{X}}_M^0}

\def\NXN           {{}_N {\mathcal{X}}_N}

\def\NXM           {{}_N {\mathcal{X}}_M}

\def\i{{\rm i}}

\title{Modular Invariants\\ and Twisted Equivariant $K$-theory II:\\ Dynkin diagram symmetries}
\author{
{\sc David E.\ Evans}\\
 {\footnotesize School of Mathematics, Cardiff University,}\\
 {\footnotesize Senghennydd Road, Cardiff CF24 4AG, Wales, U.K.}\\
 {\footnotesize e-mail: {\tt EvansDE@cf.ac.uk}}\\ \\
 {\sc Terry  Gannon }\\
 {\footnotesize Department of Mathematics, University of Alberta,}\\
{\footnotesize Edmonton, Alberta, Canada T6G 2G1}\\
{\footnotesize e-mail: {\tt tgannon@math.ualberta.ca}} }

\begin{document}
\maketitle

\begin{abstract} 

The most basic structure of {\it chiral} conformal field theory (CFT) is the Verlinde ring. Freed-Hopkins-Teleman
have expressed the Verlinde ring for the CFT's associated to loop groups, as twisted equivariant K-theory.
We build on their work to express K-theoretically the structures of {\it full} CFT. In particular, the modular invariant partition functions { (which essentially parametrise the
possible full CFTs)} have a rich interpretation within von Neumann algebras (subfactors), which has led to the developments of structures of full CFT such as the full system (fusion ring of defect lines), nimrep (cylindrical partition function), alpha-induction etc. Modular categorical interpretations for these have followed. 
For the generic families of modular invariants
 (i.e. those associated to Dynkin diagram symmetries), we 
 provide a $K$-theoretic framework for  these
 other CFT structures, and show how they relate to D-brane
charges and charge-groups.  We also study conformal
embeddings and the $\bbE_7$ modular invariant of SU(2), as well as some families of
finite group doubles. This new $K$-theoretic framework allows us to 
simplify and extend the less transparent, more ad hoc descriptions of these 
structures obtained previously within CFT.

\end{abstract}

{\footnotesize
\tableofcontents
}

\section{Introduction}

{ The  {representation ring} $R_G$ of a
 compact, connected, simply-connected Lie group $G$}
  can be realized as the equivariant (topological) $K$-group
$K_G^{0}(\hbox{pt})$ of $G$ acting trivially on a point
${\rm pt}$. 
{ For each $k\in\bbZ$, an analogue of $R_G$ for the loop group $LG$} is the 
\textit{Verlinde ring} $\hbox{Ver}_k(G)$, a finite-dimensional quotient of $R_G$ corresponding to
the positive-energy representations of $LG$ { of level $k$.}
 Freed--Hopkins--Teleman \cite{FHT,FHTi,FHTii,FHTlg}
identify the Verlinde ring $\hbox{Ver}_k(G)$ with the twisted equivariant
$K$-group ${}^\tau K^{{\rm dim}\,G}_G(G)$ for some twist
$\tau\in H^3_G(G;\bbZ)$ 
depending on $k$, where $G$ acts on itself by conjugation.
The ring structure of $\hbox{Ver}_k(G)$ is recovered from the
push-forward of group multiplication $m:G\times
G\rightarrow G$, whereas its $R_G$-module structure comes
from the push-forward of the inclusion $1\hookrightarrow G$ of the identity. 

Freed-Hopkins-Teleman interpret their theorem within the context of 3-dimensional
topological quantum field theory, namely as a $K$-theoretic expression 
{ for}  the functorial image of the torus, or equivalently for the Grothendieck ring of the image
of a circle. But it is also a theorem in 2-dimensional conformally invariant quantum  field theory 
(CFT), namely a $K$-theoretic expression for the { fusion coefficients of} the CFT 
describing strings living on $G$. 

{ Conformal field theories are a class of extremely symmetrical 2-dimensional
quantum field theories. Several  profound points of contact between them and mathematics
proper have been developed in the past 30 years or so.}
A CFT consists of two semi-autonomous \textit{chiral} halves { (each of} {
which can be thought of as a vertex operator
algebra (VOA))}; { how} { those halves}
are glued together into the \textit{full} CFT  { is described by what is called there} the 1-loop partition function
{ or} \textit{modular invariant}{. More precisely,} knowing { the} modular invariant  is equivalent to knowing { how the state-space}
(the space  Segal's functor { \cite{Seg}} associates to the circle{) 
 decomposes  into irreducibles of both
chiral VOAs. On the other hand, the Verlinde ring is chiral, being}  the ring of modules of the VOA using a
tensor product called \textit{fusion}. 

{ A} rich structure (full system, alpha-induction 
etc) {associated} to full CFT was realised mathematically in
the { sufferable or braided} subfactor { $N\subset M$} interpretation by B\"ockenhauer, Evans, Kawahigashi, 
Pinto etc (see e.g. the reviews \cite{BE4,BE5}), and transported to the categorical 
setting by Fuchs, Fr\"ohlich, Ostrik, Runkel, Schweigert etc (see e.g. the
review \cite{FFRS}). It is this subfactor formulation, rather than the vertex operator algebra one 
(see e.g. \cite{Ganmd, Gannbook}), which is the inspiration for the current work. {
In this picture, the Verlinde ring is a nondegenerately braided system of endomorphisms on the 
factor $N$,} { while the \textit{full system} (known physically as the fusion ring
of defect lines) consists of endomorphisms on the factor $M$;} { two}
{ natural maps from the Verlinde ring to the full system, built up from the
braiding, are called \textit{alpha-induction}. The \textit{neutral system} is a
nondegenerately braided
subsystem of the full system corresponding to the maximal} { possible VOA.} 
{ The \textit{nimrep} consists of maps $N\rightarrow M$, known as  boundary
states; note that both the Verlinde ring and the full system act on it.
We also have \textit{D-brane charges}, i.e. numerical assignments to the boundary states, i.e. to the
nimrep, which
respect the  Verlinde action.}

{ Interpretting Freed-Hopkins-Teleman within the CFT context suggests
several} related developments, many of which are provided in this paper. { Fix
a group $G$ and level $k$.} The different { full}
CFTs associated to 
$G$ and $k$ are largely parametrised 
by their modular invariants. For  each choice of \textit{viable} or \textit{sufferable} 
modular invariant (that is, {one} realised by a CFT),
{ we have these
 structures of the corresponding full CFT (namely, the full system etc),} { related to but going beyond} 
the Verlinde ring. We seek their $K$-theoretic descriptions.

When $G$ is a finite group, Evans \cite{Ev} { does this.}
Here the twist $\tau$ lies in the group cohomology $H^3_G(\mathrm{pt} ;\bbT)$.
The sufferable modular invariants for the twisted double 
$\mathcal{D}^\tau(G)$ are parametrized by certain pairs $(H,\psi)$ for a subgroup $H$ of
$G\times G$ and $\psi\in H^2_H(\mathrm{pt};\bbT)$ \cite{EP1,Ost}. Let
$H\times H$ act on $G\times G$ on the left and right:
$(h_{\rm L},h_{\rm R})\cdot (g,g')=h_{\rm L}(g,g')h_{\rm R}^{-1}$.
Then $^{\tau,\psi} K^0_{H\times H}(G\times G)$ can be defined and 
\hbox{identified} with the {full system}, and again $^{\tau,\psi} K^1_{H\times H}(G\times G)=0$
(at least when the twists $\tau$ and $\psi$ are trivial). Choosing $H$
to be the diagonal subgroup isomorphic to $G$ recovers the Verlinde ring. 

As warm-up, we work out explicitly here the story for certain classes of finite groups (cyclic
and dihedral). But our real interest is Lie groups.  
From this point of view, Lie groups are much more complicated than finite groups  --- for
example, there is no direct analogue of the $(H,\psi)$
parametrization of { sufferable} modular invariants --- but \cite{EG1} argues that
similar extensions of Freed--Hopkins--Teleman can be
expected. This paper establishes several of these, going far beyond \cite{EG1}.

We start with the case of an $n$-torus, which we work out in complete detail. Its Verlinde ring was 
determined $K$-theoretically in \cite{FHTii} (the bundle picture was described in \cite{EG1}); in
section 3 we  find the full system, nimrep, alpha-induction etc for any modular invariant. We find 
that all modular invariants are sufferable, with a unique { full system, nimrep etc}.  To our knowledge, this result is new.

{ Consider now} $G$ compact, connected, and simply-connected{; then
the} twist $\tau\in H^3_G(G;\bbZ){\,\cong\bbZ}$ { can be identified with} the level.
{ In section 5 we study} the most important source of modular invariants for these $G$,
 the { so-called} \textit{simple-current invariants}{, i.e. modular invariants
 associated to symmetries of the affine Dynkin diagram of $G$}. { The corresponding}
 { full CFTs}
{ describe} strings living on the non-simply-connected
groups $G/Z$ (for $Z$ a subgroup of the centre of $G$).  { For these we develop}
 a complete theory. { In particular,} the full system is given by the
twisted equivariant $K$-theory of $G\times G$ acting diagonally on
$(G/Z_0)\times (G/Z_0)$ for some subgroup $Z_0$ of $Z$. By contrast, $^\tau K_G^{{\rm
dim}\,G}(G/Z)$ for $G$ acting adjointly on $G/Z$ is
the associated nimrep ({this $K$-group  vanishes in degree} dim$(G)+1$), and
$^{\tau'} K_{G/Z_0}^{{\rm dim}\,G}(G/Z_0)$ again for the adjoint action gives the neutral
system. We give the
alpha-inductions below.
To our knowledge, these expressions have not { appeared in the literature in this
generality}.
 { For a concrete example, we work out the case $G={\rm SU}(2)$ and $Z=\{\pm I\}$
 explicitly.}  

{ A representation-theoretic expression for the} nimrep { of} these simple-current modular invariants had been conjectured
in \cite{BFS,GaGa2}; in Theorem 3 below we prove this conjecture and rewrite
{ the} nimrep in a considerably simpler way. This is { one of the main
results}  of the paper. We were led to this description
by trying to match the conjectured { expression for the}
nimrep with the expected $K$-group.

{ Fix again a compact, connected, simply-connected group $G$. The generic  modular invariants of $G$ are expected to be the simple-current invariants and their twists by outer
automorphisms (i.e. for each $G$ we expect  only finitely many additional modular invariants,
over all levels $k$).}
We give in section{s 4.2 and 5.3} a complete $K$-homological  description for the 
{ outer 
automorphism twists of simple-current invariants for any subgroup $Z$ (possibly trivial)} of the centre of $G$.
 We compare our $K$-homological descriptions
with the analogous descriptions (when they exist) coming from conformal field 
theory (see e.g. \cite{GaGa,GaGa1,GaGa2,GaGa3}). They match beautifully.
Nothing like this to our knowledge { has} appeared in the literature. { These
systematic $K$-theoretic expressions} should be contrasted with the ad hoc { 
expressions} in the
CFT literature.
 
 Twisted K-theory formulations of D-brane charge groups have long been known. We show these fit naturally within
 our framework: in particular the assignment of D-brane charges amounts to the forgetful functor from the nimrep to the charge group. To our knowledge, this has not appeared in the literature before { in this generality}.
 
 { For example, the modular invariants for $G=\,$SU(2)  fall into an A-D-E pattern.
 In this familiar case, the only modular invariants which  remain} { unaccounted for}
 are the three { nongeneric, i.e.} exceptional,
{ ones} { $\bbE_6,\bbE_7,\bbE_8$}. Of these, two (namely $\bbE_6$ and $\bbE_8$) are due to \textit{conformal 
 embeddings}. { These arise when the modular invariant is of pure extension type (these 
 are called \textit{type I}) and the maximal chiral extension corresponds to a group $H\ge G$  at level $1$.}
 In \cite{EG1} we approximated the full system for any
 {conformal embedding} by
{ $^\tau K^{{\rm dim}\,{G}}_{{G}}(H)$} where the { sub}group ${G}$ acts adjointly 
on $H$. This is clearly a part of the story here: we could see the level change in the bundles in going to the
subgroup (e.g. from $G_2$ level 1 to SU$(2)$ level 28), and we observed that McKay's
A-D-E name for the largest finite stabiliser for $H$-conjugation on $G$ matched the
Cappelli-Itzykson-Zuber name for the corresponding modular invariant. 
In this paper we propose (in the spirit of \cite{Ev}) to realize
the full system by ${}^\tau K^{G\times G}_0(H\times H)$ where now
the action is given by $(g,g').(h,h')=(ghg'{}^{-1},gh'g'^{-1})$.
This tightens the match with the full system. Finally,
we can $K$-theoretically identify most of the full system for the remaining SU(2) modular invariant
$\bbE_7$. See section 6 for the details of our treatment of these exceptional modular invariants. 

The point of our work is not merely to translate everything in CFT into $K$-theory,
but rather to demonstrate that the latter provides systematic tools for understanding
the former. For example, as explained in section 4.2, the easiest way to
prove the conjectured picture in CFT for nimreps associated to outer
automorphisms would be to relate it to our $K$-homological description. Indeed
as already mentioned, we were led to our proof ({ see} Theorem 3 below)
that the conjectured nimrep
for simple-current invariants \textit{is} a nimrep, by $K$-homological 
considerations. Or consider a 
 modular invariant $\cZ$. It may or may not correspond to an CFT; for example
in general there will not exist the necessary extra structure, e.g. even a nimrep, 
compatible with it. Subfactor theory elegantly describes these extra 
structures but there still remains the problem of showing existence. 
For example, to establish that all modular invariants for SU(3) \textit{are} 
realised by  subfactors, { Evans-Pugh}
\cite{EvPui,EvPuii} start with the nimrep graphs, construct for them the
Boltzmann weights, construct from this the subfactors, compute the nimreps 
and recover the graphs they started with. What is needed is 
something more systematic, which for instance does not need to know what the 
nimrep graphs are expected to be.
A motivation for our present and future $K$-theory work is to help 
provide tools for 
showing that certain modular invariants are realised by CFT { (i.e are sufferable)} 
and { for finding} the 
extra structure. Here by $K$-theory we really mean the concepts of vector 
bundle/Fredholm module/cycle,  not  just the equivalence classes (i.e. 
$K$-groups themselves) which would only produce graphs. This point is also
discussed at the end of \cite{EG1}. { In addition, we hope this more conceptual
approach to the structures of full CFT will help interpret the so-called exotic modular
tensor categories (associated to e.g. the Haagerup subfactor) as combinations
of loop group and finite group systems (see  \cite{EG2}).}

The modular invariant is an integral matrix indexed by the
primaries (the preferred basis in the Verlinde ring), and as such is
a linear map between $K$-groups. As proposed in \cite{EG1}, it should be understood as an element of $KK$-theory.
Likewise for alpha-induction and sigma-restriction. In the special case of the double of
finite groups (with trivial $H^3$ twist), a natural basis can be found in which the modular
group  acts by permutation matrices, and so in this case these matrices can also be
interpreted as $KK$-elements. Given the accomplishments of this paper, developing these 
pictures is  the natural next step.

\section{Review and notation}

\subsection{Groups, representations and cohomology}

See \cite{EG1} for details and references. For any compact, finite-dimensional  group $G$,
we write $R_G$ for the representation ring (equivalently, character ring) of 
the group $G$, the span over $\bbZ$ of the
(isomorphism classes of) irreducible representations ($=\,$\textit{irreps}).
For arbitrary compact, connected, simply-connected  Lie  $G$ of rank $r$,
we write $\rho_\lambda$ for the irrep with highest weight $\lambda\in\bbZ_{\ge 0}^r
:= P_+(G)$. 
Write  $\Lambda_1,\ldots,\Lambda_r$ for the fundamental weights;
then as is well-known $R_G$ is the polynomial ring $\bbZ[\rho_{\Lambda_1},\ldots,
\rho_{\Lambda_r}]$.

It is convenient to introduce special notation for three common $R_G$:
\begin{align}
R_{{\rm SU}2}&=\bbZ[\sigma]=\mathrm{span}\{\sigma_1,\sigma_2,\ldots\},\nonumber\\
R_{\textrm{O2}}&=\bbZ[\delta,\kappa]/(\delta\kappa=\kappa,\delta^2=1)=\mathrm{span}\{1,\delta,
\kappa_1,\kappa_2,\ldots\},\nonumber\\
R_{\bbT}&=\bbZ[a^{\pm
1}]=\mathrm{span}\{1,a,a^{-1},a^2,a^{-2},\ldots\},
\nonumber
\end{align}
where $\sigma_i$ is the $i$-dimensional ${\rm
SU}(2)$-representation (so $\sigma= \sigma_2$ is the
defining representation), $\delta=\hbox{det}$, $\kappa_i$
is the two-dimensional O(2)-representation with winding
number $i$ (so $\kappa =\kappa_1$ is the defining
representation), and $a^i$ is the one-dimensional
representation for the circle ${\rm
SO}(2)={\rm U}(1)=S^1=\mathbb{T}$ with winding number $i$.

It is convenient to extend the notation to $\sigma_0=0$, $\kappa_0=1+\delta$, and 
$\sigma_{-n}=-\sigma_n$ and $\kappa_{-n}=\kappa_n$ for all $n\in\bbZ_{>0}$.
Some inductions we need are Dirac induction
 D-Ind$_T^{\mathrm{SU2}}$, which sends $a^i$ to $\sigma_{i}$, and
Ind$_T^{\textrm{O2}}$, which sends $a^i$ to $\kappa_{i}$.

When the fundamental group $\pi_1(G)$ is not trivial, we can define
\textit{spinors}. Consider $G=\mathrm{SO}(2n+1)$, the
quotient of Spin$(2n+1)$ by its centre $Z=\{1,z\}$: 
the ${\mathrm{Spin}(2n+1)}$-irreps have highest weight  
$\lambda=(\lambda_1,\ldots,\lambda_n)\in\bbZ_{\ge 0}^n$; the value of
the $\lambda$-irrep on $z\in Z$ is $(-1)^{\lambda_n}$. The 
representation ring $R_{\mathrm{SO}(2n+1)}$ is the $\bbZ$-span of the irreps
with ${\lambda}_n$ even;  the spinors are
the irreps with $\lambda_n$ odd, and their $\bbZ$-span $R^-_{\mathrm{SO}(2n+1)}$
is a module for $R_{\mathrm{SO}(2n+1)}$.

Suppose $\phi:H\rightarrow G$ is a group homomorphism. Write $\Delta_H$
for the diagonal subgroup $\{(\phi(h),\phi(h)):h\in H\}$ of $G\times G$.
Write `$H^{\mathrm{ad}}$ on $G$'
for the adjoint action of $H$ on $G$, where $h.g=\phi(h)g\phi(h)^{-1}$,
and write `$H^L$ on $G$' (resp. `$H^R$ on $G$') for the left (right) action $h.g=\phi(h)g$ (resp. $h.g=g
\phi(h^{-1}))$.

For reasons explained next subsection, we are interested in equivariant cohomology, in particular
$H^1_G(X;\bbZ_2)$ and $H^3_G(X;\bbZ)$  when $G$ is a group acting on $X$. These are defined {through the Borel construction $H^n_G(X;A)=H^n((E_G\times
X)/G;A)$. Here  $EG$ is a contractible space on which $G$ acts freely, and the quotient is taken for the diagonal action, and where  $B_G=E_G/G$ is the classifying space of $G$, so that $H^*_G(point) \simeq H^*(BG)$. 
If  $N$ is a normal subgroup of $G$ acting freely on $X$, then we have 
the natural isomorphism 
\begin{equation}\label{freenormcoh}
H^*_G(X;A)\simeq H^*_{G/N}(X/N;A)\,.
\end{equation}

It is often convenient to compute 
these cohomology groups using spectral sequences via the fibration 
$X\rightarrow(E_G\times X)/G\rightarrow B_G$. Then $E_2^{p,q}=H^p_G(\mathrm{pt};
H^q(X;\bbZ))$  converges to $H^*_G(X)$. When $G$ is compact, connected, and simply-connected, 
$H^*_G(\mathrm{pt};\bbZ)\simeq \bbZ\oplus 0\oplus 0\oplus 0\oplus\bbZ\oplus\cdots$ and
$H^*(G;\bbZ)=\bbZ\oplus 0\oplus 0\oplus \bbZ\oplus\cdots$, while $H^*_G(\mathrm{pt};\bbZ_2)\simeq \bbZ_2\oplus 0\oplus 0\oplus\cdots$
and $H^*(G;\bbZ_2)\simeq \bbZ_2\oplus 0\oplus\cdots$.
Spectral sequence computations then identify both $H^3_{G^{\mathrm{ad}}}(G;\bbZ)$ and $H^3_{G^{\mathrm{ad}}}(G\omega;\bbZ)$ with $H^0_G(\mathrm{pt};H^3(G;\bbZ))=\bbZ$, where $\omega$ is an outer automorphism and 
$G\omega \subset G \rtimes \langle \omega \rangle$, with the natural adjoint action on the semi-direct product.
Similar computations show that, for $Z$ a subgroup of the centre of $G$,
\begin{equation}\label{H3H1}
H^3_{G^{\rm{ad}}}(G/Z;\bbZ)\simeq H^3_{G^{\mathrm{ad}}\times Z^L}(G;\bbZ)\simeq 
\bbZ\,, \ \ 
H^1_{G^{\rm{ad}}}(G/Z;\bbZ_2)\simeq \mathrm{Hom}(Z,\bbZ_2)\,.\end{equation} 

Given any nontrivial $\epsilon\in H^1_G(\mathrm{pt};\bbZ_2)$, i.e. a continuous group
homomorphism $\epsilon:G\rightarrow\bbZ_2$, the subgroup $H:=\mathrm{ker}\,\epsilon$
of $G$ is index 2. There are two types of $G$-irreps $\rho$ (see \cite[Sect.2.1]{EG1} for
details): \smallskip

\noindent\textit{type} ${2\atop 1}$: there is an inequivalent $G$-irrep $\rho'$ such that
$\mathrm{Res}_H^G\rho\simeq \mathrm{Res}_H^G\rho'$ is irreducible;

\smallskip\noindent\textit{type} ${1\atop 2}$: Res$_H^G\rho=\rho_1\oplus\rho_2$ where 
$\rho_i$ are $H$-irreps.\smallskip

By the \textit{graded representation ring} $^\epsilon R_G$ 
\cite[Sect.4]{FHT},\cite[Sects.1.2,2.1]{EG1} we mean the $\bbZ$-span of formal
differences $\rho_1\ominus\rho_2=-(\rho_2\ominus\rho_1)$ over all type ${1\atop 2}$ $G$-irreps $\rho$.
By $^\epsilon R^1_G$ we mean the $\bbZ$-span of formal anti-symmetrisations $\rho^-:=
(\rho-\rho')/2=-\rho^{\prime -}$ over all type ${2\atop 1}$ $G$-irreps $\rho$.

\subsection{Twisted equivariant $K$-theory}

We are interested in $K$-groups $^\tau K^*_G(X)$ and $K$-homology $^\tau K_*^G(X)$,
where  $G$ is a compact group acting on topological space  $X$
and the twist $\tau$ lies in  the equivariant graded Brauer group} $H^1_G(X;
\bbZ_2)\times H^3_G(X;\bbZ)$   which classify graded $\bbZ_2$ $G$-equivariant continuous trace algebras with spectrum $X$. Here, `$\times$' denotes Cartesian product; the group structure
  involves  the Bockstein homomorphism, see e.g. \cite[eqn.(1.4)]{EG1}.  Strictly 
 speaking, we should use the label $\tau$ for the graded $G$-equivariant bundle of compacts on $X$ rather than its equivalence class. The corresponding twisted equivariant $K$-groups for equivalent bundles are isomorphic but not canonically.  For
the definition, basic properties, and original references of twisted equivariant $K$-theory, see
\cite{AS,Kar1,EG1}. 

\textit{Bott periodicity} says ${}^\tau K^{i+2}_G(X)={}^\tau K^i_G(X)$ and 
${}^\tau K_{i+2}^G(X)={}^\tau K_i^G(X)$, for all $G,X,i$. The additive groups
 ${}^\tau K^*_G(X)$ and ${}^\tau K_*^G(X)$ carry an $R_G$-module structure.
  When $\epsilon\in 
H^1_G(\mathrm{pt};\bbZ_2)$, then ${}^\epsilon K^*_G(\mathrm{pt})={}^\epsilon R_G\oplus
{}^\epsilon R^1_G$, where the graded representation rings were defined last subsection.
On the other hand, twisting by $H^3_{G}(\mathrm{pt};\bbZ)$ introduces spinors.  For example, 
$H^3_{\mathrm{SO}(n)}(\mathrm{pt};\bbZ)\simeq\bbZ_2$, which we label as $\pm$. The for the twist labelled $-$, we have  $^-K^*_{\mathrm{SO}(n)}(\mathrm{pt})\simeq R^-_{\mathrm{SO}(n)}$.

We freely interchange $K$-theory and $K$-homology, using whichever is more
convenient or appropriate, through
\textit{Poincar\'e duality} \cite{Tu,EEK,BMRS}, which for a compact manifold $X$ says 
\begin{equation}\label{PD}
{}^\tau K_i^G(X)\, \simeq {}^{\tau'} K^{i+{\rm dim}\,X}_G(X) \,,
\end{equation}
where $\tau+\tau'=({\rm sw}_1^G(X),{\rm sw}_3^G(X))\in H^1_G(X;
\bbZ_2)\times H^3_G(X;\bbZ)$. Here $sw_1^G$ is the  $G$-equivariant Stiefel--Whitney class and 
$sw_3^G$ is the integral $G$-equivariant Stiefel-Whitney class.

For $X$ a compact manifold fixed pointwise by $G$, and $H$
a subgroup of $G$, we know $K^*_G(X\times G/H)
=K^*_H(X)$ (see \cite{Se}) and hence
\begin{equation}\label{cancel}
K_*^G(X\times G/H) \, \simeq \, {}^{\tau'} K_{*+{\rm
dim}\,G+{\rm dim}\,H}^H(X)
\end{equation}
by Poincar\'e duality, for the appropriate twist $\tau'$.
If $N$ is a normal subgroup of $G$, and $N$ acts freely on
$X$, then using \eqref{freenormcoh}
\begin{equation}\label{freenormal}
{}^\tau K_*^G(X) \, \simeq  \, {}^\tau K_*^{G/N}(X/N)\,.
\end{equation}
If instead $X$ is fixed by $H$ and  the twist $\tau\in H^1_{G\times H}(X;\bbZ_2)\times
H^3_{G\times H}(X;\bbZ)$ lies in $(H^1_G(X;\bbZ_2)\otimes 1)\times(H^3_G(X;\bbZ)\otimes 1)$
then
\begin{equation}\label{fpkth}
^\tau K_{G\times H}^*(X) \, \simeq \, R_H\otimes_\bbZ {}^\tau K_G^*(X)\,.\end{equation}

Let $U$ be a $G$-invariant open subset of $\textit{X}$, and
 $\tau^{\prime}$ and $\tau^{\prime\prime}$ be
the restrictions of the twist $\tau$ on $X$ to $U$ and $X/U$,
respectively. Then the \textit{six-term exact sequence} says:
\begin{equation}\label{six-term}
\begin{array}{ccccc}
^{\tau^{\prime}}K^{G}_0(U) & {\longleftarrow} &
^{\tau}K^{G}_{0}(X) & {\longleftarrow} &
^{\tau^{\prime\prime}}K^{G}_{0}(X/U)\\[3pt]
\downarrow & & & & \uparrow\\[2pt]
^{\tau^{\prime\prime}}K^{G}_1(X/U) & {\lori} &
^{\tau}K^{G}_1(X) &{\lori} &^{\tau^{\prime}}K^{G}_1(U)
\end{array}.
\end{equation}

The Hodgkin spectral sequence (see e.g. \cite[Sect 3.3]{EG1}) explains how to restrict
$G$-equivariance in $K$-theory to a subgroup $H$. It
starts from a closed subgroup $H$ of a Lie group $G$ with torsion-free $\pi_1$, and a $G$-action 
on space $X$, and defines $E_2^{p,q}=\textrm{Tor}^p_{R_G}(R_H,{}^\tau K_G^q(X))$;
the resulting spectral sequence of $R_G$-modules strongly converges to
$^{\tau'} K^*_H(X)$, where $\tau'$ is the image of $\tau$ under the natural map
$H^3_G(X;\bbZ)\rightarrow H^3_H(X;\bbZ)$ restricting $G$-equivariance to $H$.

It is useful to have explicit descriptions of the Dixmier-Douady bundles representing
the twists relevant to these  $K$-groups $^\tau K_G^*(X)$. 
These are  bundles over $X$ with fibre the compact operators $\cK=\cK(\cH)$ on a
 $G$-Hilbert space $\cH$, where we usually require $G$-stability $\cH\simeq \cH\otimes L^2(G)$ (i.e. each
$G$-irrep appears in $\cH$ with infinite multiplicity). Many of these 
were constructed in \cite[Sect.2.2]{EG1}; we'll quickly sketch two of them.  Consider 
first $^k K_{\bbT^{\mathrm{ad}}}^*(\bbT)$ for the circle $\bbT=\bbR/\bbZ$, where $k$ is in $\bbZ \simeq H^3_\bbT (\bbT)$.  Here
we can simply take $\cH=L^2(\bbT)$. Let $U_k\in U(\cH)$ be the unitary operator defining the equivalence
$\pi\otimes a^k\simeq\pi$, where $\pi$ is the regular representation of $\bbT$ and 
$a^k\in R_\bbT$. The sections of our bundle are maps $f:\bbR\rightarrow\cK$, continuous
on $[0,1]$, satisfying $f(t)=U_kf(t+1)U_k^*$ for all $t\in\bbR$. The $\bbT$-action is defined
by $(t.f)(s)=\mathrm{Ad}(\pi(t))f(s)$. 

Consider next $^k K_{G^{\mathrm{ad}}}^*(G)$ for $G=\mathrm{SU}(2)$, and  $k$ in $\bbZ \simeq H^3_{SU(2)}(SU(2))$. Here, $\cH=L^2(G)\otimes\ell^2$,
with regular $G$-representation $\pi\simeq \sum_n\sigma_n\otimes 1_\infty$ where $\sigma_n$ denotes the $n$-dimensional irreducible representation of $SU(2)$ and $1_\infty$, the identity operator on $\ell^2$.
The conjugacy classes of $G$ are parametrised by the Stiefel diagram $S$, the points  in half 
of a maximal torus $T$ (say diag$(e^{\i  t},e^{-\i t})$ for $0\le t\le \pi$). 
Cover $G$ with $G$-invariant open sets $D_1,D_2$ containing $I,-I$
respectively. On each fibre $(x,c)_i$, $x\in D_i$, $c\in\cK$, we define the $G$-action
$g.(x,c)_i=(gxg^{-1},\mathrm{Ad}(\pi_g)c)_i$. Define $U_k\in U(\cH)$ as in Section 2.2 of  \cite{EG1};
restricting $\pi$ to the maximal torus $T$, i.e. considering weight spaces $V_{m,n-1}$, $m=n,
n-2,\ldots,2-n,-n$ (for convenience write $V_{m,-n+1}$ when $m<0$), $U_k$ acts as the shift 
operator $V_{m,n-1}\rightarrow V_{m-k,n-1-k}$. For $x\in D_1\cap D_2\cap  S $, the gluing condition
is $(gxg^{-1},c)_1=(gxg^{-1},\mathrm{Ad}(\pi_gU_k\pi_g^{-1})c)_2$. Then
 for $x\in D_1\cap D_2\cap S $, when $g\in Z_G(x)=T$, $\pi_gU_k\pi_g^{-1}=
\lambda_kU_k$ where $\lambda_k$ is the character $a^k\in R_T$ (consistency requires
that it be a character).

Some of our \textit{nimreps} involve
the product $^\tau K_{0}^{G^{\mathrm{ad}}}(G)\times {}^{\tau'}  K_{0}^{G^{\mathrm{ad}}}(G/Z)
\rightarrow {}^{\tau'}  K_{0}^{G^{\mathrm{ad}}}(G/Z)$, for $Z$ a
subgroup of the centre of $G$ (the relation between the twists $\tau,\tau'$ is described in
section 5.2), which arises through the push-forward of the products 
$g*g'\mapsto gg'$ and $g*h\mapsto \phi(g)h$ on the
spaces (where $\phi$ is the projection $G\rightarrow G/Z$) --- these are clearly $G$-equivariant. 
It will sometimes be convenient to rewrite $ {}^{\tau'}  K_{0}^{G^{\mathrm{ad}}}(G/Z)$ as
${}^{\tau'}K_{0}^{G^{\mathrm{ad}}\times Z^L}(G)$, but the product in this case is again given
by $g*g'\mapsto gg'$. Some of our \textit{full systems}  involve a product on $^\tau 
K_{0}^{G^{\mathrm{ad}}\times Z^L\times Z^R}(G)$ or
${}^{\tau}  K_{0}^{G^{\mathrm{ad}}\times Z^L}(G)$, making them into rings;
this arises through the convolution product $g*g'\mapsto \oplus_{z\in Z}gzg'$ (at least when
$Z$ is finite). The product $^\tau K^{0}_{G_1^L\times G_2^R}(H\times H)\times
{}^{\tau'} K^0_{G_2^L\times G_3^R}(H\times H)\rightarrow
{}^{\tau''} K^0_{G_1^L\times G_3^R}(H\times H)$, for subgroups $G_i$ of
a \textit{finite} group $H\times H$, arises through the 
convolution product $(h_1,h_2)*(h_1',h_2')\mapsto \oplus_{g\in G_2}(h_1,h_2)g(h_1',h_2')$
---  the relation of the twists $\tau,\tau',\tau''$ is explained in \cite{Ev}. 
For Lie groups this product  is going to be more subtle,   see e.g. 
the work of  Jeffrey--Weitsman  \cite{JW} on the product of  conjugacy classes.

In twisted $K$-theory, the naive product would add the twists.
A twist-preserving product arises for instance when 
 the twist is transgressed \cite{TX}, i.e. lies in the image of the appropriate map 
$H_G^i(\mathrm{pt};A)\rightarrow H_G^{i-1}(X;A)$. A twist $\tau\in H^3_G(G;\bbZ)$
only determines the $K$-theory $^\tau K_G^*(G)$ as an
additive group; the ring structure comes from a choice of
lift (if it exists) of the 3-cocycle $\tau$ to
$H^4_G(\hbox{pt};\bbZ)$. But for $G$ compact, connected, 
simply-connected, transgression identifies $H^4_G(\hbox{pt};\bbZ)$
with $H^3_G(G;\bbZ)$ so $\tau$ parametrizes the full ring structure.

We can see directly the action of the $k$th roots of unity in $\bbT$ on
the $K$-group $^kK_\bbT^1(\bbT)$. For this purpose,
instead of placing a copy of the twisting unitary $U_k=U_1^k$ at each integer,
it is more convenient to place $U_1$ at each point in $\frac{1}{k}\bbZ$,
as was done at the end of \cite[Sect.3.1]{EG1}. A section now consists of
$k$ continuous maps $f_l:[0,\frac{1}{k}]\rightarrow\cK$ for $l\in\bbZ_k$, 
related by $f_l(\frac{1}{k})=\mathrm{Ad}(U_1)f_{l+1}(0)$ for all $l$. 
Recall from \cite[Sect.3.1]{EG1} the
identification of the basis of $^kK_\bbT^1(\bbT)$ as roots of unity in $\bbT$.
We get a $\bbZ_k$-action on this bundle, and hence on $^kK_\bbT^1(\bbT)$,
by sending section $(f_0,\ldots,f_{k-1})$ to $(f_l,f_{l+1},\ldots,f_{l+k-1})$, 
or equivalently
fibre $(t,c)$ to $(t+\frac{l}{k},c)$.

\subsection{Rational conformal field theory and subfactors}

Rational conformal field theories (RCFT) can be formulated within an algebraic theory of vertex operator algebras (VOA's), see e.g. \cite{Ganmd}, or conformal nets of subfactors see e.g. \cite{LR, xu, BE5, BE6, {rehren:2000}, Yasu}.
For the case of Wess-Zumino-Witten or loop group models, either formulation is equally valid or sufficient. 
For constructing possibly more exotic models, subfactor theory has thrown up interesting candidates \cite{EG2}.

In the VOA picture the chiral data of an RCFT consists of a rational vertex operator
algebra (VOA) $\cV$ whose category of modules forms a modular tensor category. In subfactor theory, a factor $N$ is obtained as a local factor $N=N(I_{\circ})$ of a conformally covariant quantum field theoretical net of factors $\{N(I)\}$ indexed by proper intervals $I \subset \mathbb{R}$ on the real line arising e.g. from current algebras defined in terms of local loop group representations. The $N$-$N$ system $\NXN$ is obtained as restrictions of Doplicher-Haag-Roberts morphisms (cf. \cite{haag:1992}) to $N$. We write   $\Sigma(\NXN)$ for the endomorphisms which decompose into a finite number of irreducibles
from $\NXN$. Typically these endomorphisms will be braided, indeed non-degenerately braided{,} and again will yield a modular tensor category.

The
finitely many irreducibles  (the simple objects in the modular tensor category) are
called the \textit{primaries} $\lambda\in\Phi$. The (commutative associative semi-simple)
Grothendieck ring of this category,
with preferred basis $\Phi$, is called the \textit{fusion ring} or \textit{Verlinde ring} Ver{.
Write} $\lambda\mu=\sum_{\nu\in\Phi}N_{\lambda,\mu}^\nu\nu$. The \textit{vacuum} 
primary (corresponding to $\cV$-module $\cV$ itself in the VOA picture or the identity endomorophism in the subfactor setting), is the unit in Ver, denoted 1. 
The \textit{simple-currents} $j\in\Phi$ are the invertible primaries, i.e. those for which there is
a $j'\in\Phi$ such that $jj'=1$. Multiplication in Ver 
by a simple-current permutes
the primaries. The simple-currents form a  group, by composition of those
permutations; if $j,j'$ are two simple-currents, then
\begin{equation} N_{j\lambda,j'\mu}^{jj'\nu}=N_{\lambda,\mu}^\nu\,.\label{scfus}
\end{equation}

The modular tensor category comes with unitary representations of the braid
groups and mapping class groups. In particular we get a unitary representation of
the modular group SL$(2,\bbZ)$ on the complexification $\bbC\otimes_\bbZ\mathrm{Ver}$:
write $\left({0\ -1\atop 1\ 0}\right)\mapsto S$ and $\left({1\ 1\atop 0\ 1}
\right)\mapsto T$, a diagonal matrix. \textit{Charge-conjugation} $C=S^2$ is an involution on
Ver, 
which permutes the primaries; then $S_{C\lambda,\mu}=
S_{\lambda,C\mu}=S_{\lambda,\mu}^*$ and $T_{C\lambda,C\mu}=T_{\lambda,\mu}$ for all 
primaries $\lambda,\mu$. For any simple-current $j$, there is a function $Q_j:\Phi\rightarrow
\bbQ$ such that $S_{j\lambda,\mu}=\exp[2\pi\i Q_j(\mu)]S_{\lambda,\mu}$. 
The simple-currents are precisely those primaries $j\in\Phi$ whose
\textit{quantum-dimension} $S_{j,1}/S_{1,1}$ equals 1. The
\textit{Verlinde formula} says
that the matrices $N_\lambda=(N_{\lambda,\mu}^\nu)$ are simultaneously diagonalised
by $S$ and have eigenvalues $S_{\lambda,\mu}/S_{1,\mu}$, of multiplicity 1 for
each $\mu\in\Phi$. The Perron-Frobenius eigenvalue of the nonnegative 
matrix $N_\lambda$ is the quantum-dimension $S_{\lambda,1}/S_{1,1}$. 

\medskip\noindent\textbf{Definition 1.} \textit{A matrix $\cZ=(
\cZ_{\lambda,\mu})_{\lambda,\mu\in\Phi}$ is called a} modular invariant
\textit{if $\cZ S=S\cZ$, $\cZ T=T\cZ$, each entry $\cZ_{\lambda,\mu}$
is a nonnegative integer, and $\cZ_{1,1}=1$.}\medskip

The \textit{modular invariant} describes how the space of states of the RCFT carries
a representation of the modular tensor category. 
For example, $\cZ=I$ and
$\cZ=C$ are both modular invariants; when $\cZ$ is a permutation matrix it is
called an \textit{automorphism invariant}. Automorphism invariants $\cZ$ are special among
modular invariants in that matrix multiplication $\cZ\cZ'$ or $\cZ'\cZ$ are modular invariants
iff $\cZ'$ is.

This data (the basis $\Phi$, the ring Ver, 
 and the matrices $S,T,\cZ$)
helps define the \textit{bulk CFT}, which describes \textit{closed} string
theory. \textit{Boundary CFT} (see e.g.
the review \cite{PZ}), describing \textit{open} strings,  starts with a finite set of \textit{boundary states} $x\in\cB$.
The Verlinde ring acts on the $\bbZ$-span of these $x$, and this module
structure is called a \textit{nimrep}, written $\lambda x=\sum_{y\in\cB}\cN_{\lambda,x}^yy$.
More precisely:

\medskip\noindent\textbf{Definition 2.} \textit{A set of matrices $\cN_\lambda
=(\cN_{\lambda,x}^y)_{x,y\in\cB}$, for all $\lambda\in\Phi$, is called a}
nimrep \textit{if all entries $\cN_{\lambda,x}^y$ are nonnegative integers,
$\cN_\lambda\cN_\mu=\sum_{\nu\in\Phi}N_{\lambda,\mu}^\nu\cN_\nu$, 
$\cN_{C\lambda}=\cN_\lambda^t$, and given any $x,y\in \cB$ there is a $\lambda
\in\Phi$ such that $\cN_{\lambda,x}^y\ne 0$. Two nimreps $\cN,\cN'$ are called}
equivalent \textit{if there is a permutation matrix $P$ such that $P\cN_\lambda
P^{-1}=\cN'_\lambda$.}\medskip

$P$ defines the identification of the boundary states $\cB$ and $\cB'$.
An analogue of the Verlinde formula holds, namely
\begin{equation}
\cN_{\lambda,x}^y=\sum_{(\mu,i)}\Psi_{x,(\mu,i)}\frac{S_{\lambda,\mu}}{S_{1,\mu}}
\Psi^*_{y,(\mu,i)}\,,\label{nimver}\end{equation}
saying that the matrices $\cN_\lambda=(\cN_{\lambda,x}^y)$ can be simultaneously 
diagonalised by some unitary matrix $\Psi$, with eigenvalues
$S_{\lambda,\mu}/S_{1,\mu}$, where now each $\mu\in\Phi$ comes with some
multiplicity $m_\mu$ (possibly 0) --- this multi-set is called the 
\textit{exponent}
of nimrep $\cN$. In particular, the vacuum has multiplicity 1; the nonnegative
matrix $\cN_\lambda$ will have Perron-Frobenius eigenvalue equal to the
quantum-dimension $S_{\lambda,1}/S_{1,1}$. For example, the matrices $\cN_j$
associated to any simple-current $j$ will be a permutation matrix (since
the only  diagonalisable nonnegative integer matrices with largest
eigenvalue 1 are permutation matrices). 

The nimrep will be compatible with the modular invariant $\cZ$,
in the sense the eigenvalue multiplicity $m_\mu$ equals $\cZ_{\mu,\mu}$ for all
$\mu\in\Phi$. For example the Verlinde ring is itself a nimrep, called
the \textit{Verlinde nimrep}, with 
$\cN_\lambda=N_\lambda$ and $\cB=\Phi$, and is compatible in this sense with the 
diagonal modular invariant $\cZ=I$.

Subfactor theory captures and enhances this data from bulk and boundary CFT.
We refer to \cite{J1,EKaw} for the basic theory of subfactors etc, 
and \cite{BEK1,BE5,BE6} for the theory of alpha-induction etc.

We start with a factor $N$ on which the Verlinde ring acts through the system of endomorphisms 
$\NXN$. The data of full and boundary CFT requires a  subfactor $N\subset M$.
Let $\iota:N\rightarrow M$ be the inclusion and $\overline{\iota}:M\rightarrow N$ its conjugate.
Then $\theta=\overline{\iota}\iota$ is called the \textit{canonical endomorphism} and $\gamma=
\iota\overline{\iota}$ the \textit{dual canonical endomorphism}. {We require
 $\theta$ to be in $\Sigma(\NXN)$.
Using the braiding 
or its opposite, 
 we can lift an endomorphism $\lambda\in\NXN$
of $N$ to one of $M$: $\alpha^{\pm }_\lambda$.  Then $\cZ_{\lambda,\mu}:=\langle \alpha^+_\lambda,\alpha^-_\mu\rangle$
is a modular invariant \cite{BEK1, E2}. The induced $\alpha^\pm(\NXN)$
generate  the \textit{full system} $\MXM$. These maps $\alpha^\pm$ from
the Verlinde ring $\NXN$ to the full system $\MXM$ are called \textit{alpha-inductions};
they preserve multiplication, addition, and (charge-)conjugation. By $\NXM$ (resp. $\MXN$) we 
mean all irreducible  sectors  appearing in any $\iota\lambda$ (resp. $\lambda\overline{\iota}$) for
$\lambda\in\NXN$. The {nimrep} is the $\NXN$ action on  $\NXM$, given by composition; it
will automatically be compatible with the given modular invariant 
\cite{BEK1,BEK2}.

Unlike $\NXN$, the full system $\MXM$ is not in general non-degenerately
braided. However its subsystem $\MXM^0$, defined to be the intersection of all
irreducible sectors in $\alpha^+(\NXN)$ with all those in $\alpha^-(\NXN)$,
 is called the \textit{neutral system}, and just as $\NXN$ captures the 
chiral data of some VOA $\cV$, $\MXM^0$ captures that
of the maximally extended VOA $\cV'\subseteq\cV$ compatible with the
modular invariant $\cZ$. \textit{Sigma-restriction} $\sigma_\beta:=\overline{\iota}\circ\beta
\circ\iota$ maps the full system $\MXM$ to $\Sigma(\NXN)$; it preserves addition and 
conjugates but not multiplication (after all, $\sigma$ takes the identity in $\MXM$ to the canonical
endomorphism $\theta$, not the identity in $\NXN$). Applied to sectors in the neutral system
$\MXM^0$, it coincides in the VOA language with the restriction (\textit{branching rules})
of $\cV'$-modules to $\cV$. } When a modular invariant $\cZ$ is realised by at least one subfactor, and hence with corresponding  full system, nimrep etc, we call
it \textit{sufferable} \cite{BE1,BE4,BE5,BE6,BEK1,BEK2,E2,EP1,EvPui,EvPuii,O,ocn4,xu}.

In the case of a conformally covariant quantum field theoretical net of factors $\{N(I)\}$ and taking two copies of such a net and placing the real axes on the light cone, this defines a local conformal net $\{A(\mathcal{O})\}$, indexed by double cones $\mathcal{O}$ on two-dimensional Minkowski space (cf. \cite{rehren:2000} for such constructions).
A braided subfactor $N \subset M$, determining in turn two subfactors $N \subset M_{\pm}$ obeying chiral locality, will provide two local nets of subfactors $\{N(I) \subset M_{\pm}(I) \}$.
Arranging $M_+(I)$ and $M_-(J)$ on the two light cone axes defines a local net of subfactors $\{A(\mathcal{O})\subset A_{\mathrm{ext}}(\mathcal{O})\}$ in Minkowski space.
The embedding $M_+\otimes M_-^{\mathrm{op}} \subset B$ gives rise to another net of subfactors $\{A_{\mathrm{ext}}(\mathcal{O}) \subset B(\mathcal{O})\}$, where the conformal net $\{B(\mathcal{O})\}$ satisfies locality.
As shown in \cite{rehren:2000}, there exist a local conformal two-dimensional quantum field theory such that the coupling matrix $Z$ describes its restriction to
the tensor products of its chiral building blocks $N(I)$.
There are chiral extensions $N(I) \subset M_+(I)$ and $N(I) \subset M_-(I)$ for left and right chiral nets which are indeed maximal and should
be regarded as the subfactor version of left- and right maximal extensions of the chiral algebra.

If $j$ is an order-$n$ simple-current (i.e. the permutation of primaries corresponding to $j$ is 
order $n$) then $T_{j,j}/T_{1,1}$ will be a $2n$-th root of unity; when $T_{j,j}/T_{1,1}$ 
is also order-$n$, then a modular invariant $\cZ_{\langle j\rangle}$ can be associated to 
the group generated by $j$, namely
$\langle j\rangle
\simeq\bbZ_n$; these are called the \textit{simple-current invariants} (see e.g. \cite{SY}).
In this case write $T_{j,j}/T_{1,1}=\exp[2\pi\i h_j]$ for $h_{j}\in\frac{1}{n}\mathbb{Z}$. Then 
  \begin{equation}\label{scinv}
(\cZ_{\langle j\rangle})_{\lambda,\mu}= \sum_{1\le l\le n}\delta^{{\mathbb{Z}}}
\left(Q_j(\lambda)-{l\,h_{j}}\right)\, \delta_{\mu\,
,j^{l}\lambda}\ ,
\end{equation}
 where $\delta^{{\mathbb{Z}}}(x)$ equals 1 or 0 depending on whether or
not $x$ is integral. $\cZ_{\langle j\rangle}$ will be an automorphism invariant iff the root of unity
$T_{j,j}/T_{1,1}$ has order exactly equal to that of the permutation $j$.

Another source of modular invariants are the \textit{conformal embeddings},
which are VOAs $\cV'$ containing $\cV$ (or conformal nets and subnets in the subfactor picture) but with the same central charge.
The corresponding modular invariant is built from the branching rules
expressing irreducible extended modules as direct sums of the original modules.
The result is a block-diagonal modular invariant.

An important example is the Drinfeld double of  finite groups $G$. Throughout this
paper we will restrict to the simpler case of trivial 3-cocycle twist, as finite groups are
not our primary interest. 
The primaries are  pairs $(g,\chi)$ where $g$ is a representative of a conjugacy class in $G$ and 
$\chi$ is an irreducible representation (\textit{irrep}) of the centraliser $Z_g(G)$.
The vacuum 1 is $(1,1)$. We'll write Ver$(G)$ for its Verlinde ring.
 The sufferable modular invariants  \cite{Ost,Ev},
are parametrized by pairs $(H,\psi)$ where $H$ is a subgroup of
$G\times G$, and $\psi\in H^2_H(\textrm{pt};\bbT)$ here plays the role of discrete 
torsion: then $\beta_g(h):=\psi(g,h)\,\psi(h,g)^{-1}$ is a 1-dimensional representation of
$Z_H(g)$ and can be used to twist a modular invariant, for example. The diagonal modular invariant corresponds to the choice
$H=\Delta=\{(g,g) :g\in G\}$ and $\psi=1$. Sections 3.2 and 5.1 include examples of
$\psi\ne 1$. No such parametrisation of
sufferable modular invariants is known in the loop group setting, though no
SU($n$) modular invariant is known to us to be insufferable (there are plenty of insufferable 
modular invariants for other Lie groups $G$ however).

The other important example for us is the loop group $LG$ at level $k\in\bbZ_{>0}$,
where $G$ is a compact, connected, simply-connected Lie group. We write Ver$_k(G)$
for its Verlinde ring and $P_+^k(G)$ for its primaries $\Phi$. In particular, for $G$ of rank $r$,
$\lambda\in P_+^k(G)$ is the affine highest-weight $\lambda=(\lambda_0;
\lambda_1,\ldots,\lambda_r)\in\bbZ_{\ge 0}^{r+1}$ where $\sum a^\vee_i\lambda_i=k$
for positive integers $a_i^\vee$ (the \textit{co-labels}) depending only on $G$. For all $G$,
$a_0^\vee=1$ and $(k;0,\ldots,0)$ denotes the vacuum 1.
Ver$_k(G)$ can be expressed as $R_G/I_k(G)$ for some
ideal $I_k(G)$ of $R_G$ called the \textit{fusion ideal}, and primary $\lambda$
is associated to class $[\rho_{\overline{\lambda}}]\in R_G/I_k(G)$ where
$\rho_{\overline{\lambda}}$ is the $G$-irrep with highest-weight 
$\overline{\lambda}=(\lambda_1,\ldots,\lambda_r)$. For example,
for $G=\mathrm{SU}(n)$, we have all $a_i^\vee=1$,
there are $n$ simple-currents, and charge-conjugation
is nontrivial iff $n>2$.

Restrict now to the loop group setting, for $G$ a Lie group of rank $r>0$, and 
fix a modular invariant $\cZ$ 
and compatible nimrep $\cN_\lambda$. Equivalence classes of 
\textit{D-branes} (extended structures in space-time on which the end-points
of open strings reside)
are parametrised by the boundary states $x\in\cB$. The dynamics of
the branes is controlled by their \textit{conserved charges}, i.e. an assignment of an integer $q_x$
to each $x\in\cB$  and a choice of integer $M$ such that
\begin{equation}\label{charges}
\mathrm{dim}(\overline{\lambda})\, q_x\equiv\sum_{y\in\cB}\mathcal{N}_{\lambda,x}^yq_y\ (\mathrm{mod}\ M)\,\end{equation}
is satisfied for each primary $\lambda=(\lambda_0;\lambda_1,\ldots,\lambda_{r})$
 and boundary state $x$, where dim$(\overline{\lambda})$
denotes the dimension of the $G$-irrep with highest-weight $(\lambda_1,\ldots,\lambda_r)$. 
The rescaled assignment $(\{nq_x\},nM)$, for any nonzero integer $n$,  is regarded as 
equivalent to  $(\{q_x\},M)$; the set of all equivalence classes of assignments $(\{q_x\},M)$
satisfying \eqref{charges} forms a $\bbZ$-module $\mathcal{M}_\cN$ called the \textit{charge-group},
which can be regarded as the universal solution to \eqref{charges} for the nimrep $\cN$.
The unfortunate restriction here to CFTs associated to loop groups is because
the analogue of dim$(\overline{\lambda})$ is not clear for other RCFT (perhaps it involves the
\textit{special subspace} of Nahm \cite{Nahm}, as that should give an upper 
bound on the 
quantum-dimensions $S_{\lambda,1}/S_{1,1}$ --- we thank Matthias Gaberdiel for
communication on this point). In section 7 we make a definite proposal for 
$G$ a finite group.

For example, consider the Verlinde nimrep $\cN_\lambda=N_\lambda$.
If we  start with $q_1=1$, then  \eqref{charges} forces 
 $q_\lambda= \mathrm{dim}(\overline{\lambda})$ for each primary $\lambda$, and the largest possible $M$ 
which works is $\widetilde{M}=\mathrm{gcd}_{\overline{\mu}\in I_k(G)}
\mathrm{dim}(\overline{\mu})$, where $I_k(G)$ is the fusion ideal. 
More generally, any solution to \eqref{charges} will satisfy  $q_\lambda= 
q_1\mathrm{dim}(\overline{\lambda})$; this implies we can rescale it so that $M=\widetilde{M}$, and that
charge-assignment is then recovered (up to equivalence) by the value $q_1\in
\bbZ_{\widetilde{M}}$. Hence the charge-group $\cM_N$ is simply 
$\bbZ_{\widetilde{M}}$. A consequence of this discussion is that the charge-group
$\cM_\cN$ for any nimrep $\cN$ of Ver$_k(G)$ must have ${\widetilde{M}}$-torsion;
i.e. ${\widetilde{M}} \cM_N  = 0$.

Consider for concreteness $G=\mathrm{SU}(2)$ at level $k$, with $k+1$ primaries 
$\lambda=(\lambda_0;\lambda_1)$ for $\lambda_0+\lambda_1=k$ and $0\le \lambda_i\le k$.
Its Verlinde ring Ver$_k(G)$ is the quotient of $R_G$ by the fusion ideal $(\sigma_{k+2})$, where 
the class $[\sigma_i]$ corresponds to primary $(k+1-i;i-1)$. The vacuum is $(k;0)$. Charge-conjugation here
is trivial and there is one nontrivial simple-current, $j=(0;k)$, corresponding to permutation $j\,
(\lambda_0;\lambda_1)=(\lambda_1;\lambda_0)$ with $Q_j(\lambda)=\lambda_1/2$ and
$h_j=k/4$. Its modular invariants are classified in \cite{CIZ} and all are sufferable; 
they  fall into  an A-D-E pattern:

\begin{itemize}

\item[(i)] For any $k$, $\bbA_{k+1}$ corresponds to the diagonal modular invariant $\cZ=I$. The 
nimrep and full system are the Verlinde ring Ver$_k(G)={}^{k+2}K^1_G(G)$ and both
alpha-inductions are the identity. 
The $\bbA_{k+1}$ diagram describes the multiplication in Ver$_k(G)$ by the 
fundamental
weight $\sigma=\sigma_2\in R_G$ (which generates Ver$_k(G)\simeq \bbZ[\sigma]/
(\sigma_{k+2})$) on the preferred basis $[\sigma_1],[\sigma_2],\ldots,
[\sigma_{k+1}]$. The
D-brane charge-group $\cM_{\bbA}$ is $\bbZ_{k+2}$, generated by the assignment $[\sigma_l]\mapsto 
q_{[\sigma_l]}=\mathrm{dim}(\sigma_l)=l$ for $\sigma_l\in R_G$ (modulo $k+2$ 
this dimension is well-defined).

\item[(ii)] For any even $k$, the Dynkin diagram $\bbD_{k/2+2}$ likewise
describes the nimrep corresponding 
to the simple-current invariant $\cZ_{\langle j\rangle}$, again with respect
to the preferred basis. When $k/2$ is odd (resp. even), the modular invariant
is an automorphism invariant (resp. in block-diagonal form), the full system is Ver$_k(G)$ 
\cite{BEK2} (resp. two copies of $\bbD_{k/2+2}$ \cite[III]{BE1}), and
the charge-group $\cM_{\bbD}$ is $\bbZ_4$ (resp. $\bbZ_2\times\bbZ_2$).

\item[(iii)] The diagram $\bbE_6$ defines the nimrep corresponding to the conformal 
embedding SU$(2)_{10}\rightarrow \mathrm{Sp}(4)_1$. Its full system etc
is given in \cite[III]{BE1}. The charge-group is $\bbZ_3$.

\item[(iv)] The diagram $\bbE_7$ defines the nimrep corresponding to a $\bbZ_2$-orbifold of 
the $\bbD_{10}$
simple-current invariant at $k=16$. We discuss this example in section 6.2 below.
 The charge-group is $\bbZ_2$.

\item[(v)] The diagram $\bbE_8$ defines the nimrep corresponding to the conformal embedding 
SU$(2)_{28}\rightarrow G_{2,1}$. Its full system etc is given in \cite[III]{BE1}.
The charge-group is trivial.

\end{itemize}

\subsection{$K$-theory and CFT}

Consider first $G$ a finite group. This case has been an indispensible
guide to the generalisations of Freed-Hopkins-Teleman. We restrict attention 
here to trivial $H^3$-twist --- see \cite{Ev} for a more complete treatment. 
The Verlinde ring is  $K^0_{G^{\mathrm{ad}}}(G)=
K^0_{\Delta_G^L\times \Delta_G^R}(G\times G)$. This recovers naturally the 
preferred basis described in section 2.3.

Fix now a sufferable modular invariant, i.e. a subgroup $H$ of $G\times G$
and some $\psi\in H^2_H(\mathrm{pt};\bbT)$. Assume for now that $\psi$ is trivial.
The $K$-theory of the full system is $K^0_{H^L\times H^R}(G\times G)$.
The boundary states $x\in\cB$ are the natural basis of $K^0_{H^L\times \Delta^R_G}
(G\times G)=K^0_H(G)$ where $(h_1,h_2).g=h_1gh_2^{-1}$; the $x\in\cB$
therefore consist of a pair $(H(g_1,g_2)\Delta_G,\phi)$ of an $H\times \Delta_G$-orbit
in $G\times G$, and an irrep $\phi$ of the stabiliser in $H\times\Delta_G$ of 
$(g_1,g_2)$. The special element in the nimrep, corresponding in the subfactor
language to the inclusion $\iota:N\subset M$ of factors, is
$\iota=(H\Delta_G,1)$. Likewise, 
the conjugate map, the surjection $\overline{\iota}:M\rightarrow N$,
corresponds to $\overline{\iota}=(\Delta_GH,1)\in K^0_{\Delta_G^L\times
H^R}(G\times G)$. When the 2-cocycle $\psi$ is 
not trivial, the $R_{H\times H}$ and $R_{H\times G}$ module structures of the
full system  $K^0_{H^L\times H^R}(G\times G)$ and nimrep 
$K^0_{H^L\times \Delta^R_G}(G\times G)$ resp. are unchanged, but the bundles,  
alpha-induction, sigma-restriction, the product in $\MXM$, the Ver$(G)$-module structure
of $\MXN$, etc will be twisted by $\psi$. We'll see examples of nontrivial 
$\psi$ later in the paper.

We are more interested in this paper in the case of Lie groups (of dimension
$>0$), where considerably less is known. The main result we build on is the 
expression of the Verlinde ring as a $K$-group.
Let $G$ be a compact group
(not necessarily connected or simply-connected), with identity connected 
component $G_1$, and fix some element
$f\in G$. Write $[f G_1]$ for the conjugacy class in $G$ containing $fG_1$,
and $L_fG$ for the $f$-twisted loop group 
consisting of all maps $\gamma:\bbR\rightarrow G$ satisfying $\gamma(t+1)=
f\gamma(t)f^{-1}$. Fix any twist $\tau\in H^1_G([fG_1];\bbZ_2)\times
H^3_G([fG_1];\bbZ)$, and let ${}^{\tau}R(L_f G)$ denote the space of admissible
representations --- 
see \cite{FHTlg} for its definition and
its interpretation using highest-weight modules for the corresponding
$f$-twisted affine algebra. Then the main theorem of Freed-Hopkins-Teleman is:

\medskip \noindent\textbf{Theorem 1.} \cite{FHTlg} \textit{There is a natural 
isomorphism
\begin{equation}
^{\tau+\sigma} K_{G^{\mathrm{ad}}}^{d}([f G_1])\simeq {}^{\tau}R(L_f G)\,,\label{fhtthm3}\end{equation}
where $d$ is the dimension of the centraliser $Z_G(f)$, and
 the twist $\sigma$ is the cocycle
for the projective action of $L_fG$ on the graded Clifford algebra Cliff$(\mathfrak{g}^*)$.
Moreover, $^{\tau+\sigma} K_{G^{\mathrm{ad}}}^{d+1}([f G_1])\simeq 
{}^{\tau}R^1(L_f G)$, the corresponding space of graded representations.}\medskip

The relation between $K$-theory and branes is reviewed in \cite{Mokt}. 
Consider a compact, connected, simply-connected Lie group $G$ of rank $r$,
at level $k$, and write $\kappa=k+h^\vee$ where $h^\vee$ is the dual Coxeter
number. This corresponds to the diagonal modular invariant 
and Verlinde nimrep. The charge-group $\bbZ_{\widetilde{M}}$ (recall
section 2.3) is related to the twisted \textit{nonequivariant} $K$-homology $^\kappa K_*(G)$
(see Conjecture 1(a) below).  The spectral sequences argument of section 2.1
identifies $H^3_{G^{\mathrm{ad}}}(G;\bbZ)\simeq\bbZ$ with $H^3(G;\bbZ)$.
These groups  $^\kappa K_*(G)$ can be computed (via Poincar\'e duality) from 
$^\kappa K^*_{G^{\mathrm{ad}}}(G)$ using the Hodgkin spectral sequence 
(recall section 2.2) with $H=1$, together with the resolution of Ver$_k(G)$ found in e.g. \cite{M}; the result
\cite{Bra,Doug} can be expressed as $^\kappa K_0(G)=\oplus_{l\ \mathrm{even}}\mathrm{Tor}_l$,
$^\kappa K_1(G)=\oplus_{l\ \mathrm{odd}}\mathrm{Tor}_l$, where Tor$_l$ is identified with
the degree-$l$ part of the exterior algebra $\bigwedge\{x_1,\ldots,x_{r-1}\}$.
More explicitly, $^\kappa K_i(G)\simeq
\bbZ_M^{2^{r-2}}$ for all $i$ (except for $G=\mathrm{SU}(2)$,
where $^\kappa K_i(\mathrm{SU}(2))\simeq \bbZ_M,0$ for $i$ even, odd respectively), 
for some integer $M=M_k(G)$ depending on $G$ and $\kappa$. For example, 
 $M_k(\mathrm{SU}(n))=\frac{\kappa}{\mathrm{gcd}(\kappa,y)}$ for  $y=\mathrm{lcm}(1,2,\ldots,n-1)$.

\medskip\noindent\textbf{Conjecture 1(a)} \textit{For all simply-connected, compact, connected
$G$ and any level $k$, the integer $M_k(G)$ appearing in the groups $^\kappa K_*(G)$
equals the gcd of the dimensions of all weights $\nu$ in the fusion ideal $I_k(G)$, or equivalently
 the charge-group $\cM_N$ for the Verlinde nimrep $\cN_\lambda=
N_\lambda$ is isomorphic to $\bbZ_{M_k(G)}$}.

\noindent\textbf{(b)} \textit{The largest number $M$ (call it $M_k^{spin}$) such that 
\begin{equation}
\mathrm{dim}(\lambda)\,\mathrm{dim}(\mu)\equiv \sum_{\nu\in\Phi}N_{\lambda,\mu}^\nu
\mathrm{dim}(\nu)\ \ (\mathrm{mod}\ M)\end{equation}
holds for all 
$\mathrm{SO}(2r+1)$ non-spinor dominant weights $\lambda$ and all level $k$ 
$\mathrm{SO}(2r+1)$-spinors $\mu$
 (i.e. $\mu_r$ is odd and $\lambda_r$ is even), is $2^r$ times the integer 
 $M_k(\mathrm{Spin}(2r+1))$ appearing in Conjecture 1(a).}\medskip

It is clear the $K$-theory $M_k(G)$ \textit{divides}  the gcd in Conjecture 1(a). The cases $G=
\mathrm{SU}(n)$ and $G=\mathrm{Sp}(2n)$ of Conjecture 1(a) have been proved in \cite{BDR}; relatively
simple expressions for the $K$-theoretic $M$ are found in \cite{Doug2}, and extremely
simple formulas for $M$ are conjectured in e.g. \cite{BDR}. Conjecture 1(b) is needed in section
4 below; it was proved in \cite{GaGa1} whenever 4 doesn't divide $M_k(\mathrm{SO}(2r+1))$,
where it was also shown $M_k^{spin}=2^iM_k(\mathrm{Spin}(2r+1))$ for some integer $0\le i\le r$.
There is considerable  numerical evidence supporting Conjecture 1.

The assignment $x\mapsto q_x$ of charges here is given by the map $^\kappa K_*^G(G)
\rightarrow {}^\kappa K_*(G)$ forgetting the $G$-equivariance. This can be calculated
through the commuting diagram:
\begin{equation}\label{Dch}
\begin{array}{ccccccc}
0\longrightarrow &I_k(G) & {\longrightarrow} &
R_G=K_0^G(1) & \stackrel{\alpha\,}{\longrightarrow} &
^{\tau}K^{G}_{0}(G)&\longrightarrow 0\\[3pt]
&\beta\downarrow & &\beta\downarrow & & \gamma\downarrow&\\[2pt]
0\longrightarrow&M\bbZ&\longrightarrow & \bbZ=K_0(1)&\stackrel{\delta\,}{\longrightarrow} &
^\kappa K_0(G)&
\end{array}.
\end{equation}
$I_k(G)\subset R_G$ is the fusion ideal, and Theorem 1 gives the
top line, where the ring homomorphism $\alpha$ is the push-forward of the embedding of
the identity 1 in $G$. The map $\beta:K_0^G(1)\rightarrow K_0(1)$, forgetting $G$-equivariance,
is of course given by dimension. Conjecture 1(a) says $\beta$ carries $I_k(G)$ 
surjectively onto $M\bbZ$. To identify the target of the map $\delta$, compare
 the Hodgkin spectral sequence for $^\kappa K_0(G)$, i.e. $E^\infty_{p,0}=
E^2_{p,0}=\mathrm{Tor}_p^{R_G}(\bbZ,{}^\kappa K_{0}^G(G))$  (and all other 
$E^\infty_{p,q}=0$), with that for $K_0(1)$, i.e. 
$E^\infty_{0,0}=E^2_{0,0}=\mathrm{Tor}_0^{R_G}(\bbZ,K_0^G(1))$ (and all other $E^\infty_{p,q}
=0$): we see that the target of $\delta$ should be  $\mathrm{Tor}_0^{R_G}(\bbZ,{}^\kappa K_{0}^G(G))=\bbZ\otimes_{R_G}\mathrm{Ver}_k(G)$, with $\delta[\rho]=
1\otimes_{R_G}\rho$.

\section{Group-like fusion rings}

{ This section considers the easiest possibility, when all primaries are simple-currents
 (recall the definition in section 2.3). We begin though  with section 3.1, which describes the modular invariants of arbitrary finite groups.}

\subsection{Finite groups}

As explained in section 2.3, the sufferable modular invariants for the 
finite group $G$ setting are parametrised by pairs $(H,\psi)$ where $H$
is a subgroup of $G\times G$ and $\psi\in H^2_H(\mathrm{pt};\bbT)$. 
Here we  explain, using the methods of  \cite{Ev} and expanded in \cite{EG1, EG2},   how
to arrive at a modular invariant $\cal {Z}$ when $\psi$ is trivial (and the twist
$\tau$ on the double $\cD(G)$ is also trivial). This is useful in sections
3.2, 4.1, 5.1. In sections 3.2 and 4.1 we also consider $H^2$-twists on the subgroups $H$.

Let $G$ be a finite group, and $H$ a subgroup of $G  \times G$. Let $\sigma : G \times G \rightarrow 
G \times G$ be the flip $(a,b) \mapsto (b,a)$, and $\pi_{\pm} : G \times G 
\rightarrow G$ be the coordinate maps $(a,b) \mapsto a$ and $b$ respectively. 
Then $K_{\pm} = \pi_{\pm}(H)$ are subgroups of $G$. Let
\begin{align}
N_{+} = &\,\mathrm{ker}\,{ \pi_{-}|}_{H} = H \cap (K_{+}  \times 1) = H \cap (G \times 1)\,,\nonumber\\
N_{-} =&\, \mathrm{ker}\, {\pi_{+}|}_{H} = H \cap (1 \times K_{-} ) = H \cap (1 \times G )\,,
\end{align}
so that $N_{\pm} \triangleleft H$ and $N_{\pm} \triangleleft K_{\pm}$. Then $H/N_{\pm} \simeq 
K_{\mp}$ and 
\begin{equation}
K_{+}/N_{+} \simeq H / N_{-}/N_{+} \simeq H / N_{+}/N_{-} \simeq K_{-}/N_{-} \,,\nonumber
\end{equation}
because $N_+$ and $N_{-}$ pairwise commute in $G \times G$. Let $H^{\pm} = 
\Delta_{K_{\pm}}(N_{\pm} \times 1)$ which are $\sigma$-invariant subgroups of 
$G \times G$.  The subgroups $H^+, H^-$ will have extended systems given by 
the doubles of $K_+/N_+$ and $K_-/N_-$, respectively. Write
$b^{\pm}$ for their branching coefficients (sigma-restrictions) from
$\cD(K_\pm/N_\pm)$ to $\cD(G)$. Then the modular invariant for the pair
$(H,0)$ is obtained as 
\begin{equation}
 \cZ_{\lambda \mu} = \Sigma_{\tau} b^+_{\tau \lambda} b^-_{\beta (\tau) \lambda}\nonumber
 \end{equation}
where $\tau$ runs over the primary fields of $K_+/N_+$, and $\beta$ is the identification of the chiral primary fields via the above isomorphism
$K_{+}/N_{+}  \simeq K_{-}/N_{-} \,.$

It therefore suffices to determine sigma-restriction for the flip invariant 
case $H = \sigma H$ where $K_+=K_-$ and $N_+=N_-$. Let $K = K_{\pm}\,,$
$N = N_{\pm}$ so that $H = \Delta_K(N \times 1)\,;$ then $N \triangleleft K$, and $K$ can be an arbitrary subgroup of $G$.

Consider first the case $\Delta_G \subset H \subset G \times G$ so that $K = G$, and $N \triangleleft G$. Then sigma-restriction 
$ K^0_{G/N}(G/N) \rightarrow K^0_{G}(G)\,,$  is described as follows 
\cite{EG1}. A primary field in the $G/N$ theory is labelled $[gN, \chi]$ for 
a conjugacy class in $G/N$ of a coset $gN$, and $\chi$ a representation of the 
centralizer ${Z}_{gN}(G/N)$. The quotient map $\pi : G \rightarrow G/N$ takes 
${Z}_{gn}(G)$ to ${Z}_{gN}(G/N)$ for any $n \in N$, and $\chi \pi$ is a representation of 
${Z}_{gn}(G)$. The preimage $\pi^{-1}$ of the conjugacy class of $gN$ will be a disjoint union
of conjugacy classes of $G$, each with a representative of the form $gn$ for $n\in N$ (since
$N$ is normal).
Then the sigma-restriction is $\sigma [gN , \chi] = \Sigma_n [gn, \chi \pi]\,$, where the sum
is over the $n\in N$ giving  a disjoint union.

Next consider the case $H \subset \Delta_G \subset G \times G$ so that $N = 1$ and $H = 
\Delta_K$. Here the extended system is $K$ on $K$ and we need sigma-restriction 
$ K^0_{K}(K) \rightarrow K^0_{G}(G)\,.$  Take a primary field $[k, \chi]$ where
$k$ describes a conjugacy class in $K$ and $\chi$ a representation of the centraliser 
${Z}_{k}(K)\,.$ Then sigma-restriction is given by  \cite{EG2} : 
$[k , \chi]  \mapsto [k, \mathrm{Ind}_{{Z}_{k}(K)}^{{Z}_{k}(G)} \chi ]\,.$

The general case $\Delta_K(N \times 1) \subset G \times G$ is a gluing of the  special cases 
$\Delta_K \subset \Delta_K (N \times 1) \subset K \times K$ and $\Delta_K \subset \Delta_G 
\subset G \times G$. Consequently we have  sigma-restriction in stages:
\begin{align} K^0_{K/N}(K/N) &\,\rightarrow K^0_{K}(K) \rightarrow K^0_{G}(G)\,,\nonumber\\ 
[kN , \chi] &\, \mapsto  \Sigma_n [kn, \mathrm{Ind}_{{Z}_{kn}(K)}^{{Z}_{kn}(G)} (\chi \pi) ]\,,
\nonumber\end{align}
where $\pi : K \rightarrow K/N$ is the quotient map, $kN$ represents a 
conjugacy class in $K/N$, and $\chi$ is a representation of the stabiliser 
${Z}_{kN}(K/N)\,.$ The sum over $n$ is as before.

\subsection{Finite abelian groups}

{ Now consider finite groups $G$  where all primaries are
simple-currents}. 
A primary $(g,\chi)$ of a finite group $G$ is a {simple-current}
iff $g$ lies in the centre of $G$ and $\chi$ is 1-dimensional. This means
all primaries will be simple-currents,
 iff $G$ is abelian. 
 
{ Assume} for concreteness for now that $G$ is the cyclic group
$\bbZ_p$ for $p$ prime.
For these $G$, $H^3_{G^{\mathrm{ad}}}(G;\bbZ)= 0$ while $H^4_G({\rm
pt};\bbZ)\simeq \bbZ_p$, so the transgression $\tau$
is trivial. As an $R_G$-module,
$^{\tau(\sigma)}K_G^0(G)\simeq R_G^{p}$ is independent
of $\sigma$, but as a ring, $^{\tau(\sigma)}K^0_G(G)$ is
the group ring $\bbZ[\bbZ_d\times\bbZ_{p^2/d}]$ where
$d={\rm gcd}(2p,\sigma)$. Fix trivial twist $\sigma=p$ for concreteness,
so we can identify the primaries with pairs $(a,b)\in\bbZ_p^2$, with the 
group-like fusion product $(a,b)(a',b')=(a+a',b+b')$. 

The modular data is
$S_{(a,b),(a',b')}=\frac{1}{p}\exp(2\pi\i\,(ab'+a'b)/p)$ and $T_{(a,b),(a,b)}=\exp(2\pi\i ab/p)$.
There are precisely $2p+2$ modular invariants: $2p-2$ of these are 
automorphism invariants,
defined for any $\ell\in\bbZ_p,\ell\ne 0$,  by
$\cZ^{(\ell)}_{(a,b),(\ell a,\ell^{-1}b)}=1=\cZ'{}^{(\ell)}_{(a,b),(\ell b,\ell^{-1}a)}$, and 
all other entries vanish. 
The remaining modular invariants are $\cZ^{m,n}$ for $m,n\in\{0,1\}$,
defined by $\cZ^{m,n}_{(ap^m,bp^n),(a'p^n,b'p^m)}=1$ for all $a,b,a',b'\in\bbZ_p$ 
(all other entries are 0). $\cZ'{}^{(\ell)},\cZ^{0,0},\cZ^{1,1}$ are
simple-current invariants $\cZ_{\langle(\ell,1)\rangle},\cZ_{\langle(1,0)\rangle},
\cZ_{\langle(0,1)\rangle}$ respectively.

All of these modular invariants are sufferable. Explicitly, the subgroups $H$ of $\bbZ_p\times\bbZ_p$ 
are $0\times 0$, $\langle(1,b)\rangle$ for any $b\in\bbZ_p$, $0\times\bbZ_p$,
and $\bbZ_p\times\bbZ_p$. For $H\simeq \bbZ_{d}\times\bbZ_{d'}$,
the K\"unneth formula 
for group cohomology \cite{Brown} says $H^2_H(\mathrm{pt};\bbT)\simeq {\bbZ_{d}}\otimes_\bbZ {\bbZ_{d'}}
\simeq\bbZ_{\mathrm{gcd}(d,d')}$. This means the choice $H=\bbZ_p\times\bbZ_p$
has a twist of $\psi\in\bbZ_p$ but the other $2+p$ choices for $H$ all come with trivial $\psi$.
The correspondence between $\cZ$ and $(H,\psi)$ is given in Table 1.

$$\mbox{\scriptsize$\begin{array}{|c||c|c|c|c|c|c|c|} \hline
&{(H,\psi)} &K^0_{H\times G}(G\times G)& \rm{neutral} & K^0_{H\times H}(G\times
G)&E&D_+&D_- \\ \hline
\cZ^{(\ell)}& (\langle (1,\ell)\rangle,0)&{\rm{Ver}}(G),\bbZ\ \rm{resp}&
{\rm{Ver}}(G)&{\rm{Ver}}(G)&G\times G,0\times 0\ \rm{resp}&0\times 0&0\times 0\\ 
&&\rm{for}\ \ell=1,\ne 1&&&\rm{for}\ \ell=1,\ne 1&&\\ \hline
\cZ'{}^{(\ell)}& (G\times G,\ell)&R_G&
{\rm{Ver}}(G)&R_{G\times G}&\langle(\ell,1)\rangle&0\times 0&0\times 0\\ \hline
\cZ^{0,0}& (G\times G,0)&R_G&\bbZ&R_{G\times G}&G\times 0&G\times 0&G\times 0\\ \hline
\cZ^{1,1}& (0\times 0,0)&\bbZ(G)&
\bbZ&\bbZ(G)\otimes\bbZ(G)&0\times G&0\times G&0\times G\\ \hline
\cZ^{0,1}& (G\times 0,0)&\bbZ&
\bbZ&R_G\otimes \bbZ(G)&0\times 0&G\times 0&0\times G\\ \hline
\cZ^{1,0}& (0\times G,0)&\bbZ&
\bbZ&\bbZ(G)\otimes R_G&0\times 0&0\times G&G\times 0\\ \hline
\end{array}$}$$
\centerline{\textbf{Table 1.} Data for the $G=\bbZ_p$ modular invariants}\medskip
 
The nimrep $K^*_{H^L\times \Delta_G^R}(G\times G)=K^*_{H^{\mathrm{ad}}}(G)$ 
is recovered as follows. In all cases $K^1_H(G)=0$; the $K$-groups
$K^0_H(G)$ are collected in Table 1.
Note that in all cases $K^0_{H\times G}(G\times G)$ is isomorphic as an
additive group to the group ring $\bbZ(E)$ where the subgroup $E$ of $G\times G$
is defined to be $\{(a,b)\in\bbZ_p^2:
\cZ_{(a,b),(a,b)}=1\}$.  The module structure can be written down as follows.
Using the nondegenerate pairing $(a,b)\cdot(a',b')=ab'+a'b\in\bbZ_p$ introduced
earlier, for any subgroup $E$ of $G\times G$ define 
$E^*=\{(a,b)\in\bbZ_p^2:(a,b)\cdot E=\{0\}\}$.
Then the boundary states of the nimrep are $\cB=(G\times G)/E^*$, which is
isomorphic to $E$ as a group, and the primary $(a,b)$
acts on $[x]\in \cB$ by addition mod $E^*$.

The full system $K^0_{H^L\times H^R}(G\times G)$ is collected in Table 1. To identify 
alpha-induction, it is more convenient to write the full system in equivalent
form as 
$(G\times G)/D_-\times (G\times G)/D_+^*$, where $D_\pm$ are given in Table
1 and are defined by $D_+=\{(a,b):\cZ_{(a,b),(0,0)}=1\}$, $D_-=\{(a,b):\cZ_{(0,0),(a,b)}=1\}$.
Then $D_\pm$ are isomorphic as groups, and both $D_\pm^*/D_\pm$ are isomorphic
to the Verlinde ring of the neutral system, also collected in Table 1 ($D_\pm\le D_\pm^*$
by $T$-invariance).
Alpha-induction is given by $\alpha^+_{(a,b)}=((ma,m^{-1}b)+D_-,(a,b)+D_+^*)$ 
and $\alpha^-_{(a,b)}=((a,b)+D_-,D_+^*)$ for all $\cZ$ (except $\alpha^-_{(a,b)}=
((b,a)+D_-,D_+^*)$ for $\cZ^{\prime(\ell)}$), where $m=\ell$ for
$\cZ^{(\ell)}$, $m=\ell^{-1}$ for $\cZ'{}^{(\ell)}$, and $m=1$ for the four other $\cZ$.

More generally, for $G=\bbZ_{p^\nu}$, the bijection of Table 1 between modular
invariants $\cZ$ and pairs $(H,\psi)$ is lost. 
For $G=\bbZ_{p^\nu}$, there are precisely $((\nu+1)p^{\nu+2}-
2p^{\nu+1}-(\nu+1)p^\nu+2)/(p-1)^2$ pairs $(H,\psi)$, each of which is either of the
form $H=\langle (\frac{n}{d},\ell\frac{n}{d}),(0,\frac{n}{d'})\rangle$ and $\psi\in\bbZ_{d'}$
for any $d'|d|n$, $0\le \ell<\frac{d}{d'}$, or of the form $H=\langle (\ell \frac{pn}{d},\frac{n}{d}),
(\frac{n}{d'},0)\rangle$ and $\psi\in \bbZ_{d'}$ for any $d'|d|n$, $d'<d$, $0\le \ell<\frac{d}{d'p}$. 
$G=\bbZ_{p^\nu}$ has precisely $2p^\nu-2p^{\nu-1}$ automorphism
invariants: namely $\cZ^{(\ell)}$ and $\cZ^{\prime(\ell)}$ for each $\ell\in
\bbZ^\times_{p^\nu}$. These correspond bijectively with the pairs $(\langle
(1,\ell)\rangle,0)$ and $(G\times G,\ell)$ as above.
 There are precisely $(\nu+1)^2$ $\cZ^{m,n}$, defined
in the obvious way. They are realised by pairs $(\langle(p^m,0),(0,p^n)\rangle,
0)$ as above, but for $\nu>1$ these can be realised by other pairs. For example
for $G=\bbZ_{p^2}$, $\cZ^{1,1}$ and $\cZ^{0,0}$ are realised by any
$(\langle (p,0),(0,p)\rangle,\bbZ_p)$ and $(G\times G,p\bbZ_p)$ respectively.
All other of the remaining $(p-1)(p+3)$ modular invariants for $G=\bbZ_{p^2}$
are also sufferable, although $2(p-1)$ of these are realised by exactly 2
pairs. At the time of writing, we do not know if all modular invariants
for cyclic $G$ are sufferable.

\subsection{Compact Lie groups}

{ In this subsection we investigate the case where $G$ is a  torus. These are
the only connected Lie groups whose primaries are all simple-currents, for all levels.
This case is simple enough that we can give the complete story, which is summarised in
 Theorem 2 below --- the first main result of the paper.}

Let $L\subset\bbR^n$ be any $n$-dimensional lattice, and write $T_L$ for the
tori $\bbR^n/L$. Of course all $n$-tori are homeomorphic, say to $T_L$,
and as explained
in \cite[Sect.2.2]{EG1} $\tau\in H^3_{T_L^{\mathrm{ad}}}(T_L;\bbZ)$-twist here 
can be taken to be $\tau\in\textrm{Hom}(L,L^*)$, where $L^*=
\{x\in\bbR^n:x\cdot L\subseteq\bbZ\}$ is the dual lattice.
In this section we will find it convenient to identify these twists $\tau$ with their image, i.e. with 
$n$-dimensional
sublattices $M$ of $L^*$, so in place of ${}^\tau K^*_{T_L}(T_L)$ we will write $K_{T_L}^*(T_M)$.
This $K$-group was computed in  \cite[Thm. 3.4.3]{FHTii}: 
\begin{equation}
{}^\tau K^i_{T_{L_1}}(T_{L_2})=:
K_{T_{L_1}}^i(T_{\tau(L_2)})\simeq\left\{\begin{matrix}\bbZ(L_1^*/\tau(L_2))&\mathrm{for}\ i\equiv \mathrm{dim}(T)\ (\mathrm{mod}\ 2)\\
0&\mathrm{otherwise}\end{matrix}\right.\,.\end{equation} 
All $K$-groups in this section can be  read off from that identity. 

We can and will choose $L$ so that the Verlinde ring is the group ring of
$L^*/L$, i.e. is given by $K_{T_L}^{\mathrm{dim}\,T}(T_L)$. 
The fusion product is written $[u][v]=[u+v]$;  note that each
primary, i.e. each coset in $L^*/L$, is a simple-current. This torus example
formally behaves similarly to the case of finite abelian groups considered
last subsection, in fact the torus is simple enough that we can work its
story out completely.

Consider for simplicity that $L$ is an \textit{even} lattice, i.e. that
$u\cdot u\in 2\bbZ$ for all $u\in L$. Then the matrices $S_{[u],[v]}={|L^*/L|}^{-1/2}
\exp(2\pi \i u\cdot v)$ and $T_{[u],[v]}=\delta_{[u],[v]}\exp(\pi \i u\cdot u-
\pi \i n/12)$ are well-defined, and generate the desired 
SL$(2,\bbZ)$-representation. The vacuum is $[0]$ and charge-conjugation is
$C[u]=[-u]$. { We begin with Proposition 1, which classifies all modular invariants
and nimreps for the torus without  checking for compatibility.}

\medskip\noindent\textbf{Proposition 1(a)} \textit{Modular invariants for this
chiral data are parametrized by even lattices 
$D_\pm$, such that $L\subseteq D_\pm\subseteq L^*$, and
an orthogonal  isomorphism $\beta:D_+^*/D_+\rightarrow D_-^*/D_-$
(defined below). The modular invariant is given by} 
\begin{equation}
\cZ_{[u],[v]}=\left\{\begin{matrix} 1&\mathrm{if}\ u\in D_+^*\ \mathrm{and}\
[v]\subset\beta([u]_+)\cr
0&\mathrm{otherwise}\end{matrix}\right.\,,\end{equation}
\textit{where we write $[u]_\pm$ for the $D_\pm$-coset $u+D_\pm$.}

\noindent\textbf{(b)} \textit{The nimreps (up to equivalence) are parametrised by lattices $E$, 
$L\subseteq E\subseteq L^*$: the boundary states $x\in\cB$ 
are the cosets $L^*/E^*$ and the module structure is given by $[u].[m]_{E*}=[u+m]_{E*}$ for
all $[u]\in L^*/L,[m]_{E*}\in L^*/E^*$. This nimrep is  diagonalised by
the matrix $\Psi_{[e],[m]_{E*}}={|E/L|}^{-1/2}\exp(2\pi \i e\cdot m)$
for $[e]\in E/L$, $[m]_{E*}\in L^*/E^*$, i.e. the exponents are the $L$-cosets 
$E/L$, and all multiplicities are 1.}

\medskip By \textit{orthogonal  isomorphism}
we mean a group isomorphism $\beta$ such that $[u]_+\cdot[u]_+\equiv
\beta([u]_+)\cdot\beta([u]_+)$ (mod 2), and hence $[u]_+\cdot [v]_+\equiv \beta([u]_+)\cdot
\beta([v]_+)$ (mod 1), for all $[u]_{+},[v]_{+}\subset D_+^*$.
Note that block-diagonal $\cZ$ correspond to $D_+=D_-$ and $\beta=id$,
 while {automorphism invariants} correspond to $D_+=D_-=L$. 
In particular, the diagonal modular invariant $\cZ=I$ is $D_\pm=L$ and $\beta=+1$,
while charge-conjugation $\cZ=C$ is $D_\pm=L$ with $\beta=-1$.
There is a simpler description in the $n=1$ case, more or less given by \cite{BE6}, but
we don't know a similar
description in higher dimension. Since the groups $E/L$ and $L^*/E^*$ 
are isomorphic (Pontryagin duality), the matrix $\Psi$ is square, as it must be.

\medskip
\noindent\textit{Proof.} Let $\cZ$ be any modular invariant. The triangle
inequality and $\cZ=S\cZ S^*$ imply
\begin{equation}|\cZ_{[u],[v]}|=\frac{1}{|{L^*/L}|}\,\left|\sum e^{2\pi\i u\cdot a}
\cZ_{[a],[b]}e^{-2\pi\i b\cdot v}
\right|\le \frac{1}{|{L^*/L}|}\sum \cZ_{[a],[b]}=\cZ_{[0],[0]}=1\,,\nonumber
\end{equation}
where the sums are over $[a],[b]\in L^*/L$, so each entry $\cZ_{[u],[v]}$ either equals 0 or 1. It also says that $\cZ_{[u],[v]}=1$ iff
$u\cdot a\equiv v\cdot b$ (mod 1) whenever $\cZ_{[a],[b]}=1$. Hence we get additivity:
if both $\cZ_{[u],[v]}=\cZ_{[u'],[v']}=1$, then $\cZ_{[u+u'],[v+v']}=1$. 

Define $\mathcal{L}=\cup([u];[v])\subset \bbR^{n,n}$, where the union is over all cosets
such that $\cZ_{[u],[v]}=1$, and where we place on it the indefinite inner product $(u;v)\cdot
(u';v')=u\cdot u'-v\cdot v'$. Then $\mathcal{L}$ is an even (by $T$-invariance)
lattice (by additivity). To prove self-duality, suppose $(u_+;u_-)\in\mathcal{L}^*$, i.e. $u_\pm
\in L^*$ satisfy $u_+\cdot u\equiv u_-\cdot v$ (mod 1) whenever $\cZ_{[u],[v]}\ne 0$. Then
from $\cZ=S\cZ S^*$ we get 
\begin{equation}\label{SMtorus}
\cZ_{[u_+],[u_-]}=\frac{1}{|L^*/L|}\sum e^{2\pi\i u_+\cdot u}\cZ_{[u],[v]}
e^{-2\pi\i u_-\cdot v}=\frac{1}{|L^*/L|}\sum\cZ_{[u],[v]}>0\,,\end{equation}
where the sums are over $[u],[v]\in L^*/L$, and thus $(u_+;u_-)\in\cL$. Conversely,
any even self-dual lattice $\mathcal{L}\subset\bbR^{n,n}$ containing $(L;L)$ 
defines a modular invariant: commutation with $T$ is immediate; self-duality requires
$\sum_{[u],[v]}\cZ_{[u],[v]}=|L^*/L|$ so \eqref{SMtorus} verifies $\cZ=S\cZ S^*$ at any
$(u_+;u_-)\in\cL=\cL^*$, while if $(u_+;u_-)\not\in\cL$ then there will be some $(u;v)\in
\cL^*=\cL$ for which $(u_+;u_-)\cdot(u;v)\not\in\cZ$ and so $(S\cZ S^*)_{[u_+],[u_-]}=0$
by the usual projection calculation.

Define $D_+=\{[d]\,:\,\cZ_{[d],[0]}=1\}$ and $D_-=\{[d]\,:\,\cZ_{[0],[d]}=1\}$. Then these are
even lattices satisfying $L\subseteq D_{\pm}\subseteq L^*$. Because
 $\cL$ is integral, $\cZ_{[u],[v]}=1$ implies both $([u];[v])\cdot
(D_+;0)\subseteq\bbZ$ and $([u];[v])\cdot(0;D_-)\subseteq\bbZ$, i.e. both
$u\in D_+^*$ and $v\in D_-^*$. Suppose without loss
of generality $|D_+/L|\le |D_-/L|$. Choose any $[u]\in D^*_+/L$ and compute
\begin{align}\sum_{[v]\in D_-^*/L}\cZ_{[u],[v]}&=&\frac{1}{|L^*/L|} \sum_{[a],[b]\in L^*/L} e^{2\pi\i u\cdot a}
\cZ_{[a],[b]}\sum_{[v]\in D^*_-/L}e^{-2\pi\i b\cdot v}\nonumber\\&=&\frac{|D_-^*/L|}{|L^*/L|}
\sum_{[a]\in L^*/L,[d]\in D_-/L} e^{2\pi\i u\cdot a}  \cZ_{[a],[d]}=|D_+/L|\,,\end{align}
using additivity and $|D_-^*/L|=|L^*/D_-|$. But the left-side
$\sum_{[v]} \cZ_{[u],[v]}$ will be a multiple of $|D_-/L|$, by additivity and
the fact that  $(0;D_-)$ is a sublattice of
$\mathcal{L}$. Therefore $|D_+/L|=|D_-/L|$ and there is a map $\beta:D_+^*/D_+
\rightarrow D_-^*/D_-$ such that $\beta([d]_+)$ is the unique coset for which $([d]_+;
\beta([d]_+))\subset\mathcal{L}$. In fact $\beta$ is a bijection: if $\beta([d]_+)=\beta([d']_+)$
then $(d-d';0)\in\mathcal{L}$ and hence $[d]_+=[d']_+$. The rest is clear. 

Fix a nimrep,  with boundary states $\cB=\{x_1,\ldots,x_N\}$. Then each $[u]\in L^*/L$ permutes the $x_i$, since it is a simple-current. 
Note that the indecomposability of the nimrep, required in Definition 2, is equivalent
 to the transitivity of this permutation representation.
 
 $[e]\in L^*/L$ is an exponent of the nimrep iff there exists a vector $v_{[e]}\ne 0$ (namely the
 corresponding column of $\Psi$ \eqref{nimver}), with components labelled by $x_i$, such that
 $(v_{[e]})_{[u].x_i}=e^{2\pi\i e\cdot u}(v_{[e]})_{x_i}$ for all $[u]\in L^*/L$, 
$x_i\in\cB$. Transitivity implies all components
 of $v_{[e]}$ will be nonzero. If $[e],[e']$ are exponents, then so must be $[e+e']$, whose vector
 $v_{[e+e']}$ has components $(v_{[e+e']})_{x_i}=(v_{[e]})_{x_i}(v_{[e']})_{x_i}$. Let
 $E=\cup[e]$ be the union of all exponents; then that additivity means that
$E$ is a lattice. Note that the kernel of the
 nimrep, i.e. the union of all cosets $[u]\in L^*/L$ which act trivially on all $x_i$, is precisely the 
dual lattice $E^*$.
 
Choose any $u\in L^*\setminus E^*$. Then there exists an exponent $[e']\in E/L$ such that $u\cdot e'\not\in
\bbZ$. Say $[e']$ has order $k$ in $L^*/L$. Then the permutation matrix $\mathcal{N}_{[u]}$
corresponding to multiplication by $[u]$ has trace 
$$\mathrm{tr}\,\mathcal{N}_{[u]}=\sum_{[  e]\in E/L}e^{2\pi\i e\cdot u}=\frac{1}{k}\sum_{[e]\in
E/L}\sum_{j=0}^{k-1}e^{2\pi\i\,(e+je')\cdot u}=0\,.$$ 
 Therefore any such $[u]$ acts with \textit{no} fixed points. This means the
full nimrep can be recovered from say the first column of every $\cN_{[u]}$.
 
For any $[u]\in L^*/L$, $[u].x_1$ depends only on the class $[u]+E^*$ (and 
$x_1$), so define $[u].x_1=x_{[u]+E*}$. Then this assignment is a bijection between
these cosets $L^*/E^*$ and the $x_i\in\cB$ (surjectivity is transitivity, and 
injectivity is fixed-point freeness),
satisfying $[u].x_{[u']+E*}=x_{[u+u']+E*}$. The rest is clear. 
 QED\medskip

{ Proposition 1 does all the hard work; it is now straightforward to sort out the
full system etc.}

\medskip
\noindent\textbf{Theorem 2.} \textit{Choose any modular invariant ${\cZ}_{[u],[v]}$, given by lattices $D_{\pm}$
and map $\beta:D_+^*/D_+\rightarrow D_-^*/D_-$. Put $E=\cup_{[u]_+\in D_+^*/D_+}\left(
[u]_+\cap\beta([u]_+)\right)$.}

\smallskip\noindent\textbf{(a)} \textit{$E$ is a lattice, and the nimrep
${ \bbZ(}L^*/E^*)=K^{\mathrm{dim}\,T}_{T_{L}}(T_{E^*})$ is the only one compatible with the modular invariant. Its $\mathrm{Ver}(T_L)$-module 
product $[u].[x]_{E*}=[u+x]_{E*}$ is the push-forward of the product in the Verlinde ring $K_{T_L}^{\mathrm{dim}\,T}(T_{L})$.}

\smallskip\noindent\textbf{(b)} \textit{The full system is ${ \bbZ(}L^*/D_-\times L^*/D^*_+)=K^0_{T_L\times T_L}(T_{D_-}\times
T_{D_+^*})$, and alpha-induction  $\alpha^\pm:{ \bbZ(}L^*/L)\rightarrow {
\bbZ(}L^*/D_-\times L^*/D^*_+)$ is given 
by $\alpha^+([v])=(\beta[v]_-,[v]_{+*})$ and $\alpha^-([v])=([v]_-,[0]_{+*})$, using
obvious notation.   The modular invariant is recovered by ${\cZ}_{[u],[v]}=\langle\alpha^+_{[u]},
\alpha^-_{[v]}\rangle$. The neutral system $\MXMo$ is ${\bbZ(}D^*_-/D_-)=K^{\mathrm{dim}\,T}_{T_{D_-}}(T_{D_-})$. 
Sigma-restriction $\MXMo\rightarrow \NXN$ is the unravelling $K^{\mathrm{dim}\,T}_{T_{D-}}(T_{D_-})\rightarrow
K^{\mathrm{dim}\,T}_{T_L}(T_{D_-})\rightarrow K^{\mathrm{dim}\,T}_{T_L}(T_L)$ of $D_-$-cosets into $L$-cosets. 
The subfactor inclusion $\iota$, its conjugate $\bar{\iota}$, and the canonical endomorphism 
$\theta$ are respectively $[0]_{E*}\in K^{\mathrm{dim}\,T}_{T_{L}}(T_{E^*})$,
$[0]\in K^{\mathrm{dim}\,T}_{T_E}(T_L)$, and $[0]_-=\oplus_{[u]\in D_-}[u]$.}

\smallskip\noindent\textbf{(c)} \textit{$K^i(T_{E^*})=\bbZ^{2^{n-1}}$ for each $i$, where $n=\mathrm{dim}
(T)$. The charge-group $\cM$ is $\bbZ$, generated by $q_{[x]_{E*}}=1$. The map $K^{T_L}_0
(T_{E^*})\rightarrow K_0(T_{E^*})$ forgetting $T$-equivariance is $1\otimes_{R_T}
K_0^{T_L}(T_{E^*})\subseteq \bbZ\otimes_{R_T}K_0^{T_L}(T_{E^*})$, as in section 2.4.}\medskip

The proof is immediate from the proposition.
Because of this uniform  $K$-theoretic description, we would expect all these to be
\textit{sufferable}, i.e. realised by subfactors. Note the similarity of
Theorem 2 to Table 1 of last subsection, describing the $G=\bbZ_p$ situation.
We see that the full system consists of  $|L^*/D_+^*|$
copies of the  nimrep $\bbZ[L^*/D_-]$ for the block-diagonal modular invariant 
corresponding to the neutral system. In these block-diagonal cases, the neutral 
system $K_{T_{D_-}}(T_{D_-})={\bbZ(}D_-^*/D_-)$ embeds in the nimrep $K_{T_L}(T_{D_-})$, 
namely through the  { expansion of cosets in} $D^*/D=D^*/E^*$ { to
cosets in}
$L^*/E^*$, or $K$-theoretically through  the induction $K^{\mathrm{dim}\,T}_{T_{D}}(T_{E*})\rightarrow K^{\mathrm{dim}\,T}_{T_L}(T_{E*})$. 
More generally that embedding happens whenever the subfactor inclusion is
\textit{type I}.  

Incidentally, this faithful parametrisation of modular invariants by lattices $\cL$ is formally
 very similar to the $(H,0)$ parametrisation for finite groups as discussed in section 3.1:
 $K_\pm$, $N_\pm$ and $\beta$  there correspond respectively to $D_\pm^*/L$, $D_\pm/L$
 and $\beta$ here.

\section{Outer automorphisms}

{ As explained in the introduction, almost all modular invariants for Lie groups
$G$ are combinations of simple-current invariants and outer automorphisms. In this
section we work out the case of outer automorphisms but trivial simple-current invariant,
i.e. the case of pure outer automorphism. We begin with the analogue for finite groups.}

\subsection{Finite groups}

Let $G$ be a finite group and fix any automorphism
$\omega\in \textrm{Aut}\,G$. It is convenient to introduce the semi-direct product $G_\omega=
G \rtimes \langle\omega\rangle$, with product $(g,\omega)(g',\omega')=(g\,\omega(g'),\omega\omega')$.
We will see shortly that inner automorphisms are invisible in the $K$-theoretic
description we seek, so really $\omega$ lives in Out$\,G$. 

Note that $\omega$ permutes the basis of the Verlinde ring, sending $(g,\chi)$ to
$(\omega(g),\chi\circ\omega^{-1})$. This permutation defines 
 an automorphism invariant $\cZ=I_\omega$. We claim that this $\cZ$ is
sufferable,  corresponding
to the pair $(H,\psi)=(\Delta_G^\omega,0)$ where $\Delta_G^\omega$ is
the $\omega$-twisted diagonal $\{(g,\omega(g))\,:\,g\in G\}$ in $G\times G$.

$K^*_{\Delta^{\omega\, L}_G\times \Delta_G^{\omega\, R}}(G\times G)\simeq
K^*_{G^{\mathrm{ad}_\omega}}(G)$, with action $g.x=\omega(g)x\omega(g^{-1})$,
is obviously isomorphic to 
$K_{G^{\mathrm{ad}}}^*(G)$. In other words, the full system (and neutral system) 
is the Verlinde ring.
 Alpha-induction is clear: $\alpha^+$ is the action of $\omega$ on Ver$(G)$,
while $\alpha^-$ is the identity. 

$K^i_{\Delta^{\omega\, L}_G\times\Delta^R_G}(G\times G)\simeq K_{(G,1)^{\mathrm{ad}}}^i((G,
\omega))$, regarding both $(G,1)$ and $(G,\omega)$ as subsets
of $G_\omega$, is the nimrep for $i=0$ and vanishes for $i=1$.  Equivalently, this is $K_{G^{\mathrm{ad}_\omega}}^*(G)$ 
where $G$ acts on itself by the $\omega$-twisted adjoint $g.x=gx\omega(g^{-1})$. The 
preferred basis $\cB$ is obtained automatically from $K$-theory, namely pairs $(g,\chi)$ where now
$g$ is a representative of an $\omega$-twisted conjugacy class in $G$ and $\chi$ is an irrep
of the twisted stabiliser $\cZ^\omega_g(G)$.

More generally, one gets graded products (\textit{twisted fusions}) 
$K_{(G,1)^{\mathrm{ad}}}^0(G,\omega)\times K_{(G,1)^{\mathrm{ad}}}^0(G,\omega')\rightarrow 
K_{(G,1)^{\mathrm{ad}}}^0(G,\omega'\circ\omega)$, so 
$\sum_{\omega'\in\langle\omega\rangle}K_{(G,1)^{\mathrm{ad}}}^0(G,\omega')$ is an associative 
graded ring with unity. All of these products (Verlinde, nimrep, twisted fusion) are simply
the push-forward of product on the space $G$: $g\in G$ acts on $(a,b)\in
K_{(G,1)^{\mathrm{ad}}}^0(G,\omega)\times K_{(G,1)^{\mathrm{ad}}}^0(G,\omega')$ by $((\omega^{-1}g)ag^{-1},
gb(\omega'g)^{-1})$, so the product $(a,b)\mapsto ab$ is clearly compatible
with this $G$-action. 

We can  use $\omega$ to twist any modular invariant: $(\cZ^\omega)_{\lambda,
\mu}=\cZ_{\lambda,\omega(\mu)}$. Suppose $\cZ$ is sufferable, corresponding
to pair $(H,\psi)$. Then $\cZ^\omega$ is also sufferable, corresponding
to pair $((1,\omega)H,\psi^\omega)$, where $(1,\omega)H=\{(h_1,\omega h_2):(
h_1,h_2)\in H\}$ and where $\psi\mapsto \psi^\omega$ is the
 isomorphism $H^2_H(\mathrm{pt};\bbT)\rightarrow H^2_{(1,\omega) H}(\mathrm{pt};\bbT)$ 
obtained from the isomorphism $(1,\omega) H\rightarrow H$.

\subsection{Compact Lie groups}

The most obvious nondiagonal modular invariant is charge-conjugation $\cZ=C$.
In the loop group $LG$ setting, charge-conjugation corresponds to an outer
automorphism of $G$.
This { sub}section develops the $K$-theoretic interpretation of the nimrep,
alpha-induction etc for the modular invariants $I_\omega$ associated to
any outer automorphism of $G$ (inner automorphisms are invisible).
Specialising to the trivial automorphism recovers the Verlinde
ring realisation of \cite{FHT}.

Let $G$ be any connected, compact, simply-connected Lie group and fix a level $k\in\bbZ_{\ge 0}$. 
Write $\kappa=k+h^\vee$ as in secton 2.4. The  group of outer automorphisms of $G$
is naturally identified with the group of symmetries of the Dynkin diagram of $G$, and as
a permutation of these vertices also permutes highest weights of $G$ in the usual 
way and through this the level $k$ primaries $\lambda\in P^k_+(G)$ ($\omega$
fixes $\lambda_0$).
Pick an outer automorphism $\omega$, say of order $d$, and lift it to an automorphism of $G$. 
For example, an automorphism of $G=\mathrm{SU}(n)$ realising 
charge-conjugation automorphism $C$ is complex conjugation.

Introduce the $\omega$-twisted diagonal $\Delta^\omega_G$ and the semi-direct product $G_\omega=
G \rtimes \langle\omega\rangle$ as in section 4.1. 
As with section 4.1, the full system and alpha-induction are easy here. The full system 
$^\kappa K_*^{\Delta_G^{\omega\, L}\times\Delta_G^{\omega\, R}}(G\times G)\simeq {}^\kappa
K_*^{G^{\mathrm{ad}}}(G)$ is just the Verlinde ring in degree $0$ and vanishes in the other.
 Implicit here is our
identification of the twist groups $H^3_{\Delta_G^{\omega\, L}\times\Delta_G^{\omega\, R}}
(G\times G;\bbZ)\simeq H^3_{G^{\mathrm{ad}}}(G;\bbZ)$, using \eqref{freenormcoh} and the 
spectral sequence of section 2.1. Take $\alpha^+$ to be $\omega$, and $\alpha^-$
to be the identity. More interesting is the nimrep: it should be ${}^\kappa 
K_*^{\Delta^{\omega\, L}_G \times\Delta^R_G}(G\times G)\simeq {}^\kappa K_*^{G}(G\omega)
={}^\kappa  K_*^{G^{\mathrm{ad}_\omega}}(G)$ where ad$_\omega$ denotes the 
$\omega$-twisted adjoint action. Again, spectral sequences permit us to identify these
twist groups. For all actions considered in this paragraph, the $H^1$-groups vanish. The 
remainder of this subsection will focus on ${}^\kappa K_*^{G^{\mathrm{ad}_\omega}}(G)$.

The orbits for the  untwisted adjoint action of $G$ on $G$ are parametrised by the Stiefel diagram. 
The analogue for the  $\omega$-twisted adjoint action $G^{\mathrm{ad}_\omega}$ is discussed in 
\cite{MW,Sta}.   For concreteness consider $G={\rm SU}(3)$ and $\omega=$ complex conjugation; 
then the twisted Stiefel diagram (see Figure 1(a)) is a segment $0\le \theta\le \pi/4$ with orbit 
representatives
diag$({\scriptsize\left(\begin{matrix}\cos(\theta)&\sin(\theta) \cr -\sin(\theta)&\cos(\theta)\end{matrix}
\right)},1)$.  The stabiliser at endpoints $\theta=0$ and
$\theta=\pi/4$ are SO$(3)=\textrm{Re}({\rm SU}(3))$ and diag$({\rm SU}(2),1)$ respectively; 
at the generic points the stabiliser is diag$({\rm SO}(2),1)$.

\medskip\epsfysize=1.75in\centerline{ \epsffile{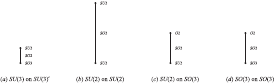}}\medskip

\centerline{Figure 1. Various one-dimensional orbit diagrams}
\medskip

The Dixmier-Douady bundle realising the $H^3$-twist $\kappa\in\bbZ$ here is constructed
exactly like that of $G^{\mathrm{ad}}$ on $G$, using now the twisted Stiefel diagram. In particular,
 the base is the connected component $G\omega$ of $G_\omega$ and the fibres are the 
 compacts $\cK(L^2(G\omega)\otimes\ell^2)$. 

Before we describe what happens in general, 
let's specialise to $G={\rm SU}(3)$ 
  and $\omega$ being complex conjugation for concreteness.
The level enters through the representation $a\mapsto a^\kappa$ of the stabiliser $\bbT\simeq
\textrm{SO(2)}$ at the overlap (the midpoint of Figure 1(a)).
Here  $\kappa=k+3$. From a Mayer-Vietoris 
calculation we see that $\kappa\in H^3_{G}(G\omega;\bbZ)$ traces to $(-1)^\kappa\in 
H^3_{\textrm{SO3}}(\textrm{SO3};\bbZ)$. We
can  calculate $^\kappa K^{G}_*(G\omega )$ using the 6-term sequence (\ref{six-term}) by 
removing the two endpoints of Figure 1(a):
\begin{equation}\label{su3cc}
\begin{array}{ccccc}
R_{\textrm{SO2}} & {\longleftarrow} &
^{\kappa}K^{G}_{0}(G\omega) & {\longleftarrow} & 0\\[3pt]
\alpha\downarrow & & & & \uparrow\\[2pt]
^\mp R_{\textrm{SO3}}\oplus R_{\textrm{SU2}}& {\lori} &
^{\kappa}K^{G}_1(G\omega) &{\lori} &0
\end{array}.
\end{equation}
Here $\mp$ denotes $(-1)^{\kappa+1}$; the only subtlety is that Poincar\'e duality  
introduces an additional $-\in H^3_{\textrm{SO3}}(\textrm{pt};\bbZ)$ twist.

The map $\alpha$ sends $p(a)\in R_{\textrm{SO2}}$ to D-Ind$_{\textrm{SO2}}^{\textrm{SO3}}
(a^{\kappa/2}p,p)$. These Dirac inductions are given explicitly
in \cite[App.A]{FHT}. We find $^\kappa K_0^G(G\omega)=0$ and
\begin{equation}
^\kappa K_1^G(G\omega)=\left\{\begin{matrix}R_{\textrm{SO3}}/(\sigma_\kappa,\sigma_{2i+\kappa}+
\sigma_{\kappa-2i})_{0<i\le (\kappa-1)/2}& \mathrm{for}\ \kappa\ \mathrm{odd} \cr
^-R_{\textrm{SO3}}/(\sigma_\kappa,\sigma_{2i+\kappa}+\sigma_{\kappa-2i})_{0< i\le \kappa/2}& \mathrm{for}\ 
\kappa\ \mathrm{even}\end{matrix}\right.\,.\end{equation}
These $K$-homology groups naturally inherit an $R_G$-module structure by restriction
(the 3-dimensional fundamental representation $\sigma_{(1,0)}\in R_G$
restricts to $\sigma_3\in R_{\textrm{SO3}}\subset R_{\textrm{SU2}}$) and through this in fact
a representation of the fusion ring $\textrm{Ver}_k(\textrm{SU3})=R_G/(\sigma_{(k+1,0)},
\sigma_{(0,k+1)})$: we see by direct calculation that this 
recovers the appropriate fusion graph (see Figure 2) and hence recovers
the nimrep.  Indeed, $\sigma_{(1,0)}\sigma_l=(\sigma_2^2-1)\sigma_l=\sigma_{l+2}
+\sigma_l+\sigma_{l-2}$ for any $l$.

\medskip\epsfysize=.4in\centerline{ \epsffile{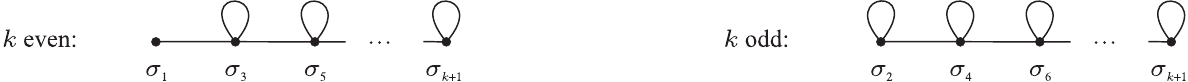}}\medskip

\centerline{Figure 2. Nimrep graphs for SU(3) with charge-conjugation}
\medskip

RCFT has a complete (though still conjectural) description for outer automorphism modular
invariants (\cite{GaGa}, building on \cite{BFS}), which
 paints a picture of the modular invariant $I_\omega$ for outer automorphisms $\omega$, 
 beautifully parallel to that of the diagonal modular invariant. See Table 2. There 
 $\mathfrak{g}^{(1)}$ is the nontwisted affine algebra of $G$ and $\frak{g}^{(\omega)}$ is the
twisted affine algebra of $\omega$-fixed points in $\mathfrak{g}^{(1)}$ \cite{kac}.
Write $P_+^k(\mathfrak{g}')$ for the level $k$ integrable highest-weights of some affine algebra 
$\frak{g}'$. The characters $\chi_x$ for $x\in P_+^k(\mathfrak{g}^{(\omega)})$ possess a
modularity: 
$$\chi_{x}(-1/\tau)=\sum_{\mu\in P_+^k(\check{\mathfrak{g}}{}^{(\omega)})}S_{x,\mu}
\chi_\mu(\tau)$$ 
for another affine algebra $\check{\frak{g}}^{(\omega)}$ called the \textit{orbit Lie algebra}. 
$\Psi$ in Table 2 refers to the matrix diagonalising the nimrep.

$$\begin{array}{|c||c|c|c|c|c|c|c|} \hline
\mathrm{twist}\ \omega&\cB&\rm{exponents}&\Psi&\rm{nimrep} & \rm{full\ system}&\alpha^+&\alpha^- 
\\ \hline
\mathrm{trivial}& P^k_+(\frak{g}^{(1)})&P^k_+(\frak{g}^{(1)}) &S\ {\rm{for}}\ \frak{g}^{(1)}
&R_G/I_k(G)&{\rm{Ver}}_k(G)&\rm{id}&\rm{id}\\ \hline
\mathrm{nontrivial}&  P^k_+(\frak{g}^{(\omega)})&P^k_+(\check{\frak{g}}^{(\omega)}) &S\ {\rm{for}}\ \frak{g}^{(\omega)}
&R^\pm_{G^{(\omega)}}/I_k^\omega(G)&{\rm{Ver}}_k(G)&\omega&{\rm{id}}\\ \hline
\end{array}$$

\centerline{\textbf{Table 2.} Comparing nontrivial $\omega$ twist to trivial twist}\medskip

Table 3 collects the data for all possibilities for nontrivial $\omega$. Here $x'$ gives the map
$P_+^k(\mathfrak{g}^{(\omega)})\rightarrow P_+^{k'}(\widehat{G^{(\omega)}})$, where 
$\widehat{G^{(\omega)}}$ here denotes the universal cover of $G^{(\omega)}$. The twisted fusion
ideal $I_k^\omega(G)$ is generated in part by the fusion ideal $I=I_{k'}(G^{(\omega)})$
or its restriction $I^{sp}=I_{k'}(\widehat{G^{(\omega)}})\cap R^-_{G^{(\omega)}}$ to the spinors.
When $x'\pm J'x'$ appears in Table 3, it is meant to include that relation for all $x'$ and all
simple-currents $J'$; when $x'=J'x'$ then $x'+J'X'$ should be replaced with $x'$. 
In the $E_6$ row, $\pi$ refers to any of the 5 nontrivial outer automorphisms of $D_4$,
and $\epsilon_\pi$ its parity; again $x'+\pi x'$ (or $x'+C'x'$) should be replaced with
$x'$ when $x'=\pi x'$ (or $x'=C'x'$ respectively). The action of Ver$_k(G)$ on
$R^\pm_{G^{(\omega)}}/I^\omega_k(G)$ is by restriction $R_G\rightarrow R^-_{G^{(\omega)}}$
in all cases, except for $G=\mathrm{SU}(2n)$ when $G^{(\omega)}$ isn't a subgroup. That
case uses the \textit{subjoining} described explicitly in \cite{GaGa1}, which uses the 
embeddings sp$(2n)\subset \mathrm{su}(2n)$ and sp$(2n)\subset \mathrm{so}(2n+1)$
to express any $G$-character restricted to $\omega$-fixed points as a virtual $G^{(\omega)}$-character.
 See \cite{GaGa} for any further clarifications on the  conventions we use here.

$$\mbox{\scriptsize$\begin{array}{|c||c|c|c|c|c|c|c|} \hline
G&\omega&\mathfrak{g}^{(\omega)}&\check{\mathfrak{g}}^{(\omega)}&R^\pm_{G^{(\omega)}} & k'&x'\in
P_+^{k'}(G^{(\omega)})& I_k^\omega(G)  \\ \hline
{\mathrm{SU}}(2n+1)& C&A_{2n}^{(2)}&A_{2n}^{(2)} &R^-_{\mathrm{SO}(2n+1)}
&k+2&(x_0+x_1+1;x_1,\ldots,x_{n-1},2x_n+1)&\langle I^{sp},x'+J'x'\rangle\\ \hline
{\mathrm{SU}}(2n)& C&A_{2n-1}^{(2)} &D_{n+1}^{(2)}
&R^-_{\mathrm{SO}(2n+1)}&{k+1}&{(x_0;x_1,\ldots,x_{n-1},2x_n+1)}&\langle I^{sp}\rangle\\ \hline
{\mathrm{Spin}}(2n)& \lambda_{n-1}\leftrightarrow\lambda_n&D_{n}^{(2)} &A_{2n-3}^{(2)}
&R_{\mathrm{Spin}(2n-1)}&k+1&(x_0+x_1+1;x_1,\ldots,x_{n-1})&\langle I,x'+J'x'\rangle\\ \hline
\mathrm{Spin}(8)& {\lambda_1\rightarrow\lambda_3\rightarrow}&D_4^{(3)} &D_{4}^{(3)}
&R_{\mathrm{SU}(3)}&k+3&{(x_0+x_1+x_2+2;x_2,x_{1}+x_{2}+1)}&\langle I,x'-J'x',\\ 
 &\lambda_4\rightarrow\lambda_1&&&&&&x'+C'x'\rangle\\ \hline
E_6& C&E_6^{(2)} &E_6^{(2)}
&R_{\mathrm{Spin}(8)}&k+6&{(x_0+\cdots+x_3+3;}&\langle I,x'-J'x',\\ 
&&&&&&x_1+x_2+x_3+2,x_4,x_3,x_{2}+x_3+1)&x'-\epsilon_\pi\pi x'\rangle\\ \hline
\end{array}$}$$

\centerline{\textbf{Table 3.} Data for  nontrivial $\omega$}\medskip

\cite{GaGa} identified the weights $P^k_+(\check{\frak{g}}^{(\omega)})$ with the $\omega$-fixed points
in $P_+^k(\mathfrak{g}^{(1)})=P_+^k(G)$, i.e. with the exponents of the desired nimrep.
For each $\lambda\in P_+^k(G)$ it defined matrices $\cN^{GG}_\lambda=\Psi\, \mathrm{diag}
(S_{\lambda,\mu}/S_{1,\mu})_{\mu\in P_+^k(\check{\frak{g}}^{(\omega)})}\Psi^{-1}$ for
the unitary matrix $\Psi$ given in the table. Their conjecture is that $\cN^{GG}_\lambda$
defines a nimrep corresponding to $\cZ=I_\omega$, and that this  nimrep product should 
correspond to $\omega$-twisted fusions. The  boundary states $\cB$ then would
 be identified with $P_+^k({\frak g}^{(\omega)})$. { The next proposition summarises what is
 known in the CFT literature.} 

\medskip\noindent\textbf{Proposition 2.} \cite{GaGa,GaGa1,GaVa} \textit{Let $G$ be a compact,
connected, simply-connected Lie group, the level $k$ be any nonnegative integer, and
$\omega$ any outer automorphism of $G$. Then the module structure entries
$\cN^{GG\ y}_{\lambda,x}$ for all $\lambda\in P_+^k(G)$ and all $x,y\in P_+^k({\frak 
g}^{(\omega)})$ are integers, and in particular coincide with the action of $\mathrm{Ver}_k(G)$
on $R^\pm_{G^{(\omega)}}/I^\omega_k(G)$, with respect to the basis $x'$, defined above. 
If all $\cN^{GG\ y}_{\lambda,x}$ are in fact nonnegative, then $\cN_\lambda^{GG}$ defines a 
nimrep corresponding to $\cZ=I_\omega$. Assuming Conjecture 1 of section 2.4, the
map $I_k^\omega(G)\rightarrow\bbZ$ given by dimension of the $G^{(\omega)}$ virtual
representation, is surjective onto $M_k(G)\,\bbZ$; the charge-group $\cM_{\cN^{GG}}$
is $\bbZ_{M_k(G)}$, and is generated by the charge assignment $q_x=\mathrm{dim}(\overline{
x'})$.}\medskip  
   
Does the proposed $K$-homology for the nimrep match the
conjectured $\cN^{GG}$? Can $K$-theoretic methods be used to prove the nonnegativity
needed in Proposition 2? Can $K$-theoretic considerations explain some of the 
arbitrariness in the RCFT description, in particular in the choice of $G^{(\omega)}$
and the presence of subjoinings for $G=\mathrm{SU}(2n)$? { We expect that it should 
be possible to solve all this; easy partial results are:}

\medskip\noindent\textbf{Proposition 3.} \textit{Let $G$ be any compact, connected, 
simply-connected Lie group, any level $k\in\bbZ_{\ge 0}$,
and $\omega$ any automorphism of $G$. Write $\kappa=k+h^\vee$ as usual. 
Then $^{\kappa}K^{G^{\mathrm{ad}}}_{1+\mathrm{dim}\,G+\mathrm{dim}\,G^{(\omega)}}
(G\omega)=0$. A natural basis for $^\kappa K^G_*(G\omega)$ can be
identified with $P_+^k({\frak g}^{(\omega)})$. We get a  natural product $\mathrm{Ver}_k(G)\times {}^\kappa 
K_*^G(G\omega)\rightarrow  K_*^G(G\omega)$, 
written say $[\lambda]*[x]=\sum_{[y]}\cN^{Kth\ [y]}_{\lambda,[x]}[y]$. These coefficients
$\cN_{\lambda,[x]}^{Kth\ [y]}$ are nonnegative integers, and the corresponding
matrices $\cN_\lambda^{Kth}$ form a representation of $\mathrm{Ver}_k(G)$. The
groups $H^3_{G^{\mathrm{ad}}}(G\omega;\bbZ)$ and $H^3(G;\bbZ)$ are naturally identified.
}\medskip

The evaluation of this $K$-homology is from Theorem 1 above. The appearance of
$P_+^k(\mathfrak{g}^{(\omega)})$ exactly matches \cite{GaGa}.
Section 16 of \cite{FHTlg} relates this product Ver$_k(G)\times {}^\kappa K^*_G
(G^{\mathrm{ad}_\omega})\rightarrow  K^*_G(G^{\mathrm{ad}_\omega})$
to Segal's fusion, hence the multiplicities $\cN^{Kth\ [y]}_{\lambda,[x]}[y]$ are nonnegative 
integers. This makes it easy to conjecture the following.

Recall the exterior algebra $\bigwedge\{x_1,\ldots,x_{r-1}\}$ from section 2.4. We can rewrite this 
as $\bigwedge\{x_1,\ldots,x_r\}/(x_1+ \cdots+x_r)$; from the argument of \cite{Doug} the
$x_i$ are associated with the fundamental weights of $G$. 
To any $n$-dimensional sublattice $L$ of the weight lattice
$P=\mathrm{span}_\bbZ\{x_1,\ldots,x_r\}$, we obtain a unique vector $v_L$ (well-defined up to $\pm 1$)
in $\bigwedge^n\{x_1,\ldots,x_r\}$ by taking any basis $y_1,\ldots,y_n$ of $L$ and
forming the wedge product $y_1\wedge\cdots\wedge y_n$.

\medskip\noindent\textbf{Conjecture 2(a)} \textit{The action of $\mathrm{Ver}_k(G)$ on 
$^{\kappa}K^{G^{\mathrm{ad}}}_{\mathrm{dim}\,G+\mathrm{dim}\,G^{(\omega)}}(G\omega)$ 
is equivalent to the conjectured nimrep $\cN^{GG}$ of \cite{GaGa}, corresponding to the modular invariant
$\cZ=I_\omega$.} 

\smallskip\noindent\textbf{(b)} \textit{The D-brane charges for the nimrep ${}^\kappa 
K_*^{G^{\mathrm{ad}}}(G\omega)$
is the natural map $q_\omega:{}^\kappa K_{\mathrm{dim}\,G+\mathrm{dim}\,G^{(\omega)}}^{G}
(G\omega)\rightarrow {}^\kappa   K_{\mathrm{dim}\,G+\mathrm{dim}\,G^{(\omega)}}
(G)$, forgetting $G$-equivariance. 
It is defined by $q_\omega([\overline{\lambda}])=(\mathrm{dim}(\overline{
\lambda})\ (\mathrm{mod}\ M))[v_L]$ where $L$ is the sublattice $(L^\omega)^\perp$
of the weight lattice $P$ orthogonal to the fixed-points $L^\omega$, and $[v_L]$ means 
the corresponding element in $\bigwedge^{\mathrm{dim}\,G-\mathrm{dim}\,G^{(\omega)}}
\{x_1,\ldots,x_r\}/(x_1+
\cdots+x_r)$. This charge can also be  defined by the commuting diagram
\begin{equation}\label{Dchom}
\begin{array}{ccccccc}
0\longrightarrow &I_k^\omega(G) & {\longrightarrow} &
R^{\pm}_{G^{(\omega)}}=K_0^{G^{(\omega)}}(1) & \stackrel{\alpha\,}{\longrightarrow} &
^{\kappa}K^{G}_{0}(G\omega)&\longrightarrow 0\\[3pt]
&\beta\downarrow & &\beta\downarrow & & \gamma\downarrow&\\[2pt]
0\longrightarrow&M\bbZ&\longrightarrow & \bbZ=K_0(1)&\stackrel{\delta\,}{\longrightarrow} &
^\kappa K_0(G)&
\end{array},
\end{equation}
where $\beta$ is dimension, and $\alpha$ is
defined by a combination of inducing $K_*\rightarrow K_*^G$ and the push-forward
of the inclusion $1\hookrightarrow G$.}\medskip

Being Ver$_k(G)$-modules, $\cN^{GG}_\lambda$ and $\cN_\lambda^{Kth}$ are
both uniquely determined by their values on the fundamental
weights $\Lambda_1,\ldots,\Lambda_r$, since the Verlinde ring is a homomorphic image of 
the polynomial ring on those fundamentals.  
Together with Theorem 3 then, we can establish the equality $\cN^{GG}_\lambda=\cN_\lambda^{Kth}$ for all $\lambda$,
and the nimrep property, both conjectured in Conjecture 2(a), if we can establish that
$\cN_{\Lambda_i}^{GG}=\cN_{\Lambda_i}^{Kth}$ for all $1\le i\le r$. 
This is perhaps the most promising way to establish the nimrep property for $\cN^{GG}_\lambda$. 
These coefficients
$\cN^{GG}_{\Lambda_i}$ can be deduced from Proposition 2, or found explicitly in 
\cite{GaGa,GaGa1,GaVa}: for example, for  SU$(2n+1)$ we have 
$\cN_{\Lambda_i, x}^{GG\ y} = \cN_{\Lambda_{ 2n+1-i}, x}^{GG\ y} = 
{N}^{\prime\ y'}_{{\Lambda}'_i,{x}'}$ for $1\le i<n$ and 
$\cN_{\Lambda_n, x}^{GG\ y} = \cN_{\Lambda_{ n+1}, x}^{GG\ y} = 
{N}^{\prime\ y'}_{2{\Lambda}'_i,{x}'}$, where  primes denote the fundamental weights and
 fusions of SO$(2n+1)$ level $k+2$.

 It should be possible to compute $^\kappa K^*_{G^{\mathrm{ad}}}(G\omega)$ and the forgetful
 map $^\kappa K^*_{G^{\mathrm{ad}}}(G\omega)\rightarrow {}^\kappa K^*(G)$ by
 generalising the methods of \cite{M} to nontrivial
$\omega$ using the twisted Stiefel diagrams of \cite{MW}, and then applying the 
 Hodgkin spectral sequence. This would go a long way toward proving Conjecture 2.

Because of this uniform $K$-theoretic description, we would expect all of these
modular invariants $I_\omega$ to be sufferable. Along these lines, 
Verrill \cite{Vrr} (with a gap filled in \cite{W3}) does for the $C$-twisted loop group associated to SU$(2n)$ what Wassermann
et al (see e.g. \cite{W}) did for the nontwisted loop groups, and realised the twisted fusions
using subfactors.

Any modular invariant $\cZ$ can be twisted by an outer automorphism, by multiplying $\cZ$
by the automorphism invariant $I_\omega$.  Moreover, in the subfactor approach the matrix product of
sufferable modular invariants  will itself be sufferable.
It is tempting then to guess one can always apply the outer automorphism $\omega$ to any 
$K$-homological description of the data associated to a modular invariant $\cZ$, to get a $K$-homological description
of $I_\omega \cZ$ or $\cZ I_\omega$. This was proved in the finite group setting at
the end of section 4.1;  the most important example of  this in the loop group setting is
given in section 5.3.

\section{Simple-current modular invariants}

{
The other ingredient of generic modular invariants are the simple-current invariants (see (2.10)
below), also
called the D-series or simple-current orbifold invariants. As before, 
the case of finite group $G$ is a baby version of what happens for Lie groups $G$.}

\subsection{Finite groups}
 { The simple-currents of a finite group $G$} consist of primaries
$(z,\phi)$ where $z\in Z(G)$ and $\phi\in\widehat{G}$ is dimension-1.
Any simple-current of a finite group gives rise to a modular
invariant by \eqref{scinv} (this fails for e.g. the odd levels of
 SU(2)); it will be an automorphism
invariant iff the root of unity $\phi(z)$ has order  equal to the
least common multiple of the orders of $\phi$ in $\widehat{G}$ and $z$ in $Z(G)$. Incidentally,
most finite groups have simple-currents; the only one of order $<168$
 which doesn't is the alternating group $A_5$.

Consider for concreteness $G=D_{2n}:=\langle r,s:r^2=s^{2n}=rsrs=1\rangle$,
the dihedral group with $4n$ elements. Let $\psi_{ij}$, $i,j\in\{0,1\}$, denote
its 4 1-dimensional irreps, defined by $\psi_{ij}(r^as^b)=(-1)^{ia+jb}$;
denote the remaining $n-1$ 2-dimensional irreps by $\chi_k$, $1\le k<n$ 
using obvious notation. Its 8 simple-currents are $z_{hij}:=(s^{hn},\psi_{ij})$ for
all $h,i,j\in\bbZ_2$; apart from the trivial $z_{000}$, they all have order
2. The remaining  $2n^2+6$ primaries are $(s^{hn},\chi_k)$, $(s^{a},\phi_{l})$,
$(rs^h,\psi'_{ij})$ where $1\le a<n$, $l\in\bbZ_{2n}$, $h,i,j\in\bbZ_2$, 
where $\psi'_{ij}((rs^h)^as^{bn})=(-1)^{ai+bj}$ etc. See 
section 3.2 of \cite{CGR} for more details, as well as the $S$ and $T$ matrices.

Recall from section 2.3 that a simple-current $j$ permutes the primaries,
and that it associates to each primary $\mu$ a rational number $Q_j(\mu)$.
The simple-current $z_{hij}$ is order-2 (unless $h=i=j=0$), so
$2Q_{z_{hij}}(\mu)$ is an integer which we'll call the parity. The permutations
and parities for $z_{hij}$ are: $(s^{nh'},\psi_{i'j'})
\mapsto(s^{n(h+h')},\psi_{i+i',j+j'})$ (with parity ${n(jh'+j'h)}$), $(s^{nh'},
\chi_{k})\mapsto(s^{n(h+h')},\chi_{nj+(-1)^jk})$ (parity ${hk+nh'j}$), 
$(s^{a},\phi_{l})\mapsto(s^{nh+(-1)^ha},\phi_{(-1)^hl+nj})$ (parity ${ja+hl}$), and 
$(rs^{h'},\psi'_{i'j'})\mapsto(rs^{\{h'+nh\}_2},\psi'_{i'+i+h'j,j'+nj})$
(parity ${j'h+i+h'j}$), where $\{k\}_2\in\{0,1\}$ is congruent mod 2 to $k$. 
The only fixed-points of nontrivial $z_{hij}$ are $(s^{h'n},\chi_k)$ and
 $(s^a,\phi_l)$ for $z_{010}$, and (when $n$ is even) 
$(rs^{h'},\psi'_{i'j'})$ for both $z_{100}$ and $z_{hh'1}$, $(s^{nh'},\chi_{n/2})$ for 
$z_{0i1}$, $(s^{n/2},\phi_{\pm n/2})$ for $z_{1i1}$, and both 
$(s^{n/2},\phi_0)$ and $(s^{n/2},\phi_n)$ for $z_{1i0}$.

Write $\cZ^{(hij)}$ for the simple-current invariant $\cZ_{\langle z_{hij}\rangle}$.
Among the $\cZ^{(hij)}$, the only nontrivial automorphism invariants occur
when $n$ is odd with $h=j=1$. These automorphism invariants permute the 
primaries, fixing those with even
parity and interchanging $\lambda\leftrightarrow z_{1i1}\lambda$ for $\lambda$ 
with odd parity. The other $\cZ^{(hij)}$   
 are all block-diagonal modular invariants: their only nonzero entries
are $\cZ^{(hij)}_{\lambda,\lambda}=2$ when $\lambda$ is a fixed-point
of $z_{hij}$ with even parity, and $\cZ^{(hij)}_{\lambda,\lambda}=\cZ^{(hij)}_{\lambda,z_{hij}\lambda}=1$
when $\lambda$ has even parity but is not a fixed-point. 
In addition, products of any two distinct $\cZ^{(hij)}$ is a new modular 
invariant.

To keep things simple, we will restrict here to the case of direct relevance
to the loop group setting (our primary interest in this paper).  
Take $H=\langle (b^n,1),(a,a),(b,b)\rangle\simeq Z\times G$ for $Z\simeq\bbZ_2$
the centre of $D_{2n}$. Then by the K\"unneth formula for group cohomology \cite{Brown}, $H^2_H(\mathrm{pt};\bbT)\simeq
\bbZ_2^4$, where one copy of $\bbZ_2$ comes from 
$H^2_{D_{2n}}(\mathrm{pt};\bbT)$, and the other three come from the 
characters of $Z$ and $G$. Again, we will restrict to trivial
$\psi\in H^2_H(\mathrm{pt};\bbT)$ although it should be clear how to modify
this discussion for nontrivial $\psi$. In this case the nimrep is
$K^0_{H^L\times \Delta_G^R}(G\times G)\simeq K^0_{Z^L\times G^{\mathrm{ad}}}(G)
\simeq K^0_{G^{\mathrm{ad}}}(G/Z)$, with basis given by pairs $(Za,\chi)$
for any irrep $\chi$ of $Z_{Za}(G)$. The full system is $K^0_{H^L\times H^R}(G\times
G)\simeq R_Z\otimes K^0_{G^{\mathrm{ad}}}(G/Z)$, i.e. two copies of the nimrep.
Then alpha-induction $\alpha^\pm:K^0_{G^{\mathrm{ad}}}(G)\rightarrow 
K^0_{G^{\mathrm{ad}}}(G/Z)$ is the map $(a,\chi)\mapsto 
(Za,\mathrm{Ind}_{Z_a}^{Z_{Za}}\chi)$ and 
$(a,\chi)\mapsto 
(Za,\mathrm{Ind}_{Z_a}^{Z_{Za}}\chi\circ s^n)$ respectively (i.e.
for $\alpha^-$ premultiply and then induce). We find that the corresponding
modular invariant is $\cZ^{(100)}$. Next subsection we find that a similar
description applies to loop groups.

\subsection{Compact Lie groups}
{ The geometry of simple-current invariants for
 the loop groups  is clear:}
they correspond to strings living on non-simply-connected groups
$G/Z$ where $Z$ is some subgroup of the centre of $G$.
These are related to section 3, except fixed points of the
simple-currents occur here, complicating things considerably.

The simple-currents in Ver$_k(G)$ for any connected, simply-connected, compact 
Lie group $G$, were classified in \cite{Fu}. The group of simple-currents for $G\times H$
is the direct product of those for $G$ and $H$, so it suffices to
consider simple $G$. All of these simple-currents correspond to extended
Dynkin diagram symmetries, with one exception (Ver$_2(E_8))$ which we can ignore 
 as  it cannot yield a modular
 invariant for $G$. For any $G$ and $k$, the simple-currents and outer
automorphisms together generate all symmetries of the extended Dynkin diagram.

For example, the group of
simple-currents for Ver$_k(\mathrm{SU}(n))$ is cyclic of order $n$, generated by
$J=(0;k,0,\ldots,0)$ which permutes $P_+^k(\mathrm{SU}(n))$ through 
 $(\lambda_0;\lambda_1,\ldots,\lambda_{n-1})\mapsto(\lambda_{n-1};\lambda_0,
\lambda_1,\ldots,\lambda_{n-2})$. Then $Q_{J^d}(\lambda)=\frac{d}{n}\sum_{i=1}^{n-1}
i\lambda_i$ and $h_{J^d}=kd(n-d)/2n$. 

First, we need to establish $K$-theoretically the relation, well appreciated
within conformal field theory (see e.g. \cite{FGK}), between the simple-currents and
the centre of $G$. Consider for concreteness
 $G={\rm SU}(2)$. Recall the  Dixmier-Douady bundle for $G$ on $G$ 
constructed in \cite{EG1} and section 2.2, and let $S$ denote the Stiefel diagram, i.e. 
diag$(e^{\i \theta},e^{-\i\theta})$ for $0\le \theta\le\pi$. For $x\in S$ and compact
 $c\in\cK$, define the action of $z=-I\in Z(G)$ on fibres by $(gxg^{-1},c)_1
\leftrightarrow (zgxg^{-1},c)_2$ for all $g\in G$. $G$-equivariance is automatic, as is
the consistency condition for $g\in Z_x(G)$ when $x\in D_1\cap D_2$, where 
$D_i$ is the cover of $S$.
To see that this action preserves the twist, note that $U_k\pi_w U_k^*= \i^{-k}\pi_w$
where $w=\left({0\ \i\atop \i\ 0}\right)$ (Ad$(w)$ moves $zx$ back into the
Stiefel diagram, and sends weightspace $V_{m,n-1}$ to $\i^n V_{-m,1-n}$).  
Therefore the full centre $\pm I$ of $G=\mathrm{SU}(2)$ acts on the $G$ on 
$G$ bundle, and hence on the $K$-group $^\kappa K^1_{G^{\mathrm{ad}}}(G)$. 

We expect that 
the same conclusion should apply to  the centre of any simply-connected,
connected, compact Lie group $G$.
For such $G$, multiplication by the centre $Z(G)$ should correspond naturally 
to the action of the simple-currents in the Verlinde ring $^\kappa 
K_{G^\mathrm{ad}}^{\mathrm{dim}\,G}(G)$, in the following sense. The primaries
$\lambda\in P_+^k(G)$ are identified with certain conjugacy classes ---
this yields a geometric picture of Ver$_k(G)$ dual to
the usual representation ring description $R_G/I_k$.
 Now, $Z(G)$ permutes these conjugacy classes by multiplication,
and this permutation agrees with the simple-current action on primaries. 
For example, for $G=\mathrm{SU}(2)$, the level $k$ primaries look like
$\lambda=(k-\lambda_1;\lambda_1)$ for integers $0\le \lambda_1\le k$; this 
corresponds to the conjugacy class intersecting $S$ at
diag$(\exp[2\pi\i(\lambda_1+1)/2\kappa],\exp[-2\pi\i(\lambda_1+1)/2\kappa])$ 
for $\kappa=k+2$.
By multiplication, the central element $z=-I$ sends the $\lambda_1$ conjugacy
class  to the $\lambda_1+\kappa$ one,
which the Weyl group identifies with $\kappa-\lambda_1-2$. This matches the 
action of the simple-current.

We can also see this identification between $Z(G)$ and simple-currents, in the 
representation ring picture. Again consider for concreteness $G=\textrm{SU}(2)$: 
we can compute its Verlinde ring Ver$_k(G)={}^\kappa K_0^G(G)$ using the
six-term sequence
(\ref{six-term}) by removing three points from Figure 1(b) (the 2 endpoints and
the midpoint):
\begin{equation}\label{su2su2b}
\begin{array}{ccccc}
0 & {\longleftarrow} &
^{\kappa}K^{G}_{0}({G}) & {\longleftarrow} &
R_G\oplus R_T\oplus R_{G}\\[3pt]
\downarrow & & & & \uparrow\beta\\[2pt]
0 & {\lori} &
^{\kappa}K^{G}_1(G) &{\lori} &R_T\oplus R_T
\end{array},
\end{equation}
where $T$ is the maximal torus.
The map $\beta$ sends the polynomials $(p(a),q(a))\in R^2_T$ to 
Dirac induction (D-Ind$_T^Gp,p+a^{-\kappa}q,$D-Ind$_T^Gq)$,
and we obtain $^\kappa K_1^G(G)=0$ and 
\begin{equation}\mathrm{Ver}_k(G)={}^\kappa K_0^G(G)
=R_G^2/\textrm{span}\{(\sigma_i,\sigma_{i+\kappa})\}_{i\in\bbZ}\,.\end{equation}
 Hence Ver$_k(G)$
is spanned by the classes $[(0,\sigma_l)]$ for $1\le l< \kappa$. The centre $Z=\{\pm I\}$
permutes the two $R_G$'s in the upper-right corner of (\ref{su2su2b}), i.e. in 
$K_G^0(Z)\simeq R_G^2$,
 and so sends $[(0,\sigma_l)]$ to $[(\sigma_l,0)]=[(0,\sigma_{\kappa-l})]$.
Again, this permutation matches the simple-current action in Ver$_k(G)$.

For any subgroup $Z$ of the centre $Z(G)$ of any connected, simply-connected, compact
$G$, $R_Z$ is an $R_G$-module with action $\rho_\lambda.\psi=\overline{\mathrm{Res}}_Z^G
(\rho_\lambda)\psi$, where  $\overline{\mathrm{Res}}_{Z}^{G}(\rho_\lambda)$ denotes
the unique  1-dimensional representation appearing in the restriction of the $G$-irrep 
$\rho_\lambda$ to the central subgroup $Z$. The same action defines a Ver$_k(G)$-module
structure, and indeed a  Ver$_k(G)$-nimrep,  on $R_Z$ for any level $k$. The exponents
of this nimrep are the simple-currents corresponding to $Z$ (each with multiplicity 1).
Moreover, for any $z\in Z$ and $\lambda\in P_+^k(G)$, $\overline{\mathrm{Res}}_Z^G(\rho_\lambda)(z)
=e^{2\pi\i Q_j(\lambda)}$, where $j$ is the simple-current corresponding to $z$ and $Q_j$
is defined early in section 2.3. This nimrep $R_Z$ will generally not have a compatible
modular invariant, but it arises indirectly in both this subsection and the next.

Now that we have identified the simple-currents with the centre, we can express
$K$-theoretically the nimreps, full systems, etc for the simple-current
modular invariants. We will do this for any $G=\mathrm{SU}(n)$ in this 
subsection, and comment briefly at the end of the subsection on the related
description for other $G$.

Consider first $G=\mathrm{SU}(2)$ for concreteness. Write $Z=\{\pm I\}$ for
its centre and $\overline{G}=G/Z\simeq\mathrm{SO}(3)$. 
Let $k$ be any nonnegative \textit{even} integer --- 
we require $k$ even here for the existence of the simple-current 
invariant $\cZ_{\langle J\rangle}$, nimrep etc (more on this shortly). 
Write $\kappa=k+2$ as usual. 
The Verlinde ring is described by 
the (unextended) $\bbA_{k+1}$ diagram, and the nimrep
by the $\bbD_{k/2+2}$ diagram, as explained in section 2.3. Apart from this, 
the theory bifurcates between $k/2$ even (where $\cZ_{\langle J\rangle}$
 is block-diagonal) and  $k/2$ odd (where it is an automorphism invariant).

For later convenience, we can redo the Ver$_k($SU(2)) calculation of \eqref{su2su2b}
(for $\kappa$ even), taking $\beta=($D-Ind$_T^Gp,a^{\kappa/2} p+a^{-\kappa/2}q,$D-Ind$_T^Gq)$
and obtaining the equivalent expression for {Ver}$_k(\mathrm{SU}(2))$: 
\begin{equation}\label{su2verl}
{}^\kappa K_0^{G}(G)=R_{T}/\textrm{span}\{a^{\kappa/2}, a^{-\kappa/2},
a^{i+\kappa/2}+a^{-i+\kappa/2},a^{j-\kappa/2}+a^{-j-\kappa/2}\}_{i,j\ge 1}\,.
\end{equation}
The obvious $\bbZ$-basis of Ver$_k(G)$ is $[a^{-\kappa/2+1}],\ldots,[1], 
\ldots,[a^{\kappa/2-1}]$. The $R_G$ action inherited from $K$-theory is
by restriction to $T$, so for example the generator $\sigma_2$ acts by
multiplication by $a+a^{-1}$. In terms of the given basis, multiplication
by $\sigma_2$  recovers the $\bbA_{2\kappa-1}$ diagram. Therefore the $R_G$-module
product descends to the Verlinde ring product.

The nimrep corresponds to ${}^\tau K^{G}_0(\overline{G})\simeq {}^\tau 
K_0^{G^{\mathrm{ad}}\times Z^R}(G)$. This $K$-homology has already been 
computed in \cite{BS}, but it is convenient to recompute it in 
order to compare below various $K$-groups. In addition, our calculation is far 
simpler. The twist $\tau$ here
belongs to $H^1_{G}(\overline{G};\bbZ_2)\times H^3_{G}(\overline{G};\bbZ)\simeq
\bbZ_2\times\bbZ$ by \eqref{H3H1}. The spectral sequence of section 2.1
identifies $H^3_G(G;\bbZ)$ with $H^3(G;\bbZ)\simeq \bbZ$, and $H^3_G(
\overline{G};\bbZ)$ with $2\bbZ\subset H^3(G;\bbZ)$. 
 This factor of 2 plays a crucial role below.

\eqref{freenormcoh} naturally identifies $H^*_{G^{\mathrm{ad}}\times Z^L}(G;A)$
with $H^*_G(\overline{G};A)$ and thus it suffices to construct the Dixmier-Douady bundle for 
$^\tau K_{G^{\mathrm{ad}}\times Z^L}^*(G)$. Focus for now on the $H^3$-twist.
Recall the bundle for $^\kappa K_G^*(G)$ given in \cite{EG1} and  section 2.2,
and use the action of $Z$ on it
(which works for any $\kappa$) to try to enhance the $G$-equivariance to 
$G\times Z$. The obstruction is $(-1)^\kappa$, because 
the forgetful map $H^3_{G\times Z}(G;\bbZ)\rightarrow H^3_G(G;\bbZ)$ is 
$\times 2$. To see this topologically, observe that the orbit diagram for $G\times Z$ on $G$ is
given by that of Figure 1(c), with centralisers $G$, the maximal torus $T$ of $G$, and 
$T\times Z$. In particular, the diagram is folded in 2 from that of $G$ on $G$ (Figure 1(b)).
Let $\kappa'$ be the level of $G\times Z$ on $G$; i.e. the twisting unitary
in the bundle is the representation 
$a^{\kappa'}$ of the centraliser $T$. The corresponding bundle for $G$ on $G$, forgetting
the $Z$ action, unfolds the orbit diagram, giving \textit{two} cuts, each with a twisting unitary
given by the same representation $a^{\kappa'}$ of the same centraliser $T$. This is equivalent
of course to a single cut, with  a twist of $2\kappa'$.

To incorporate an $H^1$-twist on this bundle, split $L^2(G)=\cH_{ns}\oplus
\cH_{sp}$ into nonspinors and spinors --- this defines a $\bbZ_2$-grading
on the space --- and use the odd automorphism $U=\left({0\atop a^{-1}}\,{a\atop 0}
\right)$ (see section 2.2 of \cite{EG1} for a similar construction).

Restrict for now to the twist $\tau=(-;\kappa/2)$ relevant to the modular
invariant $\cZ_{\langle J\rangle}$  --- we use a semi-colon to separate the
$H^1$ and $H^3$ components. The $-$ here
arises from the adjoint shift, as does the $+1$ in the $H^3$-component.
The orbits of $G$ on $\overline{G}$ are given in Figure 1(c). Now use the six-term
sequence (\ref{six-term}), removing the two endpoints of Figure 1(c):
\begin{equation}\label{su2so3}
\begin{array}{ccccc}
0 & {\longleftarrow} &
^{(-;\kappa/2)}K^{G}_{0}(\overline{G}) & {\longleftarrow} &
R_G\oplus {}^+R_{\textrm{O2}}\\[3pt]
\downarrow & & & & \uparrow\alpha\\[2pt]
0 & {\lori} &
^{(-;\kappa/2)}K^{G}_1(\overline{G}) &{\lori} &R_T
\end{array},
\end{equation}
where we have used the nonorientability of the projective plane
$G/{\rm O}(2)$ to obtain the ungraded representation ring $R_{\textrm{O2}}$ (recall the
implicit use of Poincar\'e duality (\ref{PD}) here).
The map $\alpha:R_T\rightarrow R_G\oplus R_{\textrm{O2}}$ sends
$p$ to (D-Ind$_T^Gp,\textrm{Ind}_T^{\textrm{O2}}a^{\kappa/2} p)$.
We find (for $\kappa\ne 0$) that 
${}^{(-;\kappa/2)} K_1^{\textrm{SU2}}(\textrm{SO3})=0$ and
\begin{equation}\label{so3nim}
{}^{(-;\kappa/2)} K_0^{\textrm{SU2}}(\textrm{SO3})=R_{\textrm{O2}}/\textrm{span}\{\kappa_{\kappa/2}, 
\kappa_{j+\kappa/2}+\kappa_{-j+\kappa/2}\}_{j\ge 1}\,.\end{equation}
The $R_G$-module structure arising from $K$-homology is
restriction to O(2) from $G$. For example, the generator
$\sigma_2$ restricts to
$\kappa_1$; in terms of the obvious $\bbZ$-basis of $^\tau K_0^G(\overline{G})$, 
namely $[\delta],[\kappa_0],\ldots,[\kappa_{k/2}]$, we obtain
the (unextended) Dynkin diagram $\bbD_{\kappa/2+1}=\bbD_{k/2+2}$. 
We see that this $R_G$-module multiplication
factors through the fusion ideal and is thus a well-defined action of
Ver$_k(G)$, recovering precisely the nimrep. A different expression for
this nimrep, generalising to other SU$(n)$, is obtained in Theorem 3 below.


When $k/2$ is odd, the modular invariant $\cZ_{\langle J\rangle}$ is an automorphism invariant and
the full system is simply the Verlinde ring, as is the neutral system. 
Alpha-induction 
$\alpha^{\pm}:\textrm{Ver}_k(G)\rightarrow \textrm{Ver}_k(G)$ then
is given by $\alpha^-=id$ and $\alpha^+( [\sigma_i])=[\sigma_i]$ or 
$[\sigma_{\kappa-i}]$ for $i$ odd/even respectively, recovering $\cZ_{\langle J\rangle}$.
              
Much more interesting is $k/2$ even, where $\cZ_{\langle J\rangle}$
is block-diagonal.  The neutral system is given by $^{\tau} K_0^{\overline{G}}(\overline{G})$
where $\tau=(-;+,\kappa/2)\in\bbZ_2\times(\bbZ_2\times\bbZ)$.
Again, the spectral sequence tells us the map from
$H^3_{G}({G};\bbZ)\simeq\bbZ$ to $H^3_{\overline{G}}(\overline{G};\bbZ)\simeq 
\bbZ_2\times\bbZ$ is by multiplication by 2 in the second component. This 
$K$-homology is computed explicitly in \cite{FHT} but again it is convenient to
recalculate it in a slightly different way. The orbit diagram is Figure 1(d);
using (\ref{six-term}) recovers the same sequence as in \cite{FHT}, with the 
map $\gamma:R_{T}^-\rightarrow R_{\textrm{SO3}}\oplus R_{\textrm{O2}}$ sending
the spinor $a^{1/2}p(a)$ to $(\textrm{D}$-$\textrm{Ind}_T^{\textrm{SO3}}
ap(a^2), \textrm{Ind}_T^{\textrm{O2}}a^{\kappa/4}
a^{1/2}p(a))$. The grading on $R_{\textrm{O2}}$ is lost by the implicit
application of Poincar\'e duality. We find ${}^{(-;+,\kappa/2)} K_1^{\textrm{SO3}}(\textrm{SO3})=0$
and
\begin{align}
{}^{(-;+,\kappa/2)} K_0^{\textrm{SO3}}(\textrm{SO3})=&\,\textrm{Ver}_{k/2}(\textrm{SO(3)})
\nonumber\\=&\,R_{\textrm{O2}}/\textrm{span}\{
\overline{\kappa}_{\kappa/4+i}+\overline{\kappa}_{\kappa/4-i}\}_{i=1/2,3/2,...}\,.\end{align}
The appropriate basis is $[\overline{\delta}],[\overline{\kappa}_0],\ldots,[\overline{\kappa}_{k/4}]$,
which recovers the extended fusions of SO(3) at SO(3)-level $k/2$.
This answer is different from, but equivalent to, that given in
\cite{FHT}. The bar's atop the representations here emphasise that
this O(2) lies in SO(3).

In pure extension or type I theories such as these, the neutral system should
embed into the nimrep. To recover this here, compare the corresponding 
six-term sequences term by term: we find that the map 
${}^{(-;+,\kappa/2)}K_0^{\overline{G}}(\overline{G})\rightarrow
{}^{(-;\kappa/2)}K_0^{G}(\overline{G})$ is the obvious injection
$[\overline{\delta}]\mapsto [\delta],[\overline{\kappa}_i]\mapsto[\kappa_{2i}]$.

The full system for $k/2$ even is $^{\tau'} K_0^{\Delta_{G}^L
\times \Delta_{G}^R}(\overline{G}\times \overline{G})\simeq 
{}^{\tau'} K_0^{G^{\mathrm{ad}}\times Z^L\times 
Z^R}(G)\simeq R_Z\otimes_\bbZ {}^{\tau''}K_0^{G\times Z}(G)$.
The twist groups are $H^1_{G\times Z\times Z}(G;
\bbZ_2)\simeq\bbZ_2\times\bbZ_2$ and $H^3_{G\times Z\times Z}(G;\bbZ)\simeq 
2\bbZ\times 2\bbZ$, calculated as before. The twist
relevant to $\cZ_{\langle J\rangle}$ is $\tau'=(-,-;\kappa/2,\kappa/2)$, with $\tau''=
(-;\kappa/2)$. The bundles are constructed as before --- their existence
requires the level $k$ to be even. The orbits for this action are  just those of Figure 1(c)
except with an extra factor of $Z$ everywhere. Thus the computation
of $^{(-,-;\kappa/2,\kappa/2)} K_*^{G\times Z\times Z}(G)$ mirrors 
that of (\ref{su2so3}) except for an extra  factor of $R_{\bbZ_2}$ everywhere:
the full system 
${}^{(-,-;\kappa/2,\kappa/2)} K_0^{G\times Z\times Z}(G)$ consists of two copies of 
the nimrep $\bbD_{k/2+2}$, and ${}^{(-,-;\kappa/2,\kappa/2)} K_1^{G\times Z\times Z}
(G)=0$. This is exactly what we would want.

To identify alpha-induction when $k/2$ is even, use the calculation of
the Verlinde ring (\ref{su2verl}): comparing  (\ref{su2su2b}) and 
(\ref{su2so3}), we see that the natural map $^{\kappa}K_0^G(G)\rightarrow
{}^{(-;\kappa/2)}K_0^G(\overline{G})$ is just induction from SO(2) to O(2):
$a^i\mapsto \kappa_i$. So this gives us
$\alpha^{\pm}: \textrm{Ver}_k(G)\rightarrow {}^{(-;\kappa/2)}K_0^G(\overline{G})\otimes_{R_G}
R_{Z}$. Let $\alpha^-$ be (Ind$_T^{\textrm{O2}},+1)$, and $\alpha^+$ be
(Ind$_T^{\textrm{O2}},\textrm{Res}_{Z}^T)$ (so these send $a^i$ to $(\kappa_{i},1)$ and 
$(\kappa_{i},(-1)^i)$ resp.). The neutral system is therefore the span of even $i$, as it should be.

In Table 4 we collect the data for the simple-current invariants for SU(2).
We also include, in the bottom two rows, the data for Fredenhagen's hypothetical
models \cite{Fr}. To form this table and provide the above story for SU(2), we 
imported (rather than derived) information about  the modular invariants. 
In the standard supersymmetric models associated to loop groups $LG$ (e.g. those
denoted `CIZ' in Table 4), the fermions factor off, allowing everything to
be describable by the bosons (i.e. the SU(2) part). Fredenhagen proposed
a model for SU(2) where the centre $\pm I$ also acts nontrivially on fermions.
 In this case the modular invariant for $k/2$ odd is block-diagonal
while that for $k/2$ even is an automorphism invariant, so we had to accommodate
this in the table. As Fredenhagen noticed, his theory seemed to fit into
the D-brane charge analysis of \cite{BS} if we assign to his theory the opposite
$H^1$-twist than that needed for the standard SO(3) theory.
To our knowledge it is still not yet clear whether Fredenhagen's model
is a completely consistent supersymmetric theory but, 
 at least from the point of view of $K$-homology,
Fredenhagen's model seems to provide a coherent interpretation for the other 
$H^1$-twist. In particular, his neutral system would recover
the new ring structure on ${}^{(+;-,\kappa/2)}K^1_{\overline{G}}
(\overline{G})\simeq {}^{(+;+,\kappa/2)}K_0^{\overline{G}}
(\overline{G})$ (for $\kappa/2$ even) mentioned in Remark 9.3 in \cite{FHT}.
One difference with the standard model  is that ${}^{(+;\kappa/2)}K_1^{{G}}
(\overline{G})\simeq {}^-R_{\mathrm{O2}}^1$ is 1-dimensional, not 0.

$$\begin{array}{|c||c|c|c|} \hline
&\rm{nimrep} & \rm{neutral\ system} & \rm{full\ system} 
\\ \hline
k/2\ \mathrm{even: CIZ}& {}^{(-;\kappa/2)} K_0^{\mathrm{SU2}}(\mathrm{SO(3)})&
{}^{(-;+,\kappa/2)} K_0^{\mathrm{SO3}}(\mathrm{SO(3)}) &
{}^{(-,-;\kappa/2,\kappa/2)}K_0^{\mathrm{SU2}\times Z\times Z}(\mathrm{SU(2)})\\ \hline
k/2\ \mathrm{odd:CIZ}& {}^{(-;\kappa/2)} K_0^{\mathrm{SU2}}(\mathrm{SO(3)})&
{}^\kappa K_0^{\mathrm{SU2}}(\mathrm{SU(2)})& {}^\kappa
K_0^{\mathrm{SU2}}(\mathrm{SU(2)})\\ \hline
k/2\ \mathrm{even: Fr}& {}^{(+;\kappa/2)} K_0^{\mathrm{SU2}}(\mathrm{SO(3)})&
{}^{\kappa} K_0^{\mathrm{SU2}}(\mathrm{SU(2)}) &{}^{\kappa}
K_0^{\mathrm{SU2}}(\mathrm{SU(2)})\\ \hline
k/2\ \mathrm{odd:Fr}& {}^{(+;\kappa/2)} K_0^{\mathrm{SU2}}(\mathrm{SO(3)})&
{}^{(+;+,\kappa/2)} K_0^{\mathrm{SO3}}(\mathrm{SO(3)})& {}^{(+,+;\kappa/2,\kappa/2)}
K_0^{\mathrm{SU2}\times Z\times Z}(\mathrm{SU(2)})\\ \hline
\end{array}$$

\centerline{\textbf{Table 4.} Data for the SO(3) theory, including Fredenhagen's
hypothetical one}\medskip

From this it is easy to describe what should happen in general. Write $G_n={\rm SU}(n)$.
Fix any divisor $d|n$ and level $k$, and write $n'=n/d$ and $\kappa=k+n$. 
The simple-current invariant $\cZ_{\langle J^{n'}\rangle}$ exists iff
$n'(n+1)k$ is even. Write $Z_d$ for the order-$d$ subgroup of the centre   
 $Z(G_n)\simeq \bbZ_n$. The corresponding nimrep should be 
${}^\tau K^{G_n}_0(G_n/Z_d)\simeq {}^\tau K^{G_n\times Z_d}(G_n)$ for some twist
$\tau$; the neutral system should be ${}^{\tau'}K^{G_n/Z_{d'}}_0(G_n/Z_{d'})$ for some
subgroup $Z_{d'}$ of $Z_d$ and twist $\tau'$; the full system should be
${}^{\tau''}K_0^{\Delta_{G_n}^L\times\Delta_{G_n}^R}(G_n/Z_{d'}\times G_n/Z_{d'})\simeq
{}^{\tau''}K_0^{G_n^{\mathrm{ad}}\times Z_{d'}^{L}\times Z_{d'}^{R}}(G_n)\simeq R_{Z_{d'}}
\otimes_\bbZ{}^{\tau'''}K_0^{G_n^{\mathrm{ad}}\times Z_{d'}^{L}}(G_n)$ for appropriate twists
$\tau'',\tau'''$. These $K$-homology groups all vanish in degree 1. See
Conjecture 3 below for more details.

A nimrep for $\cZ_{\langle J^{n'}\rangle}$ has already been proposed in the 
conformal field theory literature, which we now describe. Write $[\lambda]_d$ for the 
orbit $\langle J^{n'}\rangle\lambda$, and $o_d(\lambda)$ for the largest 
positive integer $o$ dividing $d$ such that $J^{n/o}\lambda=\lambda$; hence the stabiliser 
stab$_{Z_d}(\lambda)=\langle J^{n/o_d(\lambda)}\rangle\simeq\bbZ_{o_d(\lambda)}$ and
$\|[\lambda]_d\|=d/o_d(\lambda)$. For readability, we will usually drop the subscript on 
$[*]_d$ and $o_d$. Write $\cP_d^k(G_n)=\{[\nu]_d:\nu\in
P_+^k(G_n)\}$, the set of all $J^{n'}$-orbits. Let $\cF_{n,k,d}$ be the set of all
$J^{n'}$-fixed-points, i.e. all $\lambda\in P_+^k(G_n)$ such that $J^{n'}
\lambda=\lambda$. Then $\cF_{n,k,d}$ is nonempty iff $d$ divides both
$n$ and $k$. The map $\cF_{n,k,d}\rightarrow P_+^{k/d}(G_{n'})$ given by
$\nu\mapsto \nu^{(d)}=(\nu_0;\nu_1,\ldots,\nu_{n'-1})$ is a bijection.

As noted in \cite{GaGa2}, there are two possibilities which behave somewhat
differently (recall $n'(n+1)k$ must be even):

\smallskip\noindent\textbf{Case A:} \textit{Either $n'(n+1)$ is even or the power of 2 
dividing $k$ exceeds that of $n$}. Then the diagonal entries of 
$\cZ_{\langle J^{n'}\rangle}$ 
(i.e. the exponents of the corresponding nimrep) are 
precisely those $\mu\in P_+^k(G_n)$ with $Q_{J^{n'}}(\mu)\in\bbZ$,
i.e. for which $d$ divides $\sum_i i\mu_i$. Such an exponent $\mu$ has multiplicity $o_d(\mu)$.

\smallskip\noindent\textbf{Case B:} \textit{$n'(n+1)(k+1)$ is odd, 
and the power of 2 dividing $n$ is at least as large as that dividing $k$.}
Then the exponents can have $Q_{J^{n'}}(\mu)\in\frac{1}{2}\bbZ$ (see
\cite[Sect.1.2]{GaGa2} for details).\medskip

Most of the study of nimreps for $\cZ_{\langle J^{n'}\rangle}$ has been directed
at the simpler and more common Case A, and in the following \textit{we will 
restrict to Case A}. 
So for example this is automatically satisfied for $n=d=3$ and 
$n=4,d=2$, but requires 4 to divide $k$ for $n=d=2$. 
Extending the following results to Case B, and then to the
other $G$, should be a priority.

An obvious candidate for a nimrep compatible with $\cZ_{\langle J^{n'}\rangle}$
would be the $J^{n'}$-orbits $\mathrm{Ver}_k(G_n)/\langle J^{n'}\rangle$, 
or equivalently the quotient of $R_{G_n}$ by the
ideal generated by the fusion ideal $I_k(G_n)$ together with the
terms $\rho_{\overline{J^{n'}\lambda}}-\rho_{\overline{\lambda}}$ for all
$\lambda\in P_+^k(G_n)$. This Ver$_k(G_n)$-module has a basis parametrised by
the orbits $\cP_d^k(G_n)$, with coefficients 
\begin{equation}\label{naivecoeff}\cN_{\lambda,[\nu]}^{[\nu']}
=\sum_{i=1}^{d/o(\nu)}N_{\lambda,\nu}^{J^{n'i}\nu'}=\frac{o(\nu)}{\mathrm{gcd}(o(\nu),
o(\nu'))}\sum_{i=1}^{d/\mathrm{lcm}(o(\nu),o(\nu'))}N_{J^{n'i}\lambda,\nu}^{\nu'}\,,
\end{equation}
using \eqref{scfus}. But in general this is not
quite a nimrep: the asymmetry between $\nu,\nu'$ means that 
the transpose condition of Definition 2 can fail for $\lambda=C\lambda$.
Moreover, the dimension is wrong: although the exponents all have $Q_{J^{d'}}(\mu)
\in\bbZ$ (which is what we want), they all have multiplicity 1. Indeed,
$\mathrm{Ver}_k(G_n)/\langle J^{n'}\rangle$ will be a nimrep compatible with 
$\cZ_{\langle J^{n'}\rangle}$, iff there are
no fixed-points, i.e. iff gcd$(d,k)=1$. Identifying the correct
nimrep for $\cZ_{\langle J^{n'}\rangle}$, when there are fixed-points,
is more subtle: it requires \textit{resolving} the fixed-points.
The answer is finally provided in Theorem { 3} below.

\cite{BFS,GaGa2} proposed a complicated formula for a nimrep $\cN^{(n,k,d)}$ 
for  $\cZ_{\langle J^{n'}\rangle}$, using \eqref{nimver} with $\Psi$ taken
from an expression for the $S$-matrices
of nonsimply-connected groups conjectured in \cite{FSS}. The boundary states 
$a\in\cB$ consist of pairs $([\nu]_d,l)$ where $[\nu]_d\in \cP_d^k(G_n)$ 
and $l\in\bbZ_{o_d(\nu)}$. This component $l$ resolves the fixed-point $\nu$.
In particular, 
\begin{equation}
\Psi_{([\nu],l),(\mu,i)}=\frac{\sqrt{d}}{o(\nu)\sqrt{o(\mu)}}\sum_{\delta|\mathrm{gcd}(
o(\nu),o(\mu))}\xi_\delta \,S^{(\delta)}_{\nu^{(\delta)},\mu^{(\delta)}}\sum_{\ell\in
\bbZ_\delta^\times}e^{2\pi\i\ell\,(l-i)/\delta}\,,\end{equation}
where $i\in\bbZ_{o(\mu)}$, $S^{(\delta)}$ is the $S$-matrix for $G_{n/\delta}$ at level $k/\delta$, 
and $\xi_\delta$ is some root of unity depending
on $\delta,n,k,d$. By construction, the resulting $\cN^{(n,k,d)}$ will be a nimrep compatible
with $\cZ_{\langle J^{n'}\rangle}$ iff the coefficients $\cN^{(n,k,d)\ ([\nu'],l')}_{\lambda,([\nu],l)}$
arising in \eqref{nimver} are all nonnegative integers. However, in this description both
integrality and nonnegativity are highly unobvious. Using the \textit{fixed-point factorisation}  of
\cite{GaWa2}, \cite{GaGa3} found a relatively simple expression (given below) for some of these 
coefficients, and from this could prove integrality,
but nonnegativity remained out of reach. In Theorem { 3} below, we use this
 to find a  simple global description for the Ver$_k(G_n)$-module structure, 
 making its relation to $K$-homology more evident, and allowing us to finally prove nonnegativity
 and establish the nimrep property. 

We can make explicit the isomorphism of the stabiliser stab$_{Z_d}(\nu)\simeq \bbZ_{o_d(\nu)}$
by writing $\phi_{\nu}^l$ for the irrep sending $J^{hn/o(\nu)}$ to 
$e^{2\pi\i lh/o(\nu)}$ for all
$h\in\bbZ$ and each $l\in\bbZ_{o_d(\nu)}$; it is convenient to extend its 
definition to all of $Z_d$ by defining $\phi_\nu^l(J^{hn'})=0$ 
whenever $J^{hn'}\not\in \mathrm{stab}_{Z_d}(\nu)$, i.e. whenever $d/o(\nu)$ does not divide $h$.
Note that  $\phi_\nu^l=\frac{o(\nu)}{d}\sum_{i=1}^{d/o(\nu)}\psi^{l+io(\nu)}\in
\frac{o(\nu)}{d}R_{Z_d}$, where $\phi^l$ are the $d$ irreps of $Z_d$.
Identify the boundary state $([\nu],l)$ with $([\nu],\frac{d}{o(\nu)}\phi_\nu^l)$; the
strange normalisation of $\phi_\nu^l$ is needed for it to lie in $R_{Z_d}$, and will 
 absorb the extra factors appearing in \eqref{naivecoeff}.

{ Theorem 3 is one of the main results of this paper, giving a significant simplification
to the nimrep expressions appearing in the CFT literature.}

\medskip\noindent\textbf{Theorem 3.} \textit{Fix $G_n=\mathrm{SU}(n)$, level $k$, and
divisor $d|n$. Write $\kappa=k+n$, $n'=n/d$, and $Z_d$ for the order-$d$ subgroup of the
centre of $G_n$. 
Assume that Case A holds. Define $I_k^d(G_n)$
to be the ideal of $R_{G_n\times Z_d}$ generated by the fusion ideal $I_k(G_n)$, together
with the terms $\rho_{{J^{n'}\lambda}}-\rho_{{\lambda}}$ for all
$\lambda\in P_+^k(G_n)$, as well as $\rho_\lambda\otimes\phi-\rho_\lambda\otimes\phi'$, for any 
$\phi,\phi'\in R_{Z_d}$ which are equal when restricted to stab$_{Z_d}(\lambda)$. Then
$R_{G_n\times Z_d}/I_k^d(G_n)$ is a Ver$_k(G_n)$-module, where Ver$_k(G_n)$ acts  by 
multiplication in the obvious way (i.e. ignoring the $Z_d$ component).
A basis for this is $([\nu]_d,l)=([\nu]_d,\frac{d}{o_d(\nu)}\phi_\nu^l)\in\cB$; in terms of this basis we have the coefficients}
\begin{equation}\label{scnimcoeff}
\cN_{\lambda,([\nu],l)}^{\ ([\nu'],l')}=\sum_{j=1}^{d/\mathrm{lcm}(o(\nu),o(\nu'))}
\delta^{\mathrm{gcd}(o(\nu),o(\nu'))}_{l,l'}N_{J^{n'j}\lambda,\nu}^{\ \nu'}\,,
\end{equation}
\textit{where $\delta^m_{l,l'}$ equals 1 if $m$ divides $l-l'$, and 0 otherwise.
This nimrep is equivalent to the one, $\cN^{(n,k,d)}$, described above; in particular, 
$\cN^{(n,k,d)}$ is a nimrep.}

\medskip\noindent\textit{Proof.} To prove this, first derive
\eqref{scnimcoeff}, which follows straightforwardly from \eqref{naivecoeff}.
The expressions for the coefficients of $\cN^{(n,k,d)}$ is much more difficult,
but the Appendix of \cite{GaGa3} uses the fixed-point factorisation of
\cite{GaWa2} to compute
\begin{equation}
\cN_{\Lambda_m,([\nu],l)}^{(n,k,d)\ ([\nu'],l')}=\sum_{j=1}^{d/\mathrm{lcm}(o(\nu),o(\nu'))}
\delta^{\Delta}_{l,l'}N_{J^{n'j}\Lambda_{m/\Delta},\nu^{(\Delta)}}^{\ \nu^{\prime (\Delta)}}\,,
\label{fwnimsc}\end{equation}
valid for all fundamental weights $\Lambda_m$ and all $([\nu],l),([\nu'],l')\in\cB$, where 
$\Delta=\mathrm{gcd}(o(\nu),o(\nu'))$ and 
 $N^{(\Delta)}$ are the fusion coefficients for
Ver$_{k/\Delta}(G_{n/\Delta})$.  From the Pieri rule, or Claim (c) in \cite[Sect.4.1]{GaGa2},
these fusions in Ver$_{k/\Delta}(G_{n/\Delta})$ coincide with the corresponding
ones in Ver$_{k}(G_{n})$, and we find that 
$\cN_{\Lambda_m,([\nu],l)}^{(n,k,d)\ ([\nu'],l')}=
\cN_{\Lambda_m,([\nu],l)}^{\ ([\nu'],l')}$ for all $m,\nu,\nu',l,l'$.

Since these two Ver$_k(G_n)$-modules have identical formulas for multiplication
by the fundamental weights $\Lambda_m$, which generate Ver$_k(G_n)$, they are
isomorphic as Ver$_k(G_n)$-modules, and in fact $\cN_\lambda=\cN^{(n,k,d)}_\lambda$
for all $\lambda\in P_+^k(G_n)$. Therefore the coefficients of $\cN^{(n,k,d)}$
are nonnegative integers, and so $\cN=\cN^{(n,k,d)}$ is indeed a nimrep
compatible with $\cZ_{\langle J^{n'}\rangle}$. QED\medskip

Theorem { 3} provides a massive simplification of the original expressions for
$\cN^{(n,k,d)}$  in \cite{BFS,GaGa2}. Although integrality of the matrix entries 
of $\cN^{(n,k,d)}$ was established in \cite{GaGa3}, nonnegativity was out-of-reach.
We were led to Theorem { 3} by trying to match the conjectured $\cN^{(n,k,d)}$
to the $K$-homology $^\tau K_*^{G\times Z_d}(G)$. Working out this nimrep in
Case B is a natural
task. The fixed-point factorisation of \cite{GaWa2} continues to hold, so one
could extend the nimrep coefficient calculations of \cite[App.B]{GaGa3} to
Case B, although this would be technically challenging, and then try to
reinterpret the result globally as a quotient of a representation ring. Note
though that Theorem { 3} will fail as stated for Case B: in particular, 
not all exponents of $\cZ_{\langle J^{n'}\rangle}$  will now obey 
$Q_{J^{n'}}(\mu)\in\bbZ$, so the simple-current $J^{n'}\in\mathrm{Ver}_k(G_n)$
will not act trivially.

The charge-group for $\cN^{(n,k,d)}$, in Case A, was conjectured in \cite{GaGa3} to be  
\begin{equation}
\bbZ_M\oplus \bigoplus_{p|\mathrm{gcd}(d,M)}\bigoplus_{i=1}^\delta (p^i-p^{i-1})\cdot
\bbZ_{p^{\mathrm{min}(\mu,\nu-i+1)}}\label{chgpsc}\end{equation}
where $M$ is the gcd of the dimensions of all weights in the fusion ideal, and 
$p^\nu\|n$, $p^\delta\|d$ and $p^\mu\|M$. The first sum runs over all primes dividing both
$d$ and $M$. The charge assignments which should generate this, and the evidence 
supporting this conjecture, is discussed in \cite[Sect.4.2]{GaGa3}. With the new global
interpretation of $\cN^{(n,k,d)}$ coming from Theorem { 3}, it is very possible this conjecture
can now be proved. The charges and charge-groups in Case 2 are discussed in
\cite[Sect.5]{GaGa2}.

Before we can explicitly introduce the relevant $K$-homology groups, we need
to understand the appropriate cohomology groups. Note from \eqref{H3H1} that
$H^1_{G_n}(G_n/Z_d;\bbZ_2)$ is $\bbZ_2$ or 0 depending on whether or not $d$ is
even, and $H^3_{G_n}(G_n/Z_d;\bbZ)\simeq\bbZ$. The usual spectral sequence
argument gives $H^3_{G_n/Z_d}(G_n/Z_d;\bbZ)\simeq
\bbZ_d\times \bbZ$, using $H^*_{G_n/Z_d}(\mathrm{pt};\bbZ)\simeq \bbZ\oplus 0\oplus
0\oplus\bbZ_d\oplus\bbZ\oplus\cdots$ and $H^3(G_n/Z_d;\bbZ)\simeq \bbZ\oplus 0
\oplus 0\oplus\bbZ\oplus\cdots$; likewise $H^1_{G_n/Z_d}(G_n/Z_d;\bbZ_2)
\simeq \bbZ_2$ or 0, again depending on whether or not $d$ is even.
Finally, $H^3_{G_n\times Z_d\times Z_d}(G_n;\bbZ)\simeq R_{Z_d}\otimes_\bbZ
H^3_{G_n}(G_n/Z_d;\bbZ)$,  and $H^1_{G_n\times Z_d\times Z_d}(G_n;\bbZ_2)\simeq \bbZ_2(Z_d)$
or 0, depending again on whether or not $d$ is even.


\medskip\noindent\textbf{Conjecture 3.} \textit{Use the same notation as Theorem
4; we require $n'(n+1)k$ is even (necessary for the existence of $\cZ_{\langle
J^{n'}\rangle}$)  but don't assume Case A. Write $t=1$ or $2$ when 
$n'(n+1)$ is even respectively odd. Fix $s\in \bbZ_2$ (which can be dropped if $d$ is  odd).}

\smallskip\noindent\textbf{(a)} \textit{$Z(G_n)\simeq \bbZ_n$ acts on the $G_n$ on 
$G_n$ bundle, and the corresponding action on
$^\kappa K^{G_n}_0(G_n)\simeq\mathrm{Ver}_k(G_n)$ is by simple-currents.}

\smallskip\noindent\textbf{(b)} \textit{
$^{(s;\kappa/t)}K^{G_n}_1({G}_n/Z_d)=0$  and $^{(s;\kappa/t)}K^{G_n}_0({G}_n/Z_d)$ is
a nimrep for the simple-current modular invariant $\cZ_{\langle J^{n'}\rangle}$. 
This should agree with the (conjectured) nimrep  $\cN^{(n,k,d)}$.}

\smallskip\noindent\textbf{(c)} \textit{Let $d_0=\mathrm{gcd}(d,kn')$ if
$d$ is odd, and $ d_0=\mathrm{gcd}(d,kn'/2)$ if $d$ is even. Write 
$Z_0$ for the order-$d_0$ subgroup of the centre of $G_n$.
 Then the neutral system is $^{(s;\tau',\kappa /d_0)} 
K_0^{G_n/Z_0}(G_n/Z_0)$   for some choice of $\tau'\in\bbZ_{d_0}$ ($s$ can 
be dropped if $d_0$ is odd).
 The full system is $^{(s,\ldots,s;\kappa/t,\ldots,\kappa/t)}K_0^{G^{\mathrm{ad}}_n\times Z_0^L
 \times Z_0^R}(G_n)\simeq R_{Z_0}\otimes_\bbZ
{}^{(s;\kappa/t)} K_0^{G_n}(G_n/Z_0)$, i.e. $d_0$ copies of $^{(s;\kappa/t)}K_0^{G_n}
(G_n/Z_0)$ (again, $s$ can be dropped if $d_0$ is odd). 
The corresponding twisted equivariant $K_1$-homology groups all vanish. Write $\beta$ for
 the map $^\kappa K_0^{G_n}(G_n)\rightarrow {}^{(s;\kappa/t)}K_0^{G_n}(G_n/Z_0)$
coming from the projection $G_n\rightarrow G_n/Z_0$. Then $\alpha^+:\mathrm{Ver}_k(G_n)
\rightarrow R_{Z_0}\otimes_\bbZ{}^{(s;\kappa/t)} K_0^{G_n}(G_n/Z_0)$ can be chosen
to be $[\lambda]\mapsto (1,\beta(\lambda))$ while $\alpha^-$ is
$(\overline{\mathrm{Res}}_{Z_0}^{G_n},\beta)$ $(\overline{\mathrm{Res}}_{Z_0}^{G_n}$ 
was defined earlier this section). This map $\beta$ can 
also be interpreted as $^\kappa K_0^{G_n}(G_n)\rightarrow {}^{(s;\kappa)}K_0^{G_n\times Z_0}(G_n)$ where $Z_0$ 
acts trivially by adjoint; in $K$-theory alpha-induction would 
involve induction from $G_n$ to $G_n\times Z_0$.}

\smallskip\noindent\textbf{(d)} \textit{
The D-brane charge group here is $^{\kappa/t} K_0({G_n}/Z_d)$. The assignment
of charges is given by the forgetful map $^{(s;\kappa/t)} K_0^{G_n}({G_n}/Z_d)
\rightarrow{}^{\kappa/t}K_0(G_n/Z_d)$.}\medskip

The possibility of choosing either sign for $s$ is suggested by the Fredenhagen model
\cite{Fr}.

Something similar to Conjecture 3 will hold for all other compact connected groups, and this 
should be worked out. As a first step, the analogue of fixed-point factorisation has
been found for all $G$ \cite{BeGa}. That there can possibly
 be subtleties is hinted by $G=\mathrm{Spin}(8)$ at level 2 (so $\kappa=8$), 
for $Z$ the full centre $\bbZ_2\times\bbZ_2$. The corresponding modular
invariant is
\begin{equation}
\cZ=|\chi_{(2;0000)}+\chi_{(0;2000)}+\chi_{(0;0020)}+\chi_{(0;0002)}|^2+4|\chi_{(0;0100)}|^2\,,
\end{equation}
using obvious notation. The neutral system is $^{17}K_0^{E_7}(E_7)$, which is
2-dimensional (and not 5-dimensional), with sigma-restriction sending those
2 primaries to $\chi_{(2;0000)}+\chi_{(0;2000)}+\chi_{(0;0020)}+\chi_{(0;0002)}$
and $2\chi_{(0;0100)}$ (and not the naive guess $\chi_{(2;0000)}+\chi_{(0;2000)}+
\chi_{(0;0020)}+\chi_{(0;0002)}$, $\chi_{(0;0100)}$, $\chi_{(0;0100)}$, 
$\chi_{(0;0100)}$, $\chi_{(0;0100)}$). The proper way to think of this is
described at the end of  \cite[Sect.6.4]{FSS}. It would be interesting
(and not too difficult) to determine whether $^\tau K^{\mathrm{PSpin}(8)}_0(\mathrm{PSpin}(8))$
for the appropriate twist $\tau$ recovers those 2 dimensions (and not for
instance the 5).

\subsection{Mixing outer automorphisms and simple-currents}

The final generic source of modular invariants results from combining 
simple current  invariants and outer automorphisms.  Surprisingly, the analysis here is 
much simpler than in section 5.2. In this subsection we will restrict to
the modular invariant $\cZ^c_{\langle J^{n'}\rangle}=C\cZ_{\langle J^{n'}
\rangle}=\cZ_{\langle J^{n'}\rangle}C$ for $G_n=\mathrm{SU}(n)$, where
$C$ is charge-conjugation and $n'=n/d$ for some $d|n$. We retain
all notation from section 5.2. In particular, recall the Ver$_k(G_n)$-nimrep $R_{Z_d}$,
with multiplication  $\lambda.\phi_d^l=\phi_d^{l+\sum_jj\lambda_j}$.
{ We obtain the complete answer, up to an earlier conjecture.}

We'll begin by describing what appears in the CFT literature.
Consider first $d$ odd. Let $\cN^c$ be any nimrep for $\cZ=I_c$. It was proved 
in \cite[Sect.2.1]{GaGa3} that $R_{Z_d}\otimes_{R_{G_n}}\cN^c$ is a nimrep 
compatible with $\cZ^c_{\langle J^{n'}\rangle}$. 
In particular, recall the conjectured nimrep $\cN^{GG}$ of section 4.2. Then
$R_{Z_d}\otimes \cN^{GG}$ will be a nimrep for 
$\cZ^c_{\langle J^{n'}\rangle}$ iff the coefficients $\cN^{GG\ y}_{\lambda,x}$
are all nonnegative, as was conjectured. Denoting the boundary states of
$\cN=R_{Z_d}\otimes\cN^{GG}$ by $(x,l)=(x,\phi_d^l)$, for $x\in P_+^k(A_{n-1}^{(2)})$ 
and $l\in\bbZ_d$, we obtain the coefficients
$\cN_{\lambda,(x,l)}^{(y,l')}=\cN^{GG\,y}_{\lambda,x} \delta^d_{l',l+\sum_jj
\lambda_j}$ for $\delta^d$ defined in Theorem { 3}. The charge-group was proved in 
\cite[Sect.2.2.2]{GaGa3} to be \eqref{chgpsc}, subject only to the validity of 
Conjecture 1(a) for $G=\mathrm{Spin}(2n+1)$, and the charge assignments were
all identified.

When $d$ is even, this construction fails because $J^{n/2}$ is an exponent of
$I_c$, so the vacuum 1 has multiplicity 2, not 1, in $R_{Z_d}\otimes
\cN^c$. But this also suggests the cure, at least when $n'$ is even.
The Ver$_k(G_n)$-module $R_{Z_{2d}}\otimes\cN^c$ consists of two identical  copies
of what we want. Restrict for example to the submodule spanned (over $\bbZ$)
by the boundary states $(x,\phi^l)$, $x\in P_+^k(A^{(2)}_{n-1})$, where
$l\equiv\sum ix_i$ (mod 2). Assuming again that the coefficients of 
$\cN^{GG}$ are nonnegative for $G_n$, this submodule of $R_{Z_{2d}}\otimes
\cN^{GG}$ will be a nimrep for $\cZ^c_{\langle J^{n'}\rangle}$. Its
charge-group was conjectured in \cite{GaGa3} to be \eqref{chgpsc} --- see
section 3.2 there for some partial results and supporting evidence.

When $n'(n+1)$ is odd, the situation is more complicated because the exponents
of $\cZ^c_{\langle J^{n'}\rangle}$ are not directly related to those of $I_c$
(by contrast, the exponents for $\cZ^c_{\langle J^{n'}\rangle}$, when $d$
is odd say, are $J^{jn'}\mu$ of multiplicity $o(\mu)$ for any $\mu=C\mu$).
This means there won't be a direct relation between the nimreps for
$I_c$ and $\cZ^c_{\langle J^{n'}\rangle}$. Indeed, when $n'(n+1)$ is odd,
we will no longer in general have agreement between the charge-groups of
$\cZ_{\langle J^{n'}\rangle}$ and $\cZ^c_{\langle J^{n'}\rangle}$: e.g.
as mentioned in \cite{GaGa3}, for $n=k=d=4$, $\cM(\cN^{(4,4,4)})\simeq
\bbZ_4\times\bbZ_2^2$ while $\cM(\cN^{(4,4,4)c})\simeq \bbZ^2_4$.

Now let's turn to a $K$-homological interpretation, focussing on the simpler case where
$d$ is odd. The groups $H^3_{G_n^{\mathrm{ad}_c}}
(G_n/Z_d;\bbZ)\simeq H^3_{G_n^{\mathrm{ad}_c}\times Z_d^L}
(G_n;\bbZ)\simeq\bbZ$ and $H^1_{G_n^{\mathrm{ad}_c}}
(G_n/Z_d;\bbZ_2)\simeq H^1_{G_n^{\mathrm{ad}_c}\times Z_d^L}
(G_n;\bbZ_2)\simeq\bbZ_{\mathrm{gcd}(d,2)}$, calculated as before. The Dixmier-Douady
bundles should be constructable by a combination of the methods of sections 4.2 and 5.2. 

\medskip\noindent\textbf{Theorem 4.} \textit{Let $G_n=SU(n)$, $d|n$, $n'=n/d$,
$Z_d$ be the order-$d$ central subgroup, and $\kappa=k+n$ as before. Assume $d$ to be odd.
Then $^\kappa K^i_{G_n^{\mathrm{ad}_c}}(G_n/Z_d)$
is naturally isomorphic to $R_{Z_d}\otimes_{R_{G_n}}{}^\kappa K^i_{G_n^{\mathrm{ad}_c}}
(G_n)$ and vanishes for $i=(n+1)(n-2)/2$.  $^\kappa K^*_{G_n^{\mathrm{ad}_c}}(G_n/Z_d)$ 
will be a nimrep compatible
 with $\cZ^c_{\langle J^{n'}\rangle}$, provided Conjecture 2(a) holds for $G=G_n$.}\medskip
 
 \noindent\textit{Proof.} First note that, for any subgroup $Z'$ of $Z(G_n)$, $^{\kappa}
 K^*_{(G_n\times Z')^{\mathrm{ad}_c}}(G_n)\simeq R_{Z'}\otimes_{R_{G_n}}
 {}^\kappa K^*_{G_n^{\mathrm{ad}_c}}(G_n)$ as $R_{G_n}$-modules, using 
 \eqref{freenormal} and the automorphism of $G_n\times Z'$ sending $(g,z')\mapsto (gz',z')$.
$R_{G_n}$ multiplies $R_{Z'}$ through restriction, i.e. through $\overline{\mathrm{Res}}_{Z'}^{G_n}$.

Since $d$ is odd, we can write $\sqrt{z}:=z^{(d+1)/2}$ for any $z\in Z_d$. The key observation is 
that $(g,z).h=\sqrt{z}gh(c \sqrt{z}g)^{-1}$ since complex conjugation $c$ sends $\sqrt{z}$ to 
$1/\sqrt{z}$. Hence $^\kappa K^*_{G_n^{\mathrm{ad}_c}\times Z_d^L}(G_n)\simeq
{}^\kappa K^*_{(G_n\times Z_d)^{\mathrm{ad}_c}}(G_n)\simeq R_{Z_d}\otimes
{}^\kappa K^*_{G_n^{\mathrm{ad}_c}}(G_n)$ as desired.

The vanishing in degree $i=(n+1)(n-2)/2$ follows from Theorem 3. The nimrep statement
follows from \cite[Sect.2.1]{GaGa3}. QED\medskip

The full system, neutral system and sigma-induction  is independent of the twist $\omega$,
and so are as in section 5.2. Alpha-induction here is obtained by combining sections
4.2 and 5.2: $\alpha^+_\lambda=(1,\beta(\lambda))\in R_{Z_0}\otimes {}^\kappa
K_0^{G_n}(G_n/Z_0)$ and $\alpha^-_\lambda=(\overline{\mathrm{Res}}_{Z_0}^{G_n},\beta\circ c)$
where $\beta$ and $Z_0$ are as in Conjecture 3. The obvious guess for the assignment of
D-brane charges is ${}^\tau K_{G_n}^{n(n-1)/2}(G_n/Z_d)\rightarrow {}^\tau K^{n(n-1)/2}(G_n/Z_d)$.

There should be a similar story when both $d$ and $n'$ are even. Then $2d$ will also
divide $n$, so $Z_{2d}$ also exists. Of course $Z_{2d}$ must be cyclic; fix a generator $z$. Then
$(g,z^{2j}).h=(z^jg)hc(z^jg)^{-1}$ as before so we find $^{(0;\kappa)} 
K^*_{G_n^{\mathrm{ad}_c}}(G_n/Z_d)\simeq {}^{(0;\kappa)}
 K^*_{((G_n\times Z_{2d})/1\times Z_2)^{\mathrm{ad}_c}}(G_n)$. On the other hand, 
${}^{(0;\kappa)}K^*_{(G_n\times Z_{2d})^{\mathrm{ad}_c}}(G_n)$ will be isomorphic to
$R_{Z_{2d}}\otimes_{R_{G_n}}{}^\kappa K^*_{G_n^{\mathrm{ad}_c}}
(G_n)$, and as a $R_{G_n}$-module has a $\widehat{Z_2}$-grading obtained by restricting
the equivariance to $Z_2$, which always acts trivially. This is the $K$-theoretic explanation
for the aforementioned reducibility of $R_{Z_{2d}}\otimes_{R_{G_n}}\cN^c$. 
We'd expect $^{(0;\kappa)} K^*_{G_n^{\mathrm{ad}_c}}(G_n/Z_d)$ to be isomorphic to half
of $R_{Z_{2d}}\otimes_{R_{G_n}}{}^\kappa K^*_{G_n^{\mathrm{ad}_c}}(G_n)$, and thus to
describe $K$-theoretically the nimrep of \cite{GaGa3}.


\section{Exceptional modular invariants}

{ Section 2.3 reminds us that $G=\mathrm{SU}(2)$ has 3 exceptional
modular invariants. Two of these are due to \textit{conformal embeddings}
(defined in section 2.3); we study these in section 6.1.
Conformal embeddings for the finite groups 
 were described in \cite{EG1}. The final SU(2) modular invariant is considered in section 6.2.
Unlike the situation for the generic modular invariants, which is under excellent control
as we saw in the earlier sections, the theory for the exceptional modular invariants
is still much less satisfactory.}

\subsection{Conformal embeddings}

In \cite{EG1} we proposed that conformal embeddings $H_k\rightarrow G_1$ could be related to 
$K$-homology $^\tau K^{H^{\mathrm{ad}}}_*(G)$. We studied in detail e.g. the conformal 
embeddings $\hbox{SU}(2)_4\rightarrow\hbox{SU}(3)_1$ and
$\hbox{SU}(2)_{10}\rightarrow\hbox{Sp}(4)_1$, which give rise to the modular invariants
called $\bbD_4$ and $\bbE_6$, respectively, in the $\hbox{SU}(2)_k$
list of Cappelli \textit{et~al} \cite{CIZ}. Perhaps the
most interesting observation to come out of this analysis
was that the largest finite stabiliser in this adjoint action of
$\hbox{SU}(2)$ on $\hbox{SU}(3)$ resp. $\hbox{Sp}(4)$, is
called $\bbD_4$ resp. $\bbE_6$ on McKay's list \cite{McK}.
In this subsection we propose an alternative $K$-homology, namely 
$^\tau K_*^{\Delta_H^L\times \Delta_H^R}(G\times G)$
with the diagonal action $(h_L,h_R)(g_1,g_2)=(h_Lg_1h_R^{-1},h_Lg_2h_R^{-1})$.
The previous observation about finite stabilisers would persist in this new picture: if $K<H$ is the
stabiliser of $g\in G$, then the isomorphic copy $\{(k,k)\}<H\times H$ stabilizes $(g,z)\in G\times G$
for any $z\in Z(G)$. 

For convenience we restrict in the following to examples calculable in $K$-theory through 
the Hodgkin spectral sequence (recall section 2.2).

Consider first one of the simplest possible conformal embeddings: $T_2\subset {\rm SU}(2)_1$.
Here $G={\rm SU}(2)\times {\rm SU}(2)$ and $H={\rm SO}(2)\times {\rm SO}(2)$. 
Write $R_G=\bbZ[\sigma,\tau]$, $R_H=\bbZ[a^{\pm 1},b^{\pm 1}]$ using
obvious notation. The Verlinde ring Ver$_1(\textrm{SU2})$ is $R_G/(\sigma^2-1,
\sigma-\tau)$. The restriction map is $\sigma\mapsto a+a^{-1}$,
$\tau\mapsto b+b^{-1}$. The free resolution is
\begin{equation}
0\rightarrow R_G\stackrel{g}{\rightarrow} R_G^2\stackrel{f}{\rightarrow}R_G\rightarrow
R_G/(\sigma^2-1,\sigma-\tau)\rightarrow 0\,,
\end{equation}
where 
$f(a_1,a_2)=(\sigma^2-1)a_1+(\sigma-\tau)a_2$ and
$g(c)=(-\sigma+\tau,\sigma^2-1)c$.
The straightforward calculation shows that all $E^2_{p,q}=0$ except $E^2_{0,\textrm{odd}}=
\bbZ[a^{\pm 1},b^{\pm 1}]/(a^2+1+a^{-2},a+a^{-1}-b-b^{-1})$, which is 8-dimensional, with
coset representatives $1,a,b,a^2,ab,b^2,a^3,a^2b$. Thus $^\tau K^*_{\Delta_T^L\times 
\Delta_T^R}(G\times G)\simeq 0\oplus\bbZ^8$.

This result can also be obtained more elegantly from the `maximal rank' argument of section
3.3 of \cite{EG1}. From this we see $^\tau K_{T\times T}^0(G\times G)$ decomposes naturally
into 4 copies of the full system Ver$_2(T)$, where `4' is the order of the Weyl group of
$\textrm{SU}(2)\times\textrm{SU}(2)$, which on $R_{T\times T}$ acts by $a\mapsto a^{s_1}$, $b\mapsto
b^{s_2}$ for some choice of signs $s_1,s_2$. 

By comparison, the
$K$-group $^\tau K_{T^{\mathrm{ad}}}^*(\textrm{SU2})$  was calculated
in sections~3.2 and 3.3 of \cite{EG1}, and found to give 2 copies of the full system, where
`2' is the order of the  Weyl group of SU(2).

Thus it would appear that modelling this conformal embedding
by  the diagonal action rather than the adjoint action gains little, and in fact makes things
a little worse. However,  consider instead the $E_{8,1}\rightarrow \textrm{SU}(9)_1$
conformal embedding. Then by identical arguments $^\tau K_{\textrm{SU9}^{\mathrm{ad}}}^*(E_8)
\simeq 0\oplus \bbZ^{1920}$ 
while $^\tau K_{\Delta_{\text{SU}9}^L\times \Delta_{\textrm{SU}9}^R}^*(E_8\times E_8)\simeq 
0\oplus\bbZ^{3686400}$. The corresponding modular invariant is 
\begin{equation}\cZ=|\chi_{00000000}+\chi_{00100000}+\chi_{00000100}|^2\nonumber
\end{equation}
 so the full system is 9-dimensional. Thus $^\tau K_{\textrm{SU9}}^1(E_8)$
cannot consist of a number of copies of the full system, whereas $^\tau K_{\text{SU}9\times 
\textrm{SU}9}^1(E_8\times E_8)$ very possibly could. 

In \cite{EG1} we speculated that conformal embeddings would be realised $K$-homologically
using the adjoint action, by $^\tau K^H_*(S)$ for some $H$-invariant
submanifold $S_0$ of $G$, typically smaller than $G$. But at least sometimes the diagonal action 
could permit a
much more direct interpretation: its $K$-homology more often can consist of a certain number of 
copies of the full system.

Consider now the $\bbE_6$ modular invariant in the A-D-E list
of \cite{CIZ}: this comes from the conformal embedding 
SU$(2)_{10}\rightarrow {\rm Sp}(4)_1$. The full system is
12-dimensional, consisting of two copies of the $\bbE_6$ diagram.
In \cite{EG1} we considered
$^\tau K_*^{\textrm{SU2}^{\mathrm{ad}}}(\textrm{Sp4})$ in detail, determining it was $2+2$ 
dimensional. We would have preferred two copies of the $\bbE_6$ diagram.

But consider instead the subgroup $H={\rm SU}(2)\times{\rm SU}(2)$ of
$G={\rm Sp}(4)\times{\rm Sp}(4)$, where the embedding ${\rm SU}(2)\subset
{\rm Sp}(4)$ is the 4-dimensional irreducible 
SU(2)-representation given explicitly in \cite{EG1}.   The basic representation rings are $R_{\textrm{SU2}}=\bbZ[\sigma]$ and
$R_{\textrm{Sp4}}=\bbZ[s,v]$ where $\sigma=\sigma_2$, and $s,v$ are the spinor and
vector representations respectively. We write $R_H=\bbZ[\sigma,\sigma']$
and $R_G=\bbZ[s,v,s',v']$ using obvious notation.  We will use the
Hodgkin spectral sequence to compute $^\tau K^*_{\Delta_{\textrm{SU2}}^L\times 
\Delta_{\textrm{SU2}}^R}(\textrm{Sp4}\times\textrm{Sp4})$. $R_G$ multiplies $R_H$
using the restrictions $s\mapsto \sigma_4=\sigma^3-2\sigma, v\mapsto
\sigma_5=\sigma^4-3\sigma^2+1$. As an $R_G$-module,  we have
\begin{equation}
^\tau K^i_{\Delta_{\mathrm{Sp}4}^L\times
\Delta_{\mathrm{Sp}4}^R}(\mathrm{Sp}(4)\times\mathrm{Sp}(4))=\textrm{Ver}_1(\textrm{Sp4})=R_G/(w,x,y,z)
\end{equation}
where for later convenience we write $w=s^2-v-1,x=sv-s,y=v-v',z=s-s'$.
Note that these restrict to $(\sigma^2-2)(\sigma^4-3\sigma^2+1),
(\sigma^2-2)\sigma^3(\sigma^2-3),(\sigma^2-\sigma'{}^2)(\sigma^2+\sigma'{}^2-3),
(\sigma-\sigma')(\sigma^2+\sigma\sigma'+\sigma'{}^2-2)$ respectively.
A free resolution of Ver$_1(\textrm{Sp4})$ is
\begin{equation}
0\rightarrow R_G\stackrel{i}{\rightarrow} R_G^4\stackrel{h}{\rightarrow}
R_G^6\stackrel{g}{\rightarrow}R_G^4\stackrel{f}{\rightarrow}R_G
\rightarrow R_G/(w,x,y,z)\rightarrow 0\,,
\end{equation}
where $f(a_1,a_2,a_3,a_4)=a_1w+a_2x+a_3y+a_4z$, $i(d)=(z,-y,-x,w)d$, and 
\begin{align}
g(b_1,\ldots,b_6)=&\,(x,-w,0,0)b_1+(-y,0,w,0)b_2+(z,0,0,-w)b_3\nonumber\\ &\,+(0,-y,x,0)b_4
+(0,z,0,-x)b_5+(0,0,-z,y)b_6\,,\\
h(c_1,c_2,c_3,c_4)=&\,(y,x,0,-w,0,0)c_1+(z,0,-x,0,w,0)c_2+(0,z,y,0,0,w)c_3\nonumber\\ &
+(0,0,0,z,y,x)c_4\,.
\end{align}
By definition, $E^2_{p,0}=\mathrm{Tor}_p^{R_G}(R_H,R_G/(w,x,y,z))$ and
$E^2_{p,1}=0$ for all $p\ge 0$. We obtain $E^2_{0,0}=\bbZ[\sigma,\sigma']/(\sigma^2-2,
\sigma'{}^2-2)=:Q$, 
\begin{align}
E^2_{1,0}=&\,\mathrm{span}_\bbZ \mbox{\scriptsize$\{(\sigma^5-3\sigma^3,-\sigma^4+3\sigma^2-1,0,0),
((\sigma\sigma'-\sigma^2)(\sigma\sigma'+1),(\sigma-\sigma')(\sigma\sigma'+1),\sigma',1-\sigma'{}^2)
\}$}\,Q\,,\nonumber\\
E^2_{2,0}=&\,\mbox{\scriptsize$ (\sigma\sigma'{}^2-\sigma^2\sigma'+\sigma'-\sigma,3\sigma^3\sigma'-\sigma^5
\sigma',(3\sigma^3-\sigma^5)(\sigma'{}^2-1),\sigma^4\sigma'-3\sigma^2
\sigma'+\sigma',(\sigma^4-3\sigma^2+1)(\sigma'{}^2-1),0)$}\,Q\,,\nonumber\end{align}
and $E^2_{p,0}=0$ for $p\ge 3$. This sequence stabilises already at $E^2_{p,q}$, and
we find that both groups $^\tau K^i_{\Delta_{\mathrm{SU2}}^L\times
\Delta_{\mathrm{SU}2}^R}(\mathrm{Sp}(4)\times\mathrm{Sp}(4))$
are 8-dimensional, and as an $R_{\mathrm{SU2}\times\mathrm{SU}2}$-module
both are isomorphic to two copies of $\bbZ[\sigma,\sigma']/(\sigma^2-2,\sigma'{}^2-2)$.

In the examples of sections 3, 4 and 5, the action of the Verlinde ring on the nimrep
and full system comes from the push-forward of the inclusion of the identity of the group,
i.e. from the natural action of the representation ring on the equivariant $K$-groups.
Here, this $R_H$-action does not factor through to Ver$_{10}(\mathrm{SU2})\simeq R_H/(
\sigma^{11}-10\sigma^9+36\sigma^7-56\sigma^5+35\sigma^3-6\sigma,\sigma-\sigma')$
(the problem is the $\sigma-\sigma'$). This module structure should come from
the product discussed in section 2.2.

In \cite{EG1} we modeled the conformal embedding
$\hbox{SU}(2)_4\rightarrow\hbox{SU}(3)_1$ by ${}^\tau K_{\textrm{SU2}^{\mathrm{ad}}}^*(\textrm{SU3})$
and got $1+1$ dimensions (the full system is 8-dimensional, consisting of
two copies of the $\bbD_4$ graph). But we now see that a better approximation
should be ${}^\tau K_{\Delta_{\textrm{SU2}}^L\times 
\Delta_{\textrm{SU2}}^R}^*(\textrm{SU3}\times \textrm{SU3})$. It would take some 
work to compute this, since the Hodgkin spectral sequence doesn't apply here
(this embedding is of SO(3), not SU(2)). Let us make some qualitative
comments. The dimension of these $K$-groups should surely be greater
than $1+1$, as the above examples  indicate. To what
extent can we hope to see $\bbD_4$'s in these $K$-homological
groups? The distinguishing feature of the
(unextended) $\bbD_4$ diagram is the $S_3$ symmetry of the
three endpoints, which fixes the central vertex. This $S_3$
symmetry appears naturally here: on each factor space
SU(3) it is generated by multiplying by a
scalar matrix $\omega^iI$ (these form the centre of ${\rm
SU}(3)$), and by complex conjugation --- all of these
commute with the ${\rm SU}(2)$ action. So it is not impossible
to imagine $\bbD_4$ lurking here.

In hindsight it is clear why we are not getting the correct answer spot-on here:
we are ignoring the spinors. More precisely, the dimension shift dim$(G)$
in the degree, and the shift by $h^\vee$ of the level, in the identification
Ver$_k(G)\simeq {}^{k+h^\vee}K^{\mathrm{dim}\,G}_{G^{\mathrm{ad}}}(G)$ 
(for $G$ simply-connected, connected and compact) come from
the implicit presence of a Cliff$(L\mathfrak{g}^*)$-module, which comes along
for the ride. However it appears it cannot be factored off so simply in
the conformal embedding context: in particular the level shift does not
occur, as was demonstrated in the level calculations in \cite[Sect.2.3]{EG1}. Considering
for concreteness the SU$(2)_4\rightarrow\mathrm{SU}(3)_1$ conformal embedding,
 we know from \cite[Sect.2.3]{EG1} that the twist on the SU(2) and SU(3) parts really 
should be 4 and 1 respectively, and not $4+2$ and $1+3$. So this means the
spinors cannot be so easily ignored. Including them will account for at
least some of the discrepancy found above.

\subsection{The $\bbE_7$ modular invariant of SU(2)}

The preceding sections address all  modular invariants of  $G=\mathrm{SU}(2)$, except for
$\bbE_7$ at level 16. Its direct interpretation would be as a twist of the even part of
$\bbD_{10}$; in this sense it is exceptional in that no other $\bbD_n$ have such a twist.
 But in another sense the $\bbE_7$ modular invariant is not really exceptional:
it belongs to an infinite sequence at SO($n$) level 8 coming from a 
combination of \textit{level-rank duality} (which relates SO$(n)_k$ 
with SO$(k)_n$, at least when $kn$ is even) and SO(8) triality. In this subsection we
obtain a partial $K$-theoretic realisation of  $\bbE_7$ by approaching it in this way. 

The neutral system (maximal chiral extension) for $\bbE_7$
is $\bbD_{10}$, i.e. $^{8+1}K_{\mathrm{SO3}}^1(\mathrm{SO(3)})$; the branching
rules (sigma-restriction) $^{8+1}K_{\mathrm{SO3}}^1(\mathrm{SO(3)})\rightarrow
{}^{16+2}K_{\mathrm{SU2}}^1(\mathrm{SU(2)})$ is $\overline{0}\leftrightarrow
0+16$, $\overline{1}\leftrightarrow 2+14$, $\overline{2}\leftrightarrow
4+12$, $\overline{3}\leftrightarrow 6+10$, and the resolved fixed-point
$\overline{4},\overline{4}'\leftrightarrow 8$. Level-rank duality comes from
the conformal embedding SO$(24)_1\leftarrow \mathrm{SO}(8)_3\oplus\mathrm{SO}(3)_8$,
which has branching rules \cite{Vrst}
\begin{align}
0&\,\mbox{\scriptsize$\mapsto(0000)(0+16)+(2000)(2+14)+(0100)(4+12)+(1011)(6+10)+((0020)+(0002))8$}\\
v&\,\mbox{\scriptsize$\mapsto(3000)(0+16)+(1000)(2+14)+(1100)(4+12)+(0011)(6+10)+((1020)+(1002))8$}\\
s&\,\mbox{\scriptsize$\mapsto(0003)(0+16)+(0021)(2+14)+(0101)(4+12)+(1010)(6+10)+((2001)+(0001))8$}\\
c&\,\mbox{\scriptsize$\mapsto(0030)(0+16)+(0012)(2+14)+(0110)(4+12)+(1001)(6+10)+((2010)+(0010))8$}
\end{align}
using obvious notation. Thus we get the identification 
$[(0000)]\leftrightarrow\overline{0}$, $[(2000)]\leftrightarrow\overline{1}$, 
$[(0100)]\leftrightarrow\overline{2}$, $[(1011)]\leftrightarrow\overline{3}$, 
$[(0020)]\leftrightarrow \overline{4}$, $[(0002)]\leftrightarrow\overline{4}'$
between the $K$-groups $^{3+9}K^0_{\mathrm{PSpin8}}(\mathrm{PSpin(8)})$ and
$^{8+1}K^1_{\mathrm{SO}3}(\mathrm{SO}(3))$, equivalently between the
corresponding orbits of simple-currents.
The $\bbE_7$ modular invariant is recovered from the diagonal modular
invariant of SO(24)$_1$ and the triality modular invariant of SO(8)$_3$,
by \textit{contraction} \cite{DV}.

The full system of $\bbE_7$ was obtained in \cite{O,BEK2}; we reproduce
the figure below:

\medskip\epsfysize=1.75in\centerline{ \epsffile{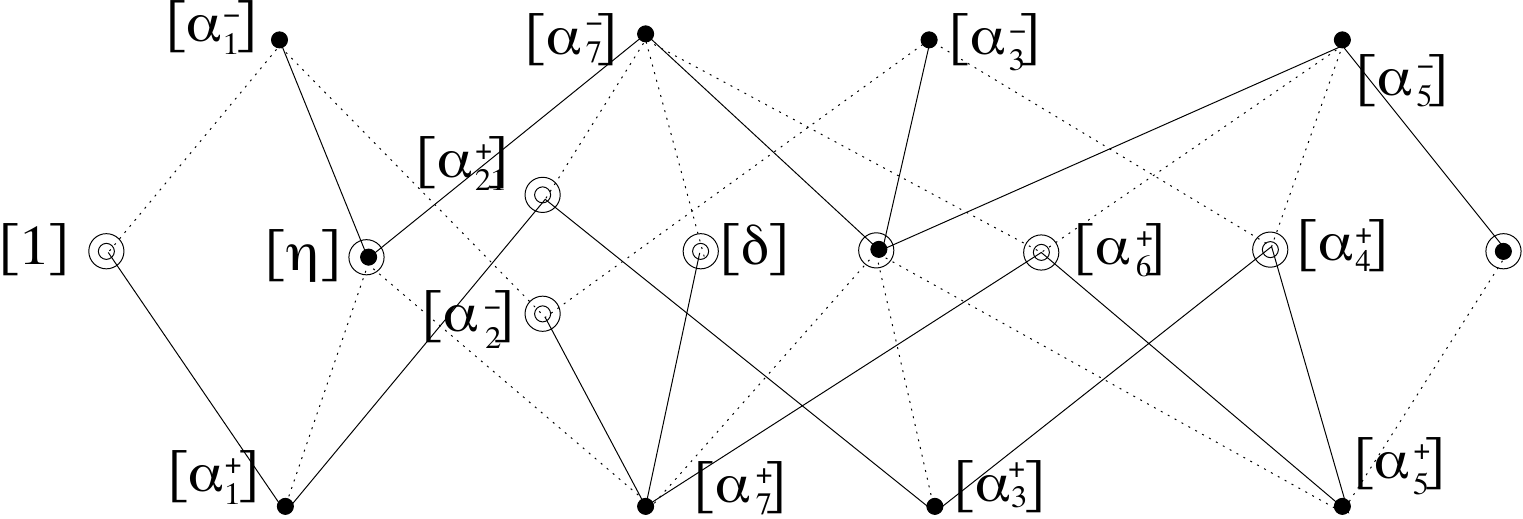}}\medskip

\centerline{Figure 3. Full system for $\bbE_7$}
\medskip

Focus on the solid lines for now (the graphs describing multiplication by
$\alpha^+_1$).
The copy of $\bbD_{10}$ is clear, as it represents the nimrep of the
neutral system.
$K$-theoretically, of course it is $^{16+2}K^1_{\mathrm{SU2}}(\mathrm{SO3})$.
Note that the 6 even vertices of $\bbD_{10}$ form the Verlinde ring 
$^{(1;0,8+1)}K^1_{\mathrm{SO3}}(\mathrm{SO3})$ of the neutral system.
The odd vertices are then the 4 spinors $^{(1;1,8+1)}K^1_{\mathrm{SO3}}(\mathrm{SO3})$.
The 3 even vertices of the $\bbE_7$ graph can be identified with the
image under level-rank duality of the nimrep for the triality modular
invariant of SO(8)$_3$, namely the $K$-group $^{3+6}K^0_{\mathrm{Spin8}}(\mathrm{Spin}(8))$
where the adjoint action is twisted by triality. 

The remaining question is how to recover the 4 remaining odd vertices of $\bbE_7$.
It is tempting to guess that these are 4 spinors for the twisted adjoint 
action of Spin(8) on itself, but apparently there are no spinors at shifted level
3+6. In some sense we already know these 4 odd vertices: they are the
odd vertices of the dotted copy of $\bbD_{10}$. In any case, most of
the $\bbE_7$ full system is clear.

A proper interpretation of the $\bbE_7$ full system, as a single $K$-group,
would presumably involve SO(24) at level 1, but although there are many
candidates, it isn't clear yet to us which is the correct one. It would
certainly help to understand $K$-theoretically level-rank duality --- this should be possible
along the lines of \cite{xu2}, but
may require understanding cosets or conformal embeddings from this framework.

\section{Outlook and speculations}

This paper is the second of a series devoted to deepening
the connection between twisted equivariant $K$-homology and
conformal field theory.  This concluding section entertains some further possibilities.

Deep connections between $K$-theory and conformal field theory$/$string theory
have been known for some time. But this paper establishes fundamental and systematic
roles $K$-theory plays in a much more extensive range of CFT structures than had been
previously appreciated. There is much more to do though.  

One of the simplest examples of an orbifold is the \textit{permutation orbifold}: let 
$H$ be any subgroup of the symmetric group $S_n$, and let $\cV$ be any VOA
(or RCFT); the permutation orbifold is $\cV^{\otimes_n}/H$, where $H$ acts by
permuting the copies of $\cV$. A basic fact (see e.g. \cite{Ban}) is  that successive
permutation orbifolds, first  by $G<S_m$ and then by $H<S_n$, is equivalent 
to a single one of $\cV$ by the wreath product $G\wr H<S_{mn}$. Now, the modular
data of the Drinfeld double of finite group $G$ can be regarded as that of
the permutation orbifold of a holomorphic VOA $\cV$ by $G$. Thus we obtain the
observation that the Verlinde ring of the permutation orbifold (by any subgroup $H$ of $S_n$) of
the  Drinfeld double of $G$,  has the expression $K^0_{G\wr
H}(G\wr H)$. This simple example is probably worth studying in more detail:
see e.g.  \cite{Wng,Gnt} for some of the rich structure present.
However the $K$-theoretic treatment of permutation orbifolds of loop group
data is still not clear to us (see the examples given in \cite{EG1}).

It would be very interesting to interpret $K$-theoretically the 
Goddard-Kent-Olive coset construction \cite{GKO}. Of course
\textit{coset} here is not meant to be taken literally --- e.g. it does not
refer physically to a string living on the homogeneous space $G/H$. 
Algebraically, 
it involves commutants: e.g.  the {coset} of a vertex operator algebra $\cV$ by a subalgebra
$\cW$ is the commutant of $\cW$ in $\cV$;  the coset of subfactor $N\subset M$ by
subfactor $S\subset N$ is the subfactor $S'\cap N\subset S'\cap M$.
A promising $K$-theoretic approach is based on \cite{GaWa}: understand 
the coset model $G/H$ as an orbifold of $G \times { H}$ by the intersection
$Z(G) \cap Z(H)$ of centres. The main difficulty seems to be the orbifold part.
This approach won't
always work: e.g. maverick cosets \cite{DJ} such as conformal embeddings
have \textit{identifications} of primaries not merely given by the simple-currents
$Z(G)\cap Z(H)$. In any case the study of coset models provides further 
motivation for developing the theory of orbifolds.

Note that the group of D-brane charges for the modular invariants of SU(2) is 
given by the centre of the corresponding A-D-E Lie group (recall the discussion
at the end of section 2.3). This surprising fact has a simple $K$-theoretic
interpretation using the trivial action of the maximal torus of each of those groups on itself.
The boundary states then are given by the simple roots, and the D-brane charges can then
be recovered as the inner products with appropriate weights. It would be interesting to
understand this elegant (though rather mysterious) 
description of the D-brane charges for SU(2), from our framework. 

As mentioned in section 2.3, it isn't obvious how to interpret \eqref{charges}
when $G$ is finite. Suppose we extrapolate the $K$-theoretic treatment
of the Verlinde D-brane charges for Lie groups $G$ as given in section 2.4, 
to the finite group setting. In both settings the Verlinde ring is given
by ${}^\tau K^{G^{\mathrm{ad}}}_0(G)$; this suggests postulating that the 
charges for the finite group should also be given by the forgetful map
$K^0_G(G)\rightarrow K^0(G)=\bbZ(G)$. This map sends the primary $(g,\chi)$
to dim$(\chi)$ times the conjugacy class of $g$. Thus the image of the
forgetful map is a free $\bbZ$-module of dimension equal to the class number of $G$. But by the
argument given at the end of section 2.3, for any Verlinde nimrep all
charges are uniquely determined from the charge of the vacuum, and thus the
charge-group of the Verlinde nimrep will always be cyclic. This means that
the assignment of charges here can be given by the forgetful map only when
$G$ has class number 1, i.e. only for the trivial case $G=1$.

So the analogy between finite groups and Lie groups is not perfect in the
context of D-brane charges. Nevertheless, one may hope that the Verlinde
charge-group
$\cM_N$ for $G$ finite is a (cyclic) subgroup of $K^0(G)=\bbZ(G)$. This 
would imply that whatever we take `dim$(g,\chi)$' to be, \eqref{charges}
will be satisfied \textit{exactly} (i.e. with $M=\infty$). There are
dim$(\mathrm{Ver}(G))$ independent ways to do this, given by the assignments 
$q_{(g,\chi)}=S_{(g,\chi),(g',\chi')}/S_{(1,1),(g',\chi')}$ for each fixed
primary $(g',\chi')$; for only one choice of $(g',\chi')$, namely $(1,1)$,
 will these `charges' be positive integers. Therefore it is tempting to
suggest that \textit{the correct choice of `$\mathrm{dim}(g,\chi)$' in 
\eqref{charges} for $G$ a finite group, is the \textit{quantum-dimension} 
$\|K_g\|\,\mathrm{dim}(\chi)$,} where $K_g$ denotes the conjugacy class of $g$.
It would be interesting to compare these quantum-dimensions with the
dimensions of Nahm's special spaces, for the fixed-point subVOAs $\cV^G$ of 
holomorphic VOAs $\cV$.

Langlands duality relates the groups ${\rm SU}(n)/\bbZ_{d}$ and ${\rm SU}
(n)/\bbZ_{n/d}$ for example. Could this be manifested  perhaps through a T-duality 
between  the corresponding modular invariants,
e.g. between the $\bbD$-series and $\bbA$-series of SU(2)?

The Verlinde ring and nimrep both come with preferred bases, in which the
structure constants are nonnegative integers. This appears to be more evident in the $K$-theory
language than in $K$-homology. To understand this point, consider the finite group case.
In either case the double cosets arise automatically from the orbit
analysis, but the appearance of representation rings of stabilisers
would require Poincar\'e duality if we use $K$-homology. For $K$-homology,
what would arise naturally presumably would be conjugacy classes of the
stabilisers --- this is the \textit{permutation basis} of the Verlinde ring,
in which SL(2,$\bbZ$) has a monomial representation. Maybe this is suggestive: 
from $K^*$ here we get the preferred basis of
Verlinde ring by primaries, but from $K_*$ we get the permutation
basis in which the modular group action is cleanest. The modular
group representation is much cloudier in the primary$=K^*$ basis, and the fusion
structure is much cloudier in the permutation$=K_*$ basis. Incidentally, there is no
analogue of the permutation basis in general for the loop group case. 


\newcommand\biba[7]   {\bibitem{#1} {#2:} {\sl #3.} {\rm #4} {\bf #5,}
                    {#6 } {#7}}
                    \newcommand\bibx[4]   {\bibitem{#1} {#2:} {\sl #3} {\rm #4}}

\def\ASENS            {Ann. Sci. \'Ec. Norm. Sup.}
\def\AM   {Acta Math.}
   \def\AnM              {Ann. Math.}
   \def\CMP              {Commun.\ Math.\ Phys.}
   \def\IJM              {Internat.\ J. Math.}
   \def\JAMS             {J. Amer. Math. Soc.}
\def\JFA              {J.\ Funct.\ Anal.}
\def\JMP              {J.\ Math.\ Phys.}
\def\JRA              {J. Reine Angew. Math.}
\def\JSP              {J.\ Stat.\ Physics}
\def\LMP              {Lett.\ Math.\ Phys.}
\def\RMP              {Rev.\ Math.\ Phys.}
\def\RNM              {Res.\ Notes\ Math.}
\def\RIMS             {Publ.\ RIMS.\ Kyoto Univ.}
\def\Inv              {Invent.\ Math.}
\def\npbp             {Nucl.\ Phys.\ {\bf B} (Proc.\ Suppl.)}
\def\nupb             {Nucl.\ Phys.\ {\bf B}}
\def\nup              {Nucl.\ Phys. }
\def\nupp             {Nucl.\ Phys.\ (Proc.\ Suppl.) }
\def\adma             {Adv.\ Math.}
\def\coma             {Con\-temp.\ Math.}
\def\PAMS             {Proc. Amer. Math. Soc.}
\def\PJM              {Pacific J. Math.}
\def\ijmp             {Int.\ J.\ Mod.\ Phys.\ {\bf A}}
\def\jpa              {J.\ Phys.\ {\bf A}}
\def\PLB              {Phys.\ Lett.\ {\bf B}}
\def\RIMS             {Publ.\ RIMS, Kyoto Univ.}
\def\Top               {Topology}
\def\TAMS             {Trans.\ Amer.\ Math.\ Soc.}

\def\Duke              {Duke Math.\ J.}
\def\K                 {K-theory}
\def\JOP               {J.\ Oper.\ Theory}

\vspace{0.2cm}\addtolength{\baselineskip}{-2pt}
\begin{footnotesize}
\noindent{\it Acknowledgement.}

The authors thank the University of Alberta Mathematics Dept, Banff International Research
Station, Cardiff School of Mathematics, Dublin Institute for Advanced Study, Erwin-Schr\"odinger
Institute (Vienna), Swansea University Dept of Computer Science, and Universit\"at  
W\"urzburg Institut f\"ur Mathematik for generous hospitality while researching this
paper. Their research was supported in part by  EPSRC
GR/S80592/01,  EU-NCG Research Training Network: MRTN-CT-2006 031962,
DAAD (Prodi Chair), and NSERC. We thank Matthias Gaberdiel, Nigel Higson and Constantin Teleman
for conversations.

\end{footnotesize}

\end{document}